\documentclass[a4paper,11pt]{article}
\usepackage{amsmath,amssymb,amsthm,MnSymbol}
\usepackage{paralist,enumitem,hyperref,url,bm,placeins}
\usepackage{color,graphicx}
\usepackage[margin=30mm]{geometry}

\newcommand{\N}{\mathbb{N}}
\newcommand{\R}{\mathbb{R}}
\newcommand{\bu}{\bm{u}}

\newcommand{\pd}{\partial}
\renewcommand{\div}{\, \mathrm{div} \,}

\newcommand{\x}{\mathbf{x}}

\theoremstyle{plain}
\newtheorem{thm}{Theorem}[section]

\newtheorem{remark}{Remark}[section]
\newtheorem{ex}{Example}[section]
\newtheorem{defn}{Definition}[section]

\numberwithin{equation}{section}

\title{Well-posedness and finite-time extinction of a PDE-ODE spatial-network model with anisotropic diffusion}
\author{Xiao Meng \footnotemark[1] \and Kei Fong Lam \footnotemark[1]} 
\date{ }

\begin{document}

\maketitle

\renewcommand{\thefootnote}{\fnsymbol{footnote}}
\footnotetext[1]{Department of Mathematics, Hong Kong Baptist University, Kowloon Tong, Hong Kong \tt(\{22482202,akflam\}@hkbu.edu.hk)}

\begin{abstract}
We study a system of reaction-diffusion equations posed on a bounded domain composed of subdomains separated by a connected network with a metric graph structure. The reaction-diffusion dynamics with anisotropic diffusion on the graph edges are coupled to well-mixed ODE dynamics occurring at the vertices by junction conditions, and to similar PDE dynamics occurring on adjacent subdomains through Robin-like boundary conditions. The resulting PDE-ODE system can be used in epidemiological and ecological settings to study population movement in between cluster centers along road-like structures and into the surrounding continuum. We employ a semi-Galerkin approximation to establish the well-posedness of weak solutions to the PDE-ODE system, and examine further properties such as regularity, boundedness and finite-time extinction.
\end{abstract}



\section{Introduction}\label{sec:intro}
Temporal evolution and spatial diffusion phenomena are prevalent in real-world applications and theoretical fields such as engineering, mathematics and epidemiology. The dynamics of diffusion within a continuous domain, such as the spreading of heat or chemicals along pipelines, the diffusion of pollutants through river networks, the movement of individuals along transportation systems, and the transport of goods by trucks on highways, are typically described by partial differential equations (PDEs). In contrast, temporal dynamics that occur at individual nodes or compartments, such as population growth, chemical reactions at junctions, or the evolution of infectious diseases within towns, are often modeled using systems of ordinary differential equations (ODEs). Both approaches have been utilised in a wide range of applications. For example, in the context of epidemiology, understanding the spatial spread of pathogens via transportation routes requires combining these models with local transmission dynamics in population centers to provide a broader and more rigorous analysis. This necessity is highlighted by numerous historical and modern examples, such as the spread of cholera through water sources in 19th-century London \cite{Nwabor} and along the fluvial system \cite{Righetto}, the transmission of the Black Death along medieval trade routes \cite{Schmid}, diffusion of dengue along the road network in an urban area \cite{Li}, and the propagation of modern pandemics through interconnected networks of roads and airways \cite{Roques}. These cases illustrate the evolution of mathematical models from specialized modules focusing on either continuous or discrete domains toward unified descriptions that integrate both spatial diffusion and local node dynamics.

To capture the structure of such coupled multiscale problems, metric graphs have emerged as a powerful and flexible mathematical framework \cite{Berkolaiko,Besse,Klopp,Roques}. A metric graph consists of a finite set of vertices connected by edges, where each edge is endowed with a metric structure that enables an identification with an open interval of the real line. This allows the edge to be parameterized and thus serve as the domain of suitable PDEs. ODEs can be posed on the vertices of the metric graph to enable interactions between edge-based and vertex-based dynamics via suitable boundary conditions \cite{Berkolaiko1,Berkolaiko,Brio}. This mixed 1D-0D structure naturally facilitates modeling processes in which spatial diffusion interacts with local responses, as observed in epidemiology, ecology, and various engineering applications. The theory of coupled PDE-ODEs on metric graphs has developed rapidly in recent decades, with advances in spectral theory \cite{Berkolaiko1,Berkolaiko}, localization phenomena \cite{Klopp}, and computational methods for eigenvalue problems \cite{Brio}. Applications span quantum chaos \cite{Klopp}, stochastic network lasers \cite{Kim,Work}, power grids \cite{Frolov}, water grids \cite{Martin}, and epidemiology \cite{Berestycki,Besse} to name a few. The latter demonstrates that metric graphs effectively combine rapid diffusion along transportation networks with population dynamics at vertices, providing a comprehensive context for studying complex epidemic dynamics.

On the other hand, metric graphs have also been used within the context of population dynamics as narrow linear elements embedded with a two-dimensional domain (termed ``matrix'') separating two or more two-dimensional subdomains (termed ``patches'') \cite{Berestycki,BerestyckiKPP,Kravitz,Roques}. The dynamics of species in such types of environment can be significantly influenced in the presences of such narrow elements, which can that facilitate movement between adjacent patches of habitat \cite{Gilbert,Haddad}, or act as barriers to hinder species movement \cite{Garcia,Klaus}. These mixed 2D-1D models couple reaction-diffusion equations on the narrow linear elements represented as a metric graph with reaction-diffusion equations on adjacent two-dimensional patches via suitable boundary conditions and exchange mechanisms.

In this work, we study a simpler variant of a model proposed in \cite{Kravitz} that combines the aforementioned elements, that is, a mixed 2D-1D-0D diffusive SIR (Susceptible-Infected-Removed) system motivated from epidemiology in which a bounded domain $\Omega$ is subdivided into multiple subdomains/patches $\{\Omega_i\}_{i=1}^M$ by a metric graph $E$, and on each patch $\Omega_i$ we couple a PDE to equations on adjacent edges $e_j \sim \Omega_i$ that make up the boundary of the patch, as well as ODEs on the corresponding vertices $v_k \sim e_j$, see e.g.~Figure \ref{fig:domain} for an illustration of the geometric setting and Section \ref{sec:model} for the $\sim$ notation. Let us state the main features of our model:
\begin{itemize}
\item In each subdomain $\Omega_i$ and edge $e_j$, we pose a second-order reaction-diffusion equation with anisotropic spatial diffusion \cite{Hillen,McKenzie} and possibly nonlinear reaction terms.
\item Material can transfer between subdomain and adjacent edges via Robin type boundary conditions. The rates of transfer from subdomain to edge and from edge to subdomain need not be the same.
\item In each vertex $v_k$, we pose a linear ordinary differential equation, and material can transfer between all edges sharing the same vertex with possibly different transfer rates by suitable junction conditions.
\end{itemize}

Our main contributions are theoretical results on the well-posedness, regularity and finite-time extinction of weak solutions to such class of mixed dimensional models. The key idea is to employ a semi-Galerkin scheme, similar in spirit to \cite{Besse}, expressing the vertex ODE solutions as functions of the edge PDE solutions, and combined with a fixed point argument that subsequently yields the existence and uniqueness of weak solutions in the most general setting. Under settings where some of the transfer coefficients are symmetric, we can establish further solution regularity and even show finite-time extinction.

The organization of the paper is as follows. In Section \ref{sec:model} we introduces the mathematical formulation for the mixed dimensional model. Section \ref{sec:main} displays the main assumptions and results. The proof of well-posedness is contained in Section \ref{sec:wellposed}, while the regularity of solutions is shown in Section \ref{sec:reg}. In Section \ref{sec:finitetime} we explore the finite-time extinction of solutions. A formal derivation of the model equations can be found in Appendix \ref{sec:derivation}, and useful regularity results for second order quasilinear elliptic systems are contained in Appendix \ref{sec:Garcke}.

\section{Model equations}\label{sec:model}
Let $\Omega \subset \R^2$ be a bounded domain with polygonal boundary $\Gamma:= \pd \Omega$, and let $\{\Omega_i\}_{i =1}^M$ denote a disjoint partition of $\Omega$ with open, bounded and polygonal subsets $\Omega_i \subset \R^2$. For a pair $i < j$, the intersection between $\pd \Omega_i$ and $\pd \Omega_j$ is a line segment. Then, the boundary of the subregion $\Omega_i$ consists of two portions: the collection of all such line segments with neighbouring subregions, and (if any) portions of the external boundary $\pd \Omega$. The collection of all such line segments forms a graph $E$ inside $\Omega$ that separate the subregions from each other, along with vertices at the location where two or more of these line segments intersect. We exclude situations where there are vertices of degree 1, so that the graph $E$ is without boundary.

\begin{figure}[h]
\centering
\includegraphics[width=0.6\textwidth]{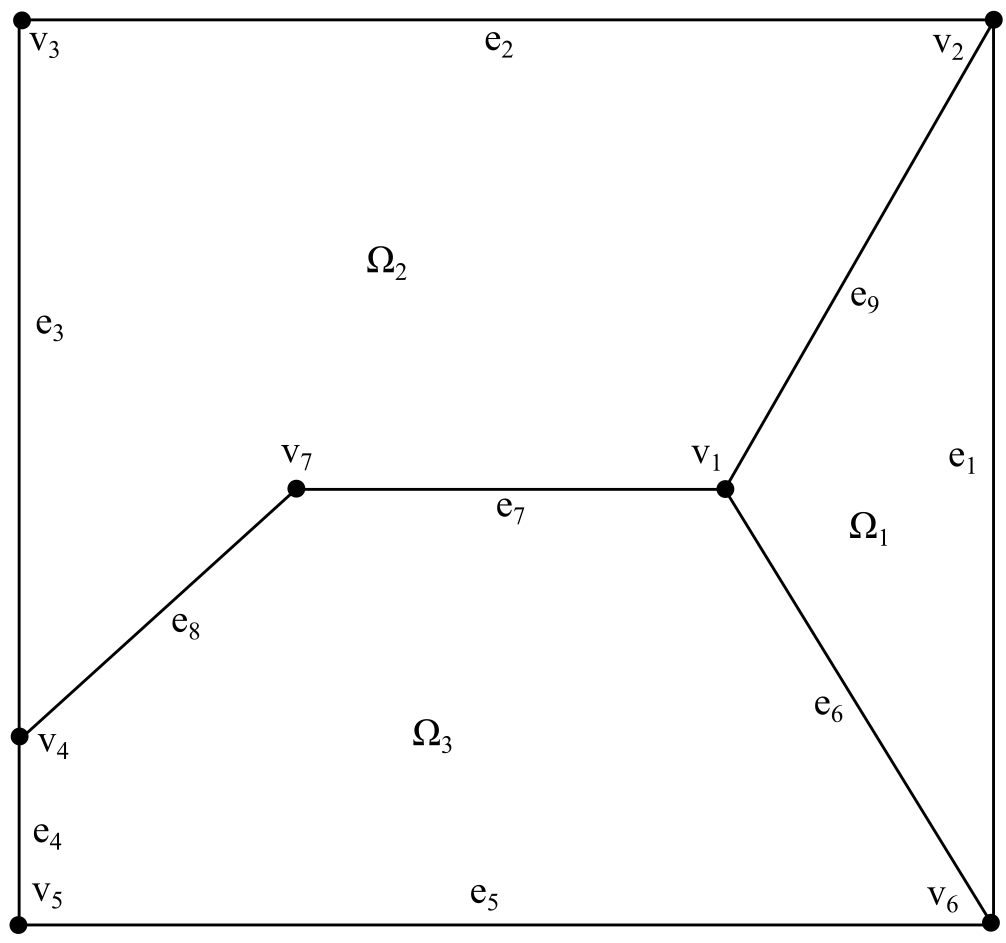}
\caption{The domain $\Omega$ consists of three subdomains, separated by nine edges that meet at seven vertices. The external boundary is comprised of edges $e_1$, $e_2$, $e_3$, $e_4$ and $e_5$. In this example, the edges adjacent to $\Omega_1$ are $e_1$, $e_6$ and $e_9$, while the edges that share vertex $v_1$ as one of their endpoints are $e_6$, $e_7$ and $e_9$.}
\label{fig:domain}
\end{figure}

We now endow the collection of line segments with a metric graph structure, and hence refer to them as \emph{edges} of the metric graph. Let $E$ be the collection of these edges and $V$ be the collection of vertices, where we use $|E|$ and $|V|$ to denote their cardinalities. Each edge $e_j$, $j = 1, \dots, |E|$, has a positive length $l_j> 0$ and an arbitrary orientation with a one-dimensional coordinate system $x \in [0,l_j]$ where we identify $x = 0$ as the source vertex of the edge $e_j$ and $x = l_j$ as the terminal vertex of the edge $e_j$. We refer to Figure~\ref{fig:domain} for an example of the domain geometry described above.

Let $T> 0$ be an arbitrary but fixed terminal time. We use the notation $\Omega_i \sim e_j$ to denote that all the subdomains $\Omega_i$ that contains edge $e_j$ as part of their boundaries, or all edges $e_j$ that are adjacent to the subdomain $\Omega_i$. Likewise we use the notation $e_j \sim v_k$ to denote all edges $e_j$ that share the vertex $v_k$ as one of their endpoints, or all vertices $v_k$ of edge $e_j$. The model that we study in this work reads as follows: 
\noindent In the subdomain $\Omega_i$ with adjacent edges $e_j$ we posed the following reaction-diffusion equation
\begin{subequations}\label{dom:equ}
\begin{alignat}{2}
\label{u:domain} \pd_t u_i - \div (\kappa(\nabla u_i)) + f(u_i) = 0 & \quad \text{ in } \Omega_i \times (0,T), \\
\label{u:bdy} \kappa(\nabla u_i) \cdot \bm{n} = \alpha_{ij} w_j - \beta_{ij} u_i & \quad \text{ on } e_j \times (0,T), 
\end{alignat}
\end{subequations}
where $\bm{n}$ is the outward point unit normal, $f$ is a reaction term, $\kappa : \R^2 \to \R^2$ is an anisotropic diffusive flux with antiderivative $\widehat{\kappa} : \R^2 \to \R$ such that $\nabla \widehat{\kappa} = \kappa$. Examples we have in mind are
\[
\widehat{\kappa}(\bm{x}) = \frac{1}{p}|\bm{x}|^p \text{ or } \quad \widehat{\kappa}(\bm{x}) = \frac{1}{2}|\bm{x}|^2 + \frac{1}{p} |\bm{x}|^p, \quad p \in [2,\infty)
\]
so that 
\[
\kappa(\nabla u) = |\nabla u|^{p-2} \nabla u \quad \text{ or } \quad \kappa(\nabla u) = (1 + |\nabla u|^{p-2}) \nabla u.
\]
The equation governing $u_i$ in $\Omega_i$ is a quasilinear reaction-diffusion equation with Robin-type boundary conditions to adjacent edges $e_j$. 

On the edge $e_j$ with adjacent subdomains $\Omega_m$ we pose the following reaction-diffusion equation
\begin{align}\label{w:edge}
 \pd_t w_j - \pd_x (\eta(\pd_x w_j)) + g(w_j) + \sum_{m \, : \, \Omega_m \sim e_j}^M (\alpha_{mj} w_j - \beta_{mj} u_m) = 0 & \quad \text{ on } e_j \times (0,T).
\end{align}
This is a similar quasilinear reaction-diffusion equation to \eqref{u:domain} but with a reaction term $g$ and an anisotropic diffusive flux $\eta: \R \to \R$ with antiderivative $\widehat{\eta} :\R \to \R$ such that $\widehat{\eta}' = \eta$, as well as additional source terms accounting for the exchange of material with all adjacent subdomains $\Omega_m \sim e_j$. For edges that are parts of the external boundary $\pd \Omega$, the summation term reduces to a single contribution. To prescribe a junction/boundary condition to \eqref{w:edge} at a (source/terminal) vertex $v_k$, we consider the balance of fluxes across the vertex $v_k$ under two scenarios:
\begin{itemize}
\item The first scenario is the case where the vertex $v_k$ does not hold a concentration of material, which we term as \emph{unpopulated}. Let $D_{k}$ denote the degree of vertex $v_k$, i.e., the total number of edges that have $v_k$ as an endpoint, and set $\bm{w}_{v_k}$ as the vector of length $D_k$ with elements $w_{j}$ such that edge $e_j \sim v_k$, as well as $\bm{\eta}(\pd_x \bm{w}_{v_k})$ as the vector of length $D_k$ with elements $\eta(\pd_x w_j)$ representing the outward flux. Then, the balance of fluxes at $v_k$ can be described by the following system:
\begin{align}\label{unpop:bc}
\bm{\eta}(\pd_x \bm{w}_{v_k}) + \bm{N}_{v_k} \bm{w}_{v_k} = \bm{0},
\end{align}
where $\bm{N}_{v_k}$ is a $D_k \times D_k$ matrix such that for $n \neq m$, the entry $(\bm{N}_{v_k})_{n,m} = -\gamma^k_{m \to n}$ encodes the constant rate of material transferring from edge $e_m$ to $e_n$ via the vertex $v_k$, while for the diagonal terms $n = m$, the entry $(\bm{N}_{v_k})_{n,n} = \sum_{m \neq n} \gamma^k_{n \to m}$ captures the rate of material transferring from edge $e_n$ to other edges connected at $v_k$. Hence, the column and row sums of $\bm{N}_{v_k}$ is always zero. The system \eqref{unpop:bc} can be viewed as a Robin type boundary condition for the edge equations \eqref{w:edge} at vertex $v_k$ that permits the transfer of material across all other edges connected at $v_k$.
\item The second scenario is the case where the vertex $v_k$ holds a concentration of material denoted as $z_k$, which we term as \emph{populated}. We thus consider an ordinary differential equation for $z_k$ of the form
\begin{align}\label{pop:equ}
\frac{d}{dt} z_k = \sum_{j \, : \, e_j \sim v_k}^{|E|} (\delta^k_{e_j} w_j(t, v_k) - \lambda^k_{e_j} z_k),
\end{align}
where $\lambda^k_{e_j}, \delta^k_{e_j} \geq 0$ are constants representing the rate of material transfer from vertex $v_k$ into edge $e_j$, and from edge $e_j$ into vertex $v_k$, respectively. We now modify \eqref{unpop:bc} to take into account source terms from the vertex $v_k$ to the edges:
\begin{align}\label{pop:bc}
\bm{\eta}(\pd_x \bm{w}_{v_k}) + (\bm{N}_{v_k} + \bm{E}_{v_k}) \bm{w}_{v_k} =  z_k \bm{\lambda}_{v_k},
\end{align}
where the diagonal matrix $\bm{E}_{v_k}$ is of size $D_k \times D_k$ with entries $(\bm{E}_{v_k})_{n,n} = \delta^k_{e_n}$ and $\bm{\lambda}_{v_k}$ denotes a vector of size $D_k$ with entries consisting of $\lambda^k_{e_j}$. 
\end{itemize}
Note that in the special case $\delta^k_{e_j} = \lambda^k_{e_j} = 0$, the junction condition \eqref{pop:bc} reduces to \eqref{unpop:bc}. Thus, if there is an unpopulated vertex $v_K$ we can set $\delta^K_{e_j} = 0$ and $\lambda^K_{e_j} = 0$ for $e_j \sim v_K$ and initialize with $z_K(t=0) = 0$ in \eqref{pop:equ} to obtain $z_K(t) = 0$ for all $t > 0$.

We furnish the above system with initial conditions:
\begin{align}\label{ini}
u_i(0,x) = u_{i,0}(x) \text{ in } \Omega_i, \quad w_j(0,x) = w_{j,0}(x) \text{ on } e_j, \quad z_k(0) = z_{k,0} \text{ for } v_k,
\end{align}
for $1 \leq i \leq M$, $1 \leq j \leq |E|$ and $1 \leq k \leq |V|$.

\begin{remark}\label{rem:Kirchhoff}
If we consider the special case where $\gamma^k_{m \to n} = \gamma^k_{n \to m}$ in \eqref{unpop:bc}, by summing the components of the system \eqref{unpop:bc} we obtain 
\begin{align}\label{kirchhoff}
\sum_{j\, : \, e_j \sim v_k}^{|E|} \eta(\pd_x w_{j}) = 0,
\end{align}
i.e., the sum of the outward fluxes of each edge associated to vertex $v_k$ is zero. Pairing this with the additional continuity condition that $w_j = w_J$ at $v_k$ for $e_j, e_J \sim v_k$, this leads to the more commonly used Kirchhoff boundary conditions on metric graphs. In this case, we may unite all edge functions $\{w_j\}_{j=1}^{|E|}$ into one single edge function $w : E \to \R$, where $w \vert_{e_j} = w_j$ for $1 \leq j \leq |E|$.
\end{remark}

We refer to Appendix \ref{sec:derivation} for a formal derivation of \eqref{w:edge} in which a reaction-diffusion equation in a thin domain in between two bulk domains shrinks to an edge, and also for a formal derivation of \eqref{pop:bc} in which a reaction-diffusion equation in a thin interval between two edges shrinks to a vertex. For the benefit of the reader, we conclude this section with a simple example to illustrate the structure of \eqref{pop:bc}.

\begin{figure}[h]
    \centering
\includegraphics[width=0.7\textwidth]{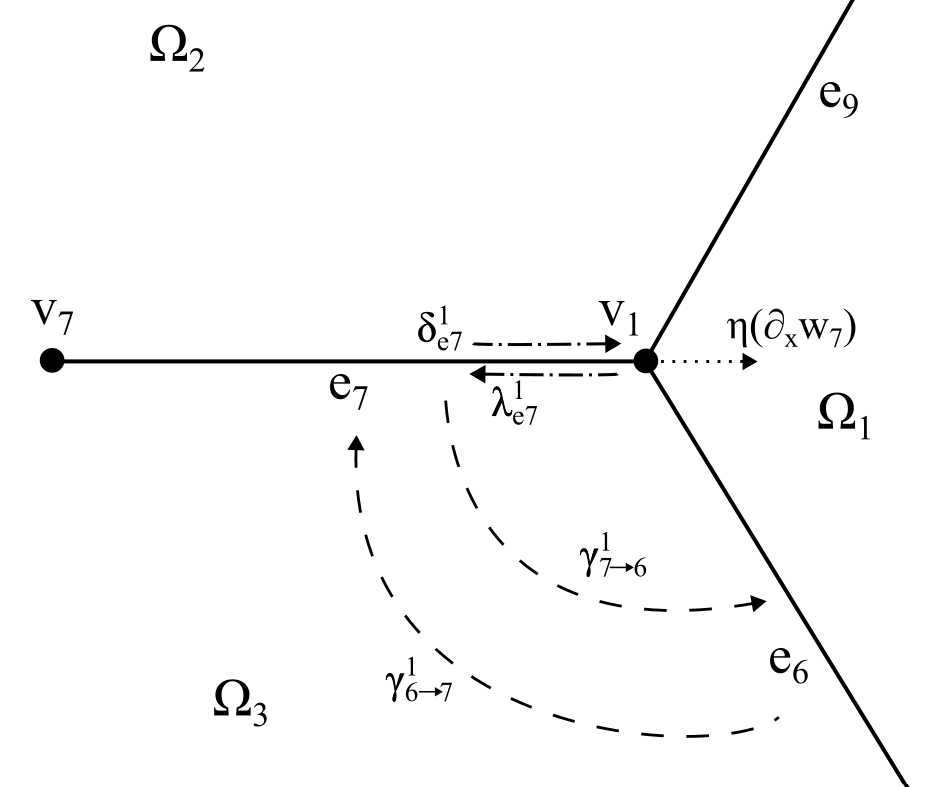}
    \caption{Balance of fluxes at vertex $v_1$ from Figure~\ref{fig:domain}.}
    \label{fig:vertex_jun}
\end{figure}  

\begin{ex}
Let us focus on vertex $v_1$ in the setting of Figure~\ref{fig:domain}, which connects edges $e_6$, $e_7$ and $e_9$. As shown in Figure~\ref{fig:vertex_jun}, at $v_1$ we have outward fluxes represented by dotted arrows, the exchange of material across $v_1$ into connected edges is represented by dashed arrows, and the exchange of material between edge and vertex is represented by dot-dashed arrows. Then, the balance of fluxes at $v_1$ for edges $e_6$, $e_7$ and $e_9$ read as the following:
\begin{align*}
\eta(\pd_x w_6) & = -\delta^1_{e_6} w_6 + \lambda^{1}_{e_6} z_1 + \gamma^1_{7 \to 6} w_7 + \gamma^1_{9 \to 6} w_9 - \gamma^1_{6 \to 7} w_6 - \gamma^1_{6 \to 9} w_6, \\
\eta(\pd_x w_7) & = -\delta^1_{e_7} w_7 + \lambda^1_{e_7} z_1 + \gamma^1_{6 \to 7} w_6 + \gamma^1_{9 \to 7} w_9 - \gamma^1_{7 \to 6} w_7 - \gamma^1_{7 \to 9} w_7, \\
\eta(\pd_x w_9) & = -\delta^1_{e_9} w_9 + \lambda^1_{e_9} z_1 + \gamma^1_{6 \to 9} w_6 + \gamma^1_{7 \to 9} w_7 - \gamma^1_{9 \to 6} w_9 + \gamma^1_{9 \to 7} w_9.
\end{align*}
Expressing in matrix-vector form:
\begin{align*}
  &  \begin{pmatrix}
        \eta(\pd_x w_6) \\
        \eta(\pd_x w_7) \\
        \eta(\pd_x w_9)
    \end{pmatrix} + \left (\begin{pmatrix}
      \gamma^1_{6 \to 7}+\gamma^1_{6 \to 9} & -\gamma^1_{6 \to 7} & -\gamma^1_{6 \to 9} \\
      -\gamma^1_{6 \to 7} & \gamma^1_{7 \to 6}+\gamma^1_{7 \to 9} & -\gamma^1_{9 \to 7}, \\
      -\gamma^1_{6 \to 9} & -\gamma^1_{7 \to 9} & \gamma^1_{9 \to 6}+\gamma^1_{9 \to 7}
    \end{pmatrix} +\begin{pmatrix}
        \delta^1_{e_6} & 0 & 0 \\
        0 & \delta^1_{e_7} & 0 \\
        0 & 0 & \delta^1_{e_9}
    \end{pmatrix} \right ) \begin{pmatrix}
        w_6 \\ w_7 \\ w_9
    \end{pmatrix} \\
    & \qquad = \bm{\eta}(\pd_x \bm{w}_{v_1}) + (\bm{N}_{v_1} + \bm{E}_{v_1}) \bm{w}_{v_1} = \begin{pmatrix}
        \lambda^1_{e_6} z_1 \\ \lambda^1_{e_7} z_1 \\ \lambda^1_{e_9}z_1
    \end{pmatrix},
\end{align*}
where we note that the column and row sums of $\bm{N}_{v_1}$ are zero.
\end{ex}

\section{Main results}\label{sec:main}
The standard Lebesgue and Sobolev spaces over $\Omega_i$ are denoted by $L^p := L^p(\Omega_i)$ and $W^{k,p} := W^{k,p}(\Omega_i)$ for any $p \in [1,\infty]$ and $k > 0$, with corresponding norms $\| \cdot \|_{L^p(\Omega_i)}$ and $\| \cdot \|_{W^{k,p}(\Omega_i)}$. In the case $p = 2$, these become Hilbert spaces and we use the notation $H^k := H^{k}(\Omega_i) = W^{k,2}(\Omega_i)$, along with the norm $\| \cdot \|_{H^k(\Omega_i)}$ and inner product $(\cdot,\cdot)_{H^k(\Omega_i)}$. Similar notations are used for Lebesgue and Sobolev spaces over edge $e_j$. The dual space of a Banach space $X$ is denoted as $X^*$ and we denote the duality pairing between $X^*$ and $X$ by $\langle \cdot, \cdot \rangle_X$.

\subsection{Assumptions}
We make the following assumptions.
\begin{enumerate}[label=$\boldsymbol{(\mathrm{A \arabic*})}$, ref=$\boldsymbol{\mathrm{A \arabic*}}$]
\item \label{ass:aniso} The function $\widehat{\kappa}: \R^2 \to \R$ is continuously differentiable, positive on $\R^2 \setminus \{ \bm{0} \}$, and for a fixed index $p \in [2,\infty)$ there exist positive constants $c_0$, $c_1$, $c_2$ and $c_3$ such that 
    \[
    c_0 |\bm{q}|^p \leq \widehat{\kappa}(\bm{q}) \leq c_1 |\bm{q}|^p, \quad \forall \bm{q} \in \R^2, 
    \]
    and for $\nabla \widehat{\kappa} = \kappa \in C^0(\R^2;\R^2)$,
    \[
   \kappa (\bm{0}) = \bm{0}, \quad  |\kappa(\bm{q})| \leq c_2 |\bm{q}|^{p-1}, \quad (\kappa(\bm{q}) - \kappa(\bm{r})) \cdot (\bm{q} - \bm{r}) \geq c_3 |\bm{q} - \bm{r}|^{p}, \quad \forall \bm{q}, \bm{r} \in \R^2.
    \]
    Consequently, $\widehat{\kappa}$ is strongly convex and thus strictly convex.
\item \label{ass:edgeaniso} The function $\widehat{\eta}: \R \to \R$ is continuously differentiable, $\widehat{\eta}(0) = 0$, $\widehat{\eta}(s) > 0$ for $s \neq 0$, and for a fixed index $p \in [2,\infty)$ there exist positive constants $c_4$, $c_5$, $c_6$ and $c_7$ such that 
    \[
    c_4 |s|^p \leq \widehat{\eta}(s) \leq c_5 |s|^p, \quad \forall s \in \R, 
    \]
    and for $\widehat{\eta}'  = \eta \in C^0(\R)$,
    \[
   \eta(0) = 0, \quad  |\eta(s)| \leq c_6 |s|^{p-1}, \quad (\eta(s) - \eta(r)))(s-r) \geq c_7 |s-r|^{p}, \quad \forall s,r \in \R.
    \]
    Consequently, $\widehat{\eta}$ is strongly convex and thus strictly convex.
\item \label{ass:f} The function $f: \R \to \R$ is monotone and is the antiderivative of a non-negative convex function $\widehat{f}: \R \to \R$ satisfying $\widehat{f} \in C^1(\R)$ and 
    \[
    \widehat{f}(s) \leq c_8 (1 + |s|^{q+1}), \quad |f(s)| \leq c_9 (1 + |s|^{q}), \quad \forall s \in \R,
    \]
for exponent 
$q < \frac{2(p-1)^2}{p} + 2$ and positive constants $c_8$ and $c_9$. Without loss of generality we further assume $f(0) = 0$.
\item \label{ass:g} The function $g: \R \to \R$ is monotone and is the antiderivative of a non-negative convex function $\widehat{g}: \R \to \R$ satisfying $\widehat{g} \in C^1(\R)$ and 
    \[
    \widehat{g}(s) \leq c_{10} (1 + |s|^{r+1}), \quad |g(s)| \leq c_{11} (1 + |s|^{r}), \quad \forall s \in \R,
    \]
for some exponent 
$r < \frac{(3p-2)(p-1)}{p} + 2$ and positive constants $c_{10}$ and $c_{11}$. Without loss of generality we further assume $g(0) = 0$.
\item \label{ass:coeff} For all $e_j \sim \Omega_i$, the coefficients $\alpha_{ij}$ and $\beta_{ij}$ are positive. For all $e_j \sim v_k$, the coefficients $\delta^k_{e_j}$ and $\lambda^k_{e_j}$ are  positive, while the coefficients $\gamma_{j \to m}^{k}$ for $e_j, e_m \sim v_k$ are non-negative and satisfy
\[
\delta^{k}_{e_j} + \sum_{e_m \sim v_k, m \neq j}^{|E|} \gamma_{j \to m}^k \geq \sum_{e_m \sim v_k, m \neq j}^{|E|} \gamma^k_{m \to j}.
\]
\item \label{ass:ini} The initial data satisfy
\[
u_{i,0} \in \begin{cases}
   L^2(\Omega_i) & \text{ if } p = 2, \\
   W^{1,p}(\Omega_i) & \text{ if } p \in (2,\infty), 
\end{cases} \quad w_{j,0} \in \begin{cases}
L^2(e_j) & \text{ if } p = 2, \\
W^{1,p}(e_j) & \text{ if } p \in (2,\infty).
\end{cases}
\]
\end{enumerate}
We point out that \eqref{ass:coeff} ensures the matrix $\bm{N}_{v_k} + \bm{E}_{v_k}$ is diagonally dominant. 
\begin{defn}\label{defn:weaksoln}
A weak solution to \eqref{dom:equ}-\eqref{w:edge}, \eqref{pop:equ}-\eqref{ini} is a collection of functions $\bm{u} = (u_1, \dots, u_M)$, $\bm{w} = (w_1, \dots, w_{|E|})$ and $\bm{z}= (z_1, \dots, z_{|V|})$ satisfying for $1 \leq i \leq M$, $1 \leq j \leq |E|$, $1 \leq k \leq |V|$,
\begin{align*}
u_i & \in L^\infty(0,T;L^2(\Omega_i)) \cap L^p(0,T;W^{1,p}(\Omega_i)), \\
\pd_{t} u_i & \in L^{p/(p-1)}(0,T;W^{1,p}(\Omega_i)^*), \\
w_j & \in L^\infty(0,T;L^2(e_j)) \cap L^p(0,T;W^{1,p}(e_j)), \\
\pd_t w_j & \in L^{p/(p-1)}(0,T;W^{1,p}(e_j)^*), \\
z_k & \in H^1(0,T),
\end{align*}
together with $z_k(0) = z_{k,0}$,
\begin{align*}
\langle u_{i}(0), \phi \rangle_{W^{1,p}(\Omega_i)} & = \langle u_{0,i}, \phi \rangle_{W^{1,p}(\Omega_i)}, \quad \forall \phi \in W^{1,p}(\Omega_i), \\
\langle w_{j}(0), \psi \rangle_{W^{1,p}(e_j)} & = \langle w_{0,j}, \psi \rangle_{W^{1,p}(e_j)}, \quad \forall \psi \in W^{1,p}(e_j),
\end{align*}
and for a.e. $t \in (0,T)$ and for all $\phi^{(i)} \in W^{1,p}(\Omega_i)$, $\psi^{(j)} \in W^{1,p}(e_j)$, it holds that
\begin{equation}\label{weakform}
\begin{aligned}
0 & = \sum_{i=1}^M \Big (  \langle \pd_t u_{i},  \phi^{(i)} \rangle_{W^{1,p}(\Omega_i)} + \int_{\Omega_i} \kappa(\nabla u_{N}) \cdot \nabla \phi^{(i)} + f(u_{i}) \phi^{(i)} \, dx \Big ) \\
& \quad + \sum_{j=1}^{|E|} \Big ( \langle \pd_t w_{j},  \psi^{(j)} \rangle_{W^{1,p}(e_j)} + \int_{e_j} \eta(\pd_x w_{j}) \pd_x \psi^{(j)} + g(w_{j}) \psi^{(j)}  \, dx \\
& \qquad  \qquad + \sum_{i \,  :\, \Omega_i \sim e_j}^M \int_{e_j} (\alpha_{ij} w_{j} - \beta_{ij} u_{i})(\psi^{(j)} - \phi^{(i)}) \, dx \Big ) \\
& \quad + \sum_{j=1}^{|E|} \sum_{k \, : \, v_k \sim e_j}^{|V|} \Big (\big [\delta^{k}_{e_j}  + \sum_{e_m \sim v_k, m \neq j}^{|E|} \gamma^k_{j \to m} \big ] w_j - \lambda^k_{e_j} z_k - \sum_{e_m \sim v_k, m \neq j}^{|E|} \gamma^k_{m \to j} w_{m} \Big ) \psi^{(j)} \Big \vert_{v_k},
\end{aligned}
\end{equation}
and
\begin{align}\label{ODE}
\frac{d}{dt} z_k = \sum_{j \, : \, e_j \sim v_k}^{|E|} (\delta^k_{e_j} w_j(t,v_k) - \lambda^k_{e_j} z_k).
\end{align}
\end{defn}

\noindent 
Let us make a few remarks:
\begin{itemize}
    \item In \eqref{weakform} and \eqref{ODE} we evaluate the edge functions $w_j$ and the edge test functions $\psi^{(j)}$ at the vertex $v_k$, which are well-defined by the use of the trace theorem for Sobolev functions;
    \item We choose to present the junction conditions \eqref{pop:bc} as an enumeration over the edges, as opposed to enumerating over the vertices, to enable a clearer perspective for the terms involving the edge functions $w_j$. But we remark that both formulations are equivalent since exchanging the order of summation yields the same result:
    \[
\sum_{j=1}^{|E|} \sum_{k \, : \, v_k \sim e_j}^{|V|} \delta^k_{e_j} w_j \psi^{(j)} \Big \vert_{v_k} = \sum_{k=1}^{|V|} \sum_{j \, : \, e_j \sim v_k}^{|E|} \delta^k_{e_j} w_j \psi^{(j)} \Big \vert_{v_k}.
    \]
\item We can combine \eqref{weakform} and \eqref{ODE} together to provide an alternate variational equality:
\begin{equation}\label{weakform:alt}
\begin{aligned}
0 & = \sum_{i=1}^M \Big (  \langle \pd_t u_{i},  \phi^{(i)} \rangle_{W^{1,p}(\Omega_i)} + \int_{\Omega_i} \kappa(\nabla u_{N}) \cdot \nabla \phi^{(i)} + f(u_{i}) \phi^{(i)} \, dx \Big ) \\
& \quad + \sum_{j=1}^{|E|} \Big ( \langle \pd_t w_{j},  \psi^{(j)} \rangle_{W^{1,p}(e_j)} + \int_{e_j} \eta(\pd_x w_{j}) \pd_x \psi^{(j)} + g(w_{j}) \psi^{(j)}  \, dx \\
& \qquad  \qquad + \sum_{i \,  :\, \Omega_i \sim e_j}^M \int_{e_j} (\alpha_{ij} w_{j} - \beta_{ij} u_{i})(\psi^{(j)} - \phi^{(i)}) \, dx \Big ) \\
& \quad + \sum_{k=1}^{|V|} \sum_{j \neq m \, : \, e_j, e_m \sim v_k}^{|E|} \Big [ (\gamma^k_{j \to m} w_j -  \gamma^k_{m \to j} w_{m} ) (\psi^{(j)} - \psi^{(m)}) \Big ] \Big \vert_{v_k} \\
& \quad + \sum_{k=1}^{|V|} \big (\frac{d}{dt} z_k \big) \theta_k + \sum_{j=1}^{|E|} \sum_{k \, : \, v_k \sim e_j}^{|V|} (\delta_{e_j}^{k} w_j \vert_{v_k} - \lambda_{e_j}^k z_k)(\psi^{(j)} \vert_{v_k} - \theta_k),
\end{aligned}
\end{equation}
for a.e. $t \in (0,T)$ and for all $\phi^{(i)} \in W^{1,p}(\Omega_i)$, $\psi^{(j)} \in W^{1,p}(e_j)$, $\theta_k \in \R$. Here we can better see the exchange mechanism between subregion and edge, as well as between adjacent edges sharing the same vertex and between edge and vertex.
\item We point out that in the special case where $\alpha_{ij} = \beta_{ij}$ for all $e_j \sim \Omega_i$, $\gamma_{j \to m}^k = \gamma_{m \to j}^k$ for all $e_j, e_m \sim v_k$ and $\delta^k_{e_j} = \lambda^k_{e_j}$ for all $e_j \sim v_k$, then by choosing $\phi^{(i)} = 1$, $\psi^{(j)} = 1$ and $\theta_k = 1$ in \eqref{weakform:alt} we derive the conservation of total mass:
\[
\frac{d}{dt} \Big ( \sum_{i=1}^M \int_{\Omega_i} u_i \, dx + \sum_{j=1}^{|E|} \int_{e_j} w_j \, dx + \sum_{k=1}^{|V|} z_k \Big ) = 0.
\]
\item Although \eqref{ass:coeff} implicitly assumes all vertices are populated, we can account for some vertices being unpopulated by setting the corresponding coefficients at that vertex (say $v_k$) to satisfy $\delta^k_{e_j} = \lambda^k_{e_j} = 0$ for all $e_j \sim v_k$ and $\gamma^k_{j \to m} = \gamma^k_{m \to j}$ for all $e_j, e_m \sim v_k$, $j \neq m$. Then, we see from \eqref{weakform} that any contributions from material transfer between $v_k$ and its associated edges vanish, leaving only inter-edge transfers through the unpopulated vertex.
\end{itemize}

\subsection{Well-posedness, temporal regularity and boundedness of solutions}
Our first result is the well-posedness of the mixed-dimensional spatial-network model.
\begin{thm}[Well-posedness of weak solutions]\label{thm:well}
Under \eqref{ass:aniso}--\eqref{ass:ini}, there exists a unique weak solution $(\bm{u}, \bm{w}, \bm{z})$ in the sense of Definition \ref{defn:weaksoln} that depend continuously on the initial data.
\end{thm}

In the setting where at all vertices $v_k$, the inter-edge transfer coefficients are symmetric, i.e., 
\begin{align}\label{ass:sym:coeff}
\gamma_{j \to m}^k = \gamma_{m \to j}^k, \quad \text{ for all } e_j, e_m \sim v_k, \ m \neq j,
\end{align}
and using the notation $\gamma^k_{j \leftrightarrow m}$ to denote this coefficient of material exchange between edges $e_j$ and $e_m$ via vertex $v_k$, we are able to establish further properties for the weak solution.

\begin{thm}[Temporal regularity under symmetric inter-edge transfer coefficients]\label{thm:reg}
Suppose the condition \eqref{ass:sym:coeff} holds at all vertices $\{v_k\}_{k=1}^{|V|}$, then the weak solution in Theorem \ref{thm:well} possesses the additional regularities:
\begin{align*}
u_i & \in L^\infty(0,T;W^{1,p}(\Omega_i)) \cap H^1(0,T;L^2(\Omega_i)) \cap L^\infty(0,T;L^\infty(\Omega_i)), \\
w_j & \in L^\infty(0,T;W^{1,p}(e_j)) \cap H^1(0,T;L^2(e_j)) \cap L^\infty(0,T;L^\infty(e_j)), \\
z_{k} & \in W^{1,\infty}(0,T), 
\end{align*}
and the initial data are attained almost everywhere in the sense that $u_{i}(0) = u_{i,0}$ a.e.~in $\Omega_i$ and $w_{j}(0) = w_{j,0}$ a.e.~in $e_j$.
\end{thm}
\begin{remark}\label{rem:bdd}
The boundedness of edge solutions arises due to the Sobolev embedding $W^{1,p}(e_j) \subset C^0(e_j)$ for $p \in [2,\infty)$. For the subdomain solutions, we use a comparison principle to establish boundedness.
\end{remark}

\subsection{Spatial regularity of solutions}
 
With regards to spatial regularity, we remark that the presence of the anisotropies $\kappa$ and $\eta$ complicate the application of standard elliptic regularity results. Although \cite{GKW} provides interior (Theorem 5.2 therein) and boundary (Theorem 5.4 therein) $H^2$-regularity in the setting $p = 2$ and under  additional Lipschitz continuity of the anisotropy function, we cannot apply these results to our present setting for the following reasons:
\begin{itemize}
    \item The regularity result of \cite{GKW} is for homogeneous Neumann boundary conditions, and the extension to non-homogeneous Neumann or Robin boundary conditions seem non-trivial.
    \item The equation \eqref{u:bdy} should be seen as a mixed boundary condition as different edge functions need not be related to each other even if these edges compose the boundary of the same subdomain $\Omega_i$.
    \item In general the boundary $\pd \Omega_i$ of a subdomain $\Omega_i$ is polygonal with possible non-convex angles, which does not meet the $C^{1,1}$ requirement of \cite[Theorem 5.4]{GKW}.
\end{itemize}

These issues can be somewhat elevated in the setting where all the vertices are unpopulated, and as remarked in Remark \ref{rem:Kirchhoff}, we may represent $w_j$ as the restriction of a global edge function $w : E \to \R$ on edge $e_j$. 

\begin{thm}[Spatial regularity in the unpopulated vertices setting]
Suppose all vertices are unpopulated and the condition \eqref{ass:sym:coeff} holds at all vertices $\{v_k\}_{k=1}^{|V|}$. In addition, let \eqref{ass:edgeaniso} hold with $p = 2$ and suppose that $\eta$ is Lipschitz continuous with Lipschitz constant $C_L>0$:
\[
|\eta(s_1) - \eta(s_2)| \leq C_L |s_1 - s_2|, \quad \forall s_1, s_2 \in \R.
\]
Then, the global edge function $w$ satisfies
\[
w \in L^2(0,T;H^2(E)).
\]
Furthermore, if $\Omega_i$ is a convex subdomain and $\kappa(\bm{s}) = \frac{1}{2}|\bm{s}|^2$, then the corresponding solution $u_i$ satisfies
\[
u_i \in L^2(0,T;H^2(\Omega_i)).
\]
\end{thm}

\subsection{Finite time extinction}
In the unpopulated vertices setting with symmetric transfer coefficients and reaction terms satisfying suitable growth conditions, we can establish finite time extinction for the subdomain and edge solutions. We note that the $H^2$-regularity from Theorem \ref{thm:reg} is not required for the following result.
\begin{thm}[Finite time extinction under symmetric coefficients and unpopulated vertices]\label{thm:finitetime}
Suppose that the transfer coefficients satisfy
\begin{align}\label{finitetime:coeff}
\alpha_{ij} = \beta_{ij} \text{ for } \Omega_i \sim e_j, \quad \gamma_{j \to m}^k = \gamma_{m \to j}^k \text{ for } e_j, e_m \sim v_k, \quad \delta^k_{e_j} = \lambda^k_{e_j} = 0 \text{ at all } v_k,
\end{align}
and the reaction terms $f$ and $g$ satisfy
\begin{align}\label{fg:finitetime}
f(s) s \geq c_{12} |s|^\sigma, \quad g(s) s \geq c_{13}|s|^{\sigma}, \quad \forall s \in \R,
\end{align}
for some positive constants $c_{12}$ and $c_{13}$ and for some exponent $1 < \sigma < 2$. Then, there exists a time $t_*$ depending on the initial data, the subdomains and edges, and the system parameters such that 
\[
u_i(t) = 0, \quad w_j(t) = 0, 
\]
for all $t > t_*$, and for all $1 \leq i \leq M$ and $1 \leq j \leq |E|$.
\end{thm}

We remark that the condition \ref{fg:finitetime} can be satisfied for instance with the choice $f(s) = c_{12}|s|^{\sigma-2}s$ and $g(s) = c_{13}|s|^{\sigma - 2}s$.

\section{Existence of weak solutions}\label{sec:wellposed}

\subsection{Semi-Galerkin scheme}
We employ a Galerkin approximation and consider orthonormal bases $\{\phi_k^{(i)}\}_{k\in \N}$ of $L^2(\Omega_i)$ for $i = 1, \dots, M$ that are orthogonal in $H^1(\Omega_i)$, as well as orthonormal bases $\{\psi_m^{(j)}\}_{m \in \N}$ of $L^2(e_j)$ for $j = 1, \dots, |E|$ that are orthogonal in $H^1(e_j)$. Commonly used examples are eigenfunctions of the Laplacian in $\Omega_i$ or on $e_j$, which by elliptic regularity theory on polygonal domains we have that $\phi^{(i)}_k \in H^2(\Omega_i) \subset W^{1,p}(\Omega_i)$ and $\psi^{(j)}_m \in H^2(e_j) \subset W^{1,p}(e_j)$.  For $N \in \N$, we introduce the finite dimensional spaces $\mathbb{U}^{(i)}_N := \mathrm{span}\{\phi_1^{(i)}, \dots, \phi_N^{(i)}\}$ and $\mathbb{W}^{(j)}_N := \mathrm{span}\{\psi_1^{(j)}, \dots, \psi_N^{(j)}\}$, along with their orthogonal projections $P_{\mathbb{U}^{(i)}_N}$ and $P_{\mathbb{W}^{(j)}_N}$, respectively. Furthermore, for fixed $N \in \N$, there exist positive constants $C_N$ depending on $N$ such that the following inverse estimates hold:
\begin{align}\label{inv:est}
\| y \|_{H^1(\Omega_i)} \leq C_N \| y \|_{L^2(\Omega_i)}, \quad \| z \|_{H^1(e_j)} \leq C_N \| z \|_{L^2(e_j)}
\end{align}
for all $y \in \mathbb{U}^{(i)}_N$, $1 \leq i \leq M$ and $z \in \mathbb{W}^{(j)}_N$, $1 \leq j \leq |E|$. 

\subsubsection{Solving the vertex ODEs}
We define 
\[
\Lambda^k := \sum_{j \, : \, e_j \sim v_k}^{|E|} \lambda^k_{e_j}, \quad k = 1, \dots, |V|.
\]
Given $W_j \in C^0([0,T];\mathbb{V}_{N}^{(j)})$ for $j = 1, \dots, |E|$, we solve the ordinary differential equation 
\begin{align}\label{Gal:ODE}
\frac{d}{dt} z_{k,W} = - \Lambda^k z_{k,W} + \sum_{j \, : \, e_j \sim v_k}^{|E|} \delta^k_{e_j} W_j(t, v_k),
\end{align}
furnished with initial value $z_{k,W}(0) = z_{k,0}$ for $k = 1, \dots, |V|$.  The unique solution is
\begin{align*}
    z_{k,W}(t) = z_{k,0}e^{ - \Lambda^k t} +  \sum_{j \, : \, e_j \sim v_k}^{|E|} \delta^k_{e_j} \int_0^t W_j(s, v_k)e^{ \Lambda^k (s-t) }\, ds, \quad t \in [0,T],
\end{align*}
and satisfies (recalling that $D_k$ is the degree of vertex $v_k$)
\begin{align}\label{zkW:est}
|z_{k,W}(t)|^2 \leq 2|z_{k,0}|^2e^{- 2\Lambda^k t} + 2 D_k \sum_{j\, : \, e_j \sim v_k}^{|E|} \frac{(\delta^k_{e_j})^2}{2 \Lambda^k} \int_0^t | W_j (s,v_k)|^2 \, ds.
\end{align}
It is clear that $z_{k,W} \in C^1([0,T])$ for $k = 1, \dots, |V|$, and we denote by $\bm{z}_W$ the vector of ODE solutions $\bm{z}_{W} = (z_{1, W}, \dots, z_{|V|, W}) \in C^1([0,T];\R^{|V|})$.

\subsubsection{Solving the subdomain-edge PDE system}
For fixed $N \in \N$ we consider Galerkin solutions
\[
u_{i,N} = \sum_{k=1}^N \mathrm{a}_k^{(i)}(t) \phi_k^{(i)}(x), \quad w_{j,N} = \sum_{m=1}^N \mathrm{b}_m^{(j)}(t) \psi_m^{(j)}(x)
\]
satisfying for all $1 \leq r,s \leq N$
\begin{subequations}\label{aux:prob:Galerkin}
\begin{alignat}{2}
\label{Gal:u} 0 & = \sum_{i=1}^M \int_{\Omega_i} \pd_t u_{i,N} \phi_r^{(i)} + \kappa(\nabla u_{i,N}) \cdot \nabla \phi_r^{(i)} + f(u_{i,N}) \phi_r^{(i)} \, dx \\
\notag & \qquad \qquad - \sum_{j \, : \, e_j \sim \Omega_i}^{|E|}\int_{e_j} (\alpha_{ij} w_{j,N} - \beta_{ij} u_{i,N}) \phi_r^{(i)} \, dx \\
\label{Gal:w} 0 & = \sum_{j=1}^{|E|} \Big ( \int_{e_j} \pd_t w_{j,N} \psi_{s}^{(j)} + \eta(\pd_x w_{j,N}) \pd_x \psi_{s}^{(j)} + g(w_{j,N}) \psi_{s}^{(j)}  \, dx \\
\notag & \qquad  \qquad + \sum_{i \, : \, \Omega_i \sim e_j}^{M} \int_{e_j} (\alpha_{ij} w_{j,N} - \beta_{ij} u_{i,N})\psi_{s}^{(j)} \, dx \Big ) \\
\notag & \quad + \sum_{j=1}^{|E|} \sum_{k \, : \, v_k \sim e_j}^{|V|} \Big (\big [\delta^{k}_{e_j} + \sum_{e_m \sim v_k, m \neq j} \gamma^k_{j \to m} \big ]  w_{j,N} - \lambda^k_{e_j} z_{k,W} - \sum_{e_m \sim v_k, m \neq j} \gamma^k_{m \to j} w_{m,N} \Big ) \psi_s^{(j)} \Big \vert_{v_k},
\end{alignat}
\end{subequations}
and furnished with initial data
\[
u_{i,N}(0) = P_{\mathbb{U}^{(i)}_N}(u_{i,0}), \quad w_{j,N}(0) = P_{\mathbb{W}^{(j)}_N}(w_{j,0}).
\]
We can express \eqref{aux:prob:Galerkin} as a system of ordinary differential equations for the coefficients $\mathrm{a}_k^{(i)}$ and $\mathrm{b}_m^{(j)}$, with right-hand sides depending continuously on the coefficients due to the continuity of $\kappa$, $\eta$, $f$ and $g$.  Hence, by the Cauchy--Lipschitz theorem, there exists a time $T_N$ and solutions $\mathrm{a}_k^{(i)} \in C^1([0,T_N])$ and $\mathrm{b}_m^{(j)} \in C^1([0,T_N])$ for $1 \leq k,m  \leq N$, $1 \leq i \leq N$ and $1 \leq j \leq |E|$, which in turn provides the existence of unique Galerkin solutions $u_{i,N} \in C^1([0,T_N]; \mathbb{U}^{(i)}_N)$ and $w_{j,N} \in C^1([0,T_N]; \mathbb{W}^{(j)}_N)$. 

\subsubsection{Existence of semi-Galerkin solution}
We now introduce the mapping 
\[
\mathcal{F}_N : \bm{W} = (W_1, \dots, W_{|E|}) \to \bm{w}_N = (w_{1,N}, \dots, w_{|E|,N}),
\]
where $\{w_{j,N}\}_{j=1}^{|E|}$ is the edge Galerkin solutions to \eqref{aux:prob:Galerkin} corresponding to the vector $\bm{z}_W$, the latter itself as the unique solution to the ordinary differential systems \eqref{Gal:ODE} corresponding to the input vector $\bm{W}$.  Our aim is to apply Schauder's fixed point theorem to deduce the existence of a fixed point $\hat{\bm{w}}_N \in C^0([0,T];\prod_{j=1}^{|E|} \mathbb{W}_{N}^{(j)})$ of $\mathcal{F}_N$, which generates a solution $(\hat{\bu}_N, \hat{\bm{w}}_N, \hat{\bm{z}}_N)$ to the following system expressed in strong formulation:
\begin{subequations}\label{fixedpoint:sys}
\begin{alignat}{2}
\label{fp:u} \pd_t \hat{u}_{i,N} - P_{\mathbb{U}^{(i)}_N}[\div (\kappa(\nabla \hat{u}_{i,N})) - f(\hat{u}_{i,N})] = 0 & \quad \text{ in } \Omega_i \times (0,T), \\
\kappa(\nabla \hat{u}_{i,N}) \cdot \bm{n} = \alpha_{ij} \hat{w}_{j,N} - \beta_{ij} \hat{u}_{i,N} & \quad \text{ on } e_j \times (0,T), \quad \Omega_i \sim e_j, \\
\notag \pd_t \hat{w}_{j,N} - P_{\mathbb{W}^{(j)}_N}[\pd_x (\eta(\pd_x \hat{w}_{j,N})) - g(\hat{w}_{j,N})] \qquad  & \\
\label{fp:w} + \sum_{\Omega_m \sim e_j} \alpha_{mj} \hat{w}_{j,N} - \beta_{mj} \hat{u}_{m,N} = 0 & \quad \text{ on } e_j \times (0,T), \\
\bm{\eta}(\pd_x \hat{\bm{w}}_{v_k,N}) + (\bm{N}_{v_k} + \bm{E}_{v_k}) \hat{\bm{w}}_{v_k,N} =  \hat{z}_{k,N} \bm{\lambda}_{v_k} & \quad \text{ on } v_k \times (0,T), \quad e_j \sim v_k, \\
\label{fp:z} \frac{d}{dt} \hat{z}_{k,N} = \sum_{j \, : \, e_j \sim v_k}^{|E|} (\delta^k_{e_j} \hat{w}_{j,N}(t, v_k) - \lambda^k_{e_j} \hat{z}_{k,N}) & \quad \text{ in } (0,T).
\end{alignat}
\end{subequations}
This is achieved by demonstrating
\begin{enumerate}
\item[(i)] there exists a suitable set $\mathcal{S}_N \subset C^0([0,T_N];\prod_{j=1}^{|E|} \mathbb{W}^{(j)}_N)$ where $\mathcal{F}_N : \mathcal{S}_N \to \mathcal{S}_N$; 
\item[(ii)] $\overline{\mathcal{F}_N(\mathcal{S}_N)}$ is compact in $C^0([0,T_N];\prod_{j=1}^{|E|} \mathbb{W}^{(j)}_N)$; \item[(iii)] $\mathcal{F}_N$ is continuous. 
\end{enumerate}
In the sequel we use the symbol $C$ to denote positive constants whose values may change line to line and even within the same line. Multiplying \eqref{Gal:u} by $\mathrm{a}_r^{(i)}$ and \eqref{Gal:w} by $\mathrm{b}_s^{(j)}$ and summing over $1 \leq r, s \leq N$ leads to
\begin{equation}\label{Gal:est}
\begin{aligned}
& \sum_{i=1}^M \Big (\frac{d}{dt} \frac{1}{2} \| u_{i,N} \|^2_{L^2(\Omega_i)} + \int_{\Omega_i} \kappa(\nabla u_{i,N}) \cdot \nabla u_{i,N} + f(u_{i,N}) u_{i,N} \, dx \Big ) \\
& \quad + \sum_{j=1}^{|E|} \Big (\frac{d}{dt} \frac{1}{2} \| w_{j,N} \|^2_{L^2(e_j)} + \int_{e_j} \eta(\pd_x w_{j,N}) \pd_x w_{j,N} + g(w_{j,N}) w_{j,N} \, dx \\
& \qquad \qquad + \sum_{i \, : \, \Omega_i \sim e_j}^{M} \int_{e_j} \alpha_j |w_{j,N}|^2 + \beta_{ij} |u_{i,N}|^2 - (\alpha_{ij} + \beta_{ij}) w_{j,N} u_{i,N}\, dx \Big ) \\
& \quad + \sum_{j=1}^{|E|} \sum_{k \, : \, v_k \sim e_j}^{|V|} \Big ( \big [\delta^{k}_{e_j} + \sum_{e_m \sim v_k, m \neq j} \gamma_{j \to m}^k \big ]  |w_{j,N}(v_k)|^2 \\
& \qquad \qquad \qquad - \sum_{e_m \sim v_k, m \neq j} \gamma_{m \to j}^k w_{m,N}(v_k) w_{j,N}(v_k) \Big ) \\
& = \sum_{j=1}^{|E|}\sum_{k \, : \, v_k \sim e_j}^{|V|} \lambda^k_{e_j} z_{k,W} w_{j,N}(v_k).
\end{aligned}
\end{equation}
Invoking the coercivity of $\kappa$ and $\eta$, as well as the monotonicity of $f$ and $g$ we have
\begin{align*}
\int_{\Omega_i} \kappa(\nabla u_{i,N}) \cdot \nabla u_{i,N} + f(u_{i,N}) u_{i,N} \, dx & \geq c_3 \| \nabla u_{i,N} \|_{L^p(\Omega_i)}^p, \\
\int_{e_j} \eta(\pd_x w_{j,N}) \pd_x w_{j,N} + g(w_{j,N}) w_{j,N} \, dx & \geq c_7 \| \pd_x w_{j,N} \|_{L^p(e_j)}^p.
\end{align*}
For the indefinite term in the third line of \eqref{Gal:est}, we apply Young's inequality 
\begin{align*}
\int_{e_j} (\alpha_{ij} + \beta_{ij}) w_{j,N} u_{i,N} \, dx \leq \frac{\beta_{ij}}{2} \| u_{i,N} \|_{L^2(e_j)}^2 + \frac{(\alpha_{ij} + \beta_{ij})^2}{2 \beta_{ij}} \| w_{j,N} \|_{L^2(e_j)}^2.
\end{align*}
Meanwhile, the Cauchy--Schwarz inequality and a short calculation show that the fourth line of \eqref{Gal:est} can be estimated from below as
\begin{equation}\label{delta:gamma:term}
\begin{aligned}
& \sum_{j=1}^{|E|}\sum_{k \, : \, v_k \sim e_j}^{|V|} \Big ( \big [\delta^{k}_{e_j} + \sum_{e_m \sim v_k, m \neq j} \gamma_{j \to m}^k \big ]  |w_{j,N}(v_k)|^2 - \sum_{e_m \sim v_k, m \neq j} \gamma_{m \to j}^k w_{m,N}(v_k) w_{j,N}(v_k) \Big )  \\
& \quad \geq \sum_{j=1}^{|E|} \sum_{k \, : \, v_k \sim e_j}^{|V|} \Big ( \big [\delta^{k}_{e_j} + \sum_{e_m \sim v_k, m \neq j} (\gamma_{j \to m}^k - \frac{1}{2} \gamma_{m \to j}^k) \big ]  |w_{j,N}(v_k)|^2 \\
 & \qquad \qquad \qquad - \sum_{e_m \sim v_k, m \neq j} \frac{1}{2} \gamma_{m \to j}^k | w_{m,N}(v_k)|^2 \Big ) \\
& \quad = \sum_{j=1}^{|E|} \sum_{k \, : \, v_k \sim e_j}^{|V|} \Big [ \delta^{k}_{e_j} + \frac{1}{2} \sum_{e_m \sim v_k, m \neq j} (\gamma_{j \to m}^k - \gamma_{m \to j}^k) \Big ] |w_{j,N}(v_k)|^2 \\
& \quad \geq \frac{1}{2} \sum_{j=1}^{|E|}\sum_{k \, : \, v_k \sim e_j}^{|V|} \delta^{k}_{e_j} |w_{j,N}(v_k)|^2,
\end{aligned}
\end{equation}
where the last inequality follows from \eqref{ass:coeff}. Lastly, the right-hand side of \eqref{Gal:est} can be estimated from above by Young's inequality
\[
\sum_{j=1}^{|E|}\sum_{k \, : \, v_k \sim e_j}^{|V|} \lambda_{e_j}^k z_{k,W} w_{j,N}(v_k) \leq \frac{1}{4}\sum_{j=1}^{|E|}\sum_{k \, : \, v_k \sim e_j}^{|V|} \delta^{k}_{e_j} |w_{j,N}(v_k)|^2 + \sum_{j=1}^{|E|} \sum_{k \, : \, v_k \sim e_j}^{|V|} \frac{(\lambda^k_{e_j})^2}{\delta_{e_j}^{v_k}} |z_{k,W}|^2.
\]
Hence, we deduce from \eqref{Gal:est} the differential inequality
\begin{equation}\label{diff:ineq}
    \begin{aligned}
& \sum_{i=1}^M \Big (\frac{d}{dt} \frac{1}{2} \| u_{i,N} \|^2_{L^2(\Omega_i)} + c_3 \|\nabla u_{i,N}\|_{L^p(\Omega_i)}^{p} \Big ) \\
& \qquad + \sum_{j=1}^{|E|} \Big (\frac{d}{dt} \frac{1}{2} \| w_{j,N} \|^2_{L^2(e_j)} + c_7 \| \pd_x w_{j,N} \|_{L^p(e_j)}^p \Big ) + \frac{1}{4} \sum_{j=1}^{|E|} \sum_{k \, : \, v_k \sim e_j}^{|V|} \delta^{k}_{e_j} |w_{j,N}(v_k)|^2\\
& \qquad \qquad + \sum_{j=1}^{|E|} \sum_{i \, : \, \Omega_i \sim e_j}^{M} \Big (\frac{\alpha_{ij}}{2} \| w_{j,N} \|_{L^2(e_j)}^2 + \frac{\beta_i}{2} \| u_{i,N} \|_{L^2(e_j)}^2 \Big )  \\
& \quad \leq  \sum_{j=1}^{|E|} \sum_{k \, : \, v_k \sim e_j}^{|V|} \frac{(\lambda^k_{e_j})^2}{\delta_{e_j}^{v_k}} |z_{k,W}|^2  + \sum_{j=1}^{|E|} \sum_{i \, : \, \Omega_i \sim e_j}^{M} \frac{(\alpha_{ij} + \beta_{ij})^2}{2 \beta_{ij}}\| w_{j,N} \|_{L^2(e_j)}^2.
    \end{aligned}
\end{equation}
Let 
\[
G := \max_{j=1,\dots, |E|} \Big ( M \max_{i \, : \, \Omega_i \sim e_j} \frac{(\alpha_{ij} + \beta_{ij})^2}{ \beta_{ij}} \Big ),
\]
then upon applying Gronwall's inequality to \eqref{diff:ineq} and recalling 
\begin{align}\label{ini:bdd}
\| u_{i,N}(0) \|_{L^2(\Omega_i)} \leq \| u_{i,0} \|_{L^2(\Omega_i)}, \quad \| w_{j,N}(0) \|_{L^2(e_j)} \leq \| w_{j,0} \|_{L^2(e_j)},
\end{align}
as well as \eqref{zkW:est} we find that
\begin{align}
\notag & \Big ( \sum_{i=1}^{M} \| u_{i,N}(t) \|_{L^2(\Omega_i)}^2 + \sum_{j=1}^{|E|} \| w_{j,N}(t) \|_{L^2(e_j)}^2 \Big ) \\
\notag & \quad + \int_0^t \frac{e^{G(t-s)}}{2} \sum_{j=1}^{|E|} \sum_{i \, : \, \Omega_i \sim e_j}^{M} \Big (\alpha_{ij} \| w_{j,N} \|_{L^2(e_j)}^2  + \beta_{ij} \| u_{i,N} \|_{L^2(e_j)}^2 \Big ) \, dx \\
\notag & \quad + \int_0^t \frac{e^{G(t-s)}}{2} \sum_{j=1}^{|E|}\sum_{k \, : \, v_k \sim e_j}^{|V|} \delta_{e_j}^{v_k} |w_{j,N}(s,v_k)|^2  \, ds \\
\notag & \leq \Big ( \sum_{i=1}^{M} \| u_{i,0} \|_{L^2(\Omega_i)}^2 + \sum_{j=1}^{|E|} \| w_{j,0} \|_{L^2(e_j)}^2 \Big ) e^{Gt} + \int_0^t 2e^{G(t-s)} \sum_{j=1}^{|E|} \sum_{k \, : \, v_k \sim e_j}^{|V|} \frac{(\lambda^k_{e_j})^2}{\delta_{e_j}^{v_k}} |z_{k,W}(s)|^2 \, ds \\
\notag &  \leq \Big ( \sum_{i=1}^{M} \| u_{i,0} \|_{L^2(\Omega_i)}^2 + \sum_{j=1}^{|E|} \| w_{j,0} \|_{L^2(e_j)}^2 \Big ) e^{Gt} + \int_0^t 4e^{G(t-s)} \sum_{j=1}^{|E|} \sum_{k \, : \, v_k \sim e_j}^{|V|} \frac{(\lambda^k_{e_j})^2}{\delta_{e_j}^{v_k}} |z_{k,0}|^2 e^{- 2\Lambda^k s} \, ds \\
\label{Est:main} & \quad + \int_0^t 4e^{G(t-s)} \sum_{j=1}^{|E|} \sum_{k \, : \, v_k \sim e_j}^{|V|} \frac{(\lambda^k_{e_j})^2}{\delta_{e_j}^{v_k}} \Big ( D_k \sum_{j\, : \, e_j \sim v_k}^{|E|} \frac{(\delta^k_{e_j})^2}{2 \Lambda^k} \int_0^s | W_j (r,v_k)|^2 \, dr \Big )  \, ds.
\end{align}
Setting 
\[
H:= \max_{j=1, \dots, |E|} \Big (|V| \max_{k \, : \, v_k \sim e_j} \frac{(\lambda^k_{e_j})^2}{\delta_{e_j}^{k}} \Big ),
\]
we estimate the second term on the right-hand side of \eqref{Est:main} as
\begin{equation}\label{zk:est1}
\begin{aligned}
& \int_0^t 4e^{G(t-s)} \sum_{j=1}^{|E|} \sum_{k \, : \, v_k \sim e_j}^{|V|} \frac{(\lambda^k_{e_j})^2}{\delta_{e_j}^{v_k}} |z_{k,0}|^2 e^{- 2\Lambda^k s} \, ds\\
& \quad \leq 4e^{Gt}H \sum_{k=1}^{|V|} |z_{k,0}|^2 \int_0^t e^{-(G+2 \Lambda^k)s} \, ds = 4 e^{Gt}H \sum_{k=1}^{|V|} |z_{k,0}|^2 \frac{(1-e^{-(G+2\Lambda^k)t})}{G + 2 \Lambda^k } \\
& \quad \leq 4 e^{Gt} \frac{H}{G} \sum_{k=1}^{|V|} |z_{k,0}|^2.
\end{aligned}
\end{equation}
To handle the third term on the right-hand side of \eqref{Est:main}, we first observe by the use of the inverse estimate \eqref{inv:est} and the Sobolev embedding theorem (with constant $C_S$)
\[
|W_{j}(r,v_k)| \leq C_{S} \| W_j(r) \|_{H^1(e_j)} \leq C_{S} C_N \| W_j(r) \|_{L^2(e_j)}.
\]
Then, setting 
\[
J := \max_{j=1, \dots, |E|} \Big ( |V| \max_{k \, :\, v_k \sim e_j} \frac{(\delta^k_{e_j})^2}{2 \Lambda^k} \Big ),
\]
and recalling that the summation of all vertex degree $\sum_{k=1}^{|V|} D_k$ is equal to twice the number of edges $2|E|$, we estimate the last term of \eqref{Est:main} as follows:
\begin{equation}\label{zk:est2}
\begin{aligned}
& \int_0^t 4e^{G(t-s)} \sum_{j=1}^{|E|} \sum_{k \, : \, v_k \sim e_j}^{|V|} \frac{(\lambda^k_{e_j})^2}{\delta_{e_j}^{v_k}} \Big ( D_k \sum_{j\, : \, e_j \sim v_k}^{|E|} \frac{(\delta^k_{e_j})^2}{2 \Lambda^k} \int_0^s | W_j (r,v_k)|^2 \, dr \Big )  \, ds \\
& \quad \leq 8 |E| e^{Gt} H J C_S C_N \int_0^t  e^{-Gs}  \int_0^s \sum_{j=1}^{|E|} \| W_j(r) \|_{L^2(e_j)}^2 \, dr \, ds.
\end{aligned}
\end{equation}
Let us define the constants 
\begin{align*}
Q_1 := \Big ( \sum_{i=1}^{M} \| u_{i,0} \|_{L^2(\Omega_i)}^2 + \sum_{j=1}^{|E|} \| w_{j,0} \|_{L^2(e_j)}^2 + 4 \frac{H}{G}\sum_{k=1}^{|V|} |z_{k,0}|^2 \Big ), \quad 
Q_2 := 8|E| H J C_S C_N,
\end{align*}
so that with the above estimates we infer from \eqref{Est:main} that 
\begin{equation}\label{Est:main:fp}
\begin{aligned}
& \Big ( \sum_{i=1}^{M} \| u_{i,N}(t) \|_{L^2(\Omega_i)}^2 + \sum_{j=1}^{|E|} \| w_{j,N}(t) \|_{L^2(e_j)}^2 \Big ) \\
& \quad \leq \Big ( Q_1 + Q_2  \int_0^t e^{-Gs} \int_0^s \sum_{j=1}^{|E|} \| W_j(r) \|_{L^2(e_j)}^2 \, dr \, ds \Big ) e^{Gt}, \quad \forall t \in [0,T_N].
\end{aligned}
\end{equation}
We introduce the set
\begin{align*}
\mathcal{S}_N & := \left \{ \bm{W} \in C^0 \Big ([0,T_N] ; \prod_{j=1}^{|E|} \mathbb{W}^{(j)}_N \Big ) \text{ such that } \right.  \\
& \qquad \qquad \left. \int_0^s \sum_{j=1}^{|E|} \| W_j(r) \|_{L^2(e_j)}^2 \, dr \leq A e^{Bs} \text{ for } s \in [0,T_N] \right \},
\end{align*}
for positive constants $A$ and $B$ to be determined.  We claim that for suitably chosen $A$ and $B$, it holds that the mapping $\mathcal{F}_N : \bm{W} \mapsto \bm{w}_N$ satisfies $\mathcal{F}_N : \mathcal{S}_N \to \mathcal{S}_N$. Integrating \eqref{Est:main:fp} from $0$ to $t$ yields
\begin{align*}
\int_0^t \sum_{j=1}^{|E|} \| w_{j,N}(s) \|_{L^2(e_j)}^2 \, ds & \leq \int_0^t \Big (Q_1 e^{Gs} + Q_2 e^{Gs} \int_0^s e^{-Gr} A e^{Br} \, dr \Big )\, ds \\
& = \frac{Q_1}{G} [e^{Gt} - 1] + \frac{Q_2 A}{B-G} \Big ( \frac{e^{Bt} - 1}{B} - \frac{e^{Gt}-1}{G} \Big ) \\
& = \frac{Q_2}{B(B-G)} A e^{Bt} + \frac{e^{Gt}}{G} \Big ( Q_1 - \frac{Q_2 A}{B-G} \Big ) + \frac{Q_2 A}{B-G} \Big ( \frac{1}{G} - \frac{1}{B} \Big ) - \frac{Q_1}{G}.
\end{align*}
Fixing $B > G$ such that $B(B-G) > Q_{2}$ and then choosing $A \in ( \frac{Q_1(B-G)}{Q_2}, \frac{Q_1 B}{Q_2})$ ensures that
\[
\frac{Q_2}{B(B-G)} < 1, \quad Q_1 - \frac{Q_2 A}{B-G} < 0, \quad \frac{Q_2 A}{B-G} \Big (\frac{1}{G} - \frac{1}{B} \Big ) - \frac{Q_1}{G} < 0,
\]
which in turn implies that
\[
\int_0^t \sum_{j=1}^{|E|} \| w_{j,N}(s) \|_{L^2(e_j)}^2 \leq A e^{Bt}, \quad \forall t \in [0,T_N].
\]
Next, we derive further estimates for the Galerkin solutions to show that $\overline{\mathcal{F}_N(\mathcal{S}_N)}$ is a compact subset of $C^0([0,T_N];\prod_{j=1}^{|E|} \mathbb{W}^{(j)}_N)$.  Multiplying \eqref{Gal:u} by $\frac{d}{dt}\mathrm{a}_r^{(i)}$ and \eqref{Gal:w} by $\frac{d}{dt}\mathrm{b}_s^{(j)}$ and summing over $1 \leq r, s \leq N$ yields
\begin{equation}\label{fixedpoint:compact1}
\begin{aligned}
& \sum_{i=1}^M \Big ( \| \pd_t u_{i,N} \|_{L^2(\Omega_i)}^2 + \frac{d}{dt} \int_{\Omega_i} \widehat{\kappa} (\nabla u_{i,N}) + \widehat{f}(u_{i,N}) \, dx\Big ) \\
& \qquad + \sum_{j=1}^{|E|} \Big ( \| \pd_t w_{j,N} \|_{L^2(e_j)}^2 + \frac{d}{dt} \int_{e_j} \widehat{\eta}(\pd_x w_{j,N}) + \widehat{g}(w_{j,N}) \, dx  \\
& \qquad  + \sum_{i \, : \, \Omega_i \sim e_j}^{M} \int_{e_j} \frac{\alpha_{ij}}{2} \frac{d}{dt}|w_{j,N}|^2 + \frac{\beta_{ij}}{2} \frac{d}{dt}|u_{i,N}|^2 - \alpha_{ij} w_{j,N} \pd_t u_{i,N} - \beta_{ij} u_{i,N} \pd_t w_{j,N} \, dx \Big ) \\
& \quad + \sum_{j=1}^{|E|} \sum_{k \, : \, v_k \sim e_j}^{|V|}\Big ( \big [\delta^{k}_{e_j} + \sum_{e_m \sim v_k, m \neq j} \gamma_{j \to m}^k \big ]  \frac{1}{2}  \frac{d}{dt} |w_{j,N}(v_k)|^2 \\
& \qquad \qquad \qquad - \sum_{e_m \sim v_k, m \neq j} \gamma_{m \to j}^k w_{m,N}(v_k)\pd_t w_{j,N}(v_k) \Big ) \\
& = \sum_{j=1}^{|E|}\sum_{k \, : \, v_k \sim e_j}^{|V|} \lambda^k_{e_j} z_{k,W} \pd_t w_{j,N}(v_k).
\end{aligned}
\end{equation}
For the indefinite terms on the third line of \eqref{fixedpoint:compact1} we apply Young's inequality, the inverse estimate \eqref{inv:est} and trace theorem (with constant $C_{\mathrm{tr}}$):
\begin{align*}
& \sum_{j=1}^{|E|} \sum_{i \, : \, \Omega_i \sim e_j}^{M} \left \vert \int_{e_j}  \alpha_{ij} w_{j,N} \pd_t u_{i,N} + \beta_{ij} u_{i,N} \pd_t w_{j,N} \, dx \right \vert  \\
& \quad \leq \sum_{j=1}^{|E|} \sum_{i \, : \, \Omega_i \sim e_j}^{M} \Big ( \alpha_{ij} C_\mathrm{tr} \| w_{j,N} \|_{L^2(e_j)} \| \pd_t u_{i,N} \|_{H^1(\Omega_i)} + \beta_{ij} \| u_{i,N} \|_{L^2(e_j)} \| \pd_t w_{j,N} \|_{L^2(e_j)} \Big ) \\
& \quad \leq \sum_{j=1}^{|E|}  \sum_{i \, : \, \Omega_i \sim e_j}^{M} \Big ( \alpha_{ij} C_N \| w_{j,N} \|_{L^2(e_j)} \| \pd_t u_{i,N} \|_{L^2(\Omega_i)} + \beta_{ij} \| u_{i,N} \|_{L^2(e_j)} \| \pd_t w_{j,N} \|_{L^2(e_j)} \Big ) \\
& \quad \leq \frac{1}{2} \sum_{i=1}^{M} \| \pd_t u_{i,N} \|_{L^2(\Omega_i)}^2 + \frac{1}{6} \sum_{j=1}^{|E|}  \| \pd_t w_{j,N} \|_{L^2(e_j)}^2 \\
& \qquad + C_N \sum_{j=1}^{|E|}  \sum_{i \, : \, \Omega_i \sim e_j}^{M}  \Big (\| w_{j,N} \|_{L^2(e_j)}^2 + \| u_{i,N} \|_{L^2(e_j)}^2 \Big ).
\end{align*}
Similarly, for the indefinite terms on the fourth line of \eqref{fixedpoint:compact1} we estimate as
\begin{align*}
& \sum_{j=1}^{|E|} \sum_{k \, : \, v_k \sim e_j}^{|V|} \sum_{e_m \sim v_k, m \neq j} \gamma_{m \to j}^k \vert w_{m,N}(v_k) \pd_t w_{j,N}(v_k) \vert \\
& \quad \leq \sum_{j=1}^{|E|} \sum_{k \, : \, v_k \sim e_j}^{|V|} \sum_{e_m \sim v_k, m \neq j} C_S  \| w_{m,N} \|_{H^1(e_m)} \| \pd_t w_{j,N} \|_{H^1(e_j)} \\
& \quad \leq \sum_{j=1}^{|E|} C_N \| w_{j,N} \|_{L^2(e_j)}^2 + \frac{1}{6} \sum_{j=1}^{|E|} \| \pd_t w_{j,N} \|_{L^2(e_j)}^2,
\end{align*}
and for the right-hand side of \eqref{fixedpoint:compact1} we have
\begin{align*}
& \sum_{j=1}^{|E|} \sum_{k \, : \, v_k \sim e_j}^{|V|} \lambda^k_{e_j} | z_{k,W} \pd_t w_{j,N}(v_k)| \leq \sum_{j=1}^{|E|} \sum_{k \, : \, v_k \sim e_j}^{|V|}  C_N |z_{k,W}| \| \pd_t w_{j,N} \|_{L^2(e_j)} \\
& \quad \leq \frac{1}{6} \sum_{j=1}^{|E|} \| \pd_t w_{j,N} \|_{L^2(e_j)}^2 + C_N \sum_{k=1}^{|V|}  |z_{k,W}|^2.
\end{align*}
Then, we infer from \eqref{fixedpoint:compact1} that 
\begin{equation}\label{fixedpoint:compact2}
    \begin{aligned}
 & \sum_{i=1}^M \Big ( \frac{1}{2} \| \pd_t u_{i,N} \|_{L^2(\Omega_i)}^2 + \frac{d}{dt} \int_{\Omega_i} \widehat{\kappa} (\nabla u_{i,N}) + F(u_{i,N}) \, dx\Big ) \\
& \qquad + \sum_{j=1}^{|E|} \Big ( \frac{1}{2} \| \pd_t w_{j,N} \|_{L^2(e_j)}^2 + \frac{d}{dt} \int_{e_j} \widehat{\eta}(\pd_x w_{j,N}) + F(w_{j,N}) \, dx \Big )  \\
& \qquad \qquad  + \sum_{i \, : \, \Omega_i \sim e_j}^{M} \frac{d}{dt} \Big ( \frac{\alpha_{ij}}{2} \| w_{j,N} \|_{L^2(e_j)}^2 + \frac{\beta_{ij}}{2} \| u_{i,N} \|_{L^2(e_j)}^2  \Big ) \\
& \qquad + \sum_{j=1}^{|E|} \sum_{k \, : \, v_k \sim e_j}^{|V|} \frac{d}{dt} \frac{1}{2} \Big ( \big [\delta^{k}_{e_j} + \sum_{e_m \sim v_k, m \neq j} \gamma_{j \to m}^k \big ]  |w_{j,N}(v_k)|^2 \Big ) \\
& \quad \leq C_N \sum_{j=1}^{|E|}\Big ( \| w_{j,N} \|_{L^2(e_j)}^2 + \sum_{i \, : \, \Omega_i \sim e_j}^{M} \| u_{i,N} \|_{L^2(e_j)}^2 \Big ) + C_N \sum_{k=1}^{|V|} |z_{k,W} |^2.
    \end{aligned}
\end{equation}
Using the assumptions \eqref{ass:aniso}, \eqref{ass:edgeaniso}, \eqref{ass:f}, \eqref{ass:g} and the Sobolev embedding theorem we infer that
\begin{align*}
& \int_{\Omega_i} \widehat{\kappa}(\nabla u_{i,N}(0)) + \widehat{f}(u_{i,N}(0)) \, dx \\
& \quad \leq C \| u_{i,0} \|_{W^{1,p}(\Omega_i)}^p + C \| u_{i,0} \|_{L^{q+1}(\Omega_i)}^{q+1} \leq C \| u_{i,0} \|_{W^{1,p}(\Omega_i)}^{\max(p,q+1)}, \\
& \int_{e_j} \widehat{\eta}(\pd_x w_{j,N}(0)) + \widehat{g}(w_{i,N}(0)) \, dx \\
& \quad \leq C \| \pd_{x} w_{j,0} \|_{W^{1,p}(e_j)}^p +  C \| w_{j,0} \|_{L^{r+1}(e_j)}^{r+1} \leq C \| w_{i,0} \|_{W^{1,p}(e_j)}^{\max(p,r+1)}.
\end{align*}
Hence, upon integrating \eqref{fixedpoint:compact2} over $(0,t)$, as well as invoking the estimates \eqref{Est:main}, \eqref{zk:est1} and \eqref{zk:est2} we obtain
\begin{align*}
& \Big (\sum_{i=1}^M \| \nabla u_{i,N}(t) \|_{L^p(\Omega_i)}^p + \sum_{j=1}^{|E|} \| \pd_x w_{j,N}(t) \|_{L^p(e_j)}^p \Big ) \\
& \qquad + \int_0^t \sum_{i=1}^{M} \| \pd_t u_{i,N} \|_{L^2(\Omega_i)}^2 + \sum_{j=1}^{|E|} \| \pd_{t} w_{j,N} \|_{L^2(e_j)}^2 \, ds \\
& \quad \leq C \sum_{i=1}^{M} \| u_{i,0} \|_{W^{1,p}(\Omega_i)}^{\max(p,q+1)} + C \sum_{j=1}^{|E|} \| w_{j,0} \|_{W^{1,p}(e_j)}^{\max(p,r+1)} \\
& \qquad + C_N \int_0^{t} \sum_{j=1}^{|E|}\Big ( \| w_{j,N} \|_{L^2(e_j)}^2 + \sum_{i \, : \, \Omega_i \sim e_j}^{M} \| u_{i,N} \|_{L^2(e_j)}^2 \Big ) \, ds + C_N \int_0^t \sum_{k=1}^{|V|} |z_{k,W} |^2 \, ds \\
& \quad \leq C_N \Big ( 1 + \int_0^t \sum_{j=1}^{|E|} \| W_j(s) \|_{L^2(e_j)}^2 \, ds \Big ) e^{Gt} \leq C_N \Big ( 1 + A e^{BT_N} \Big ) e^{G T_N} =: K_N.
\end{align*}
Hence, we see that
\[
\mathcal{F}_N(\mathcal{S}_N) = \left \{ \bm{W} \in \mathcal{S}_N \, : \, \int_0^{T_N} \sum_{j=1}^{|E|} \| \pd_t W_j \|_{L^2(e_j)}^2 \, ds \leq K_N \right \},
\]
and as $\mathbb{W}^{(j)}_N$ is finite dimensional for $j = 1, \dots, |E|$, by the Aubin--Lions lemma we have the compact embedding
\[
C^0([0,T_N]; \mathbb{W}^{(j)}_N) \cap H^1(0,T_N; \mathbb{W}^{(j)}_N) \Subset C^0([0,T_N];\mathbb{W}^{(j)}_N),
\]
which yields that $\overline{\mathcal{F}_N(\mathcal{S}_N)}$ is compact in $C^0([0,T_N]; \prod_{j=1}^{|E|} \mathbb{W}^{(j)}_N)$.

It remains to establish the continuity of $\mathcal{F}_N$ in order to invoke Schauder's fixed point theorem. We consider a sequence $(\bm{W}_n)_{n \in \N} \subset \mathcal{S}_N$ such that 
\[
\bm{W}_n \to \bm{W} \text{ strongly in } C^0\Big ([0,T];\prod_{j=1}^{|E|} \mathbb{W}^{(j)}_N\Big).
\]
Let $\bm{w}_{n} = \mathcal{F}_N(\bm{W}_n)$ and $\bm{w} = \mathcal{F}_N(\bm{W})$, and we denote by $(\bu_{n}, \bm{z}_{n})$ and $(\bu, \bm{z})$ as the associated subdomain and vertex variables satisfying  the analogue of \eqref{Gal:ODE} and \eqref{aux:prob:Galerkin} corresponding to $\bm{W}_n$ and $\bm{W}$, respectively, furnished with the same initial data. Then, the differences
\[
\overline{u}_i = u_{n,i} - u_i, \quad \overline{w}_j = w_{n,j} - w_j, \quad \overline{z}_k = z_{n,k} - z_k, \quad \overline{W}_j = W_{n,j} - W_j
\]
for $1 \leq i \leq M$, $1 \leq j \leq |E|$ and $1 \leq k \leq |V|$ satisfy for all $1 \leq r, s \leq N$
\begin{subequations}\label{aux:prob:Galerkin:continuity}
\begin{alignat}{2}
\label{Gal:cont:u} 0 & = \sum_{i=1}^M \int_{\Omega_i} \pd_t \overline{u}_{i} \phi_r^{(i)} + (\kappa(\nabla u_{n,i}) - \kappa(\nabla u_i)) \cdot \nabla \phi_r^{(i)} + (f(u_{n,i}) - f(u_i)) \phi_r^{(i)} \, dx \\
\notag & \qquad \qquad - \sum_{j \, : \, e_j \sim \Omega_i}^{|E|}\int_{e_j} (\alpha_{ij} \overline{w}_{j} - \beta_{ij} \overline{u}_{i}) \phi_r^{(i)} \, dx \\
\label{Gal:cont:w} 0 & = \sum_{j=1}^{|E|} \Big ( \int_{e_j} \pd_t \overline{w}_{j} \psi_{s}^{(j)} + (\eta(\pd_x w_{n,j}) - \eta(\pd_x w_j)) \pd_x \psi_{s}^{(j)} + (g(w_{n,j}) - g(w_j)) \psi_{s}^{(j)}  \, dx \\
\notag & \qquad  \qquad + \sum_{i \, : \, \Omega_i \sim e_j}^{M} \int_{e_j} (\alpha_{ij} \overline{w}_{j} - \beta_{ij} \overline{u}_{i})\psi_{s}^{(j)} \, dx \Big ) \\
\notag & \quad + \sum_{j=1}^{|E|} \sum_{k \, : \, v_k \sim e_j}^{|V|} \Big (\big [\delta^{k}_{e_j} + \sum_{e_m \sim v_k, m \neq j} \gamma^k_{j \to m} \big ]  \overline{w}_{j} - \lambda^k_{e_j} \overline{z}_{k} - \sum_{e_m \sim v_k, m \neq j} \gamma^k_{m \to j} \overline{w}_{m} \Big ) \psi_s^{(j)} \Big \vert_{v_k},
\end{alignat}
\end{subequations}
and
\begin{align}\label{Gal:ODE:continuity}
\frac{d}{dt} \overline{z}_{k} = - \Lambda^k \overline{z}_{k} + \sum_{j \, : \, e_j \sim v_k}^{|E|} \delta^k_{e_j} \overline{W}_j(t, v_k).
\end{align}
Multiplying \eqref{Gal:cont:u} by the coefficients of $u_{n,i} - u_i$ and \eqref{Gal:cont:w} by the coefficients of $w_{n,j} - w_j$, we obtain upon summing over $1 \leq r, s \leq N$
\begin{align*}
& \sum_{i=1}^M \Big (\frac{d}{dt} \frac{1}{2} \| \overline{u}_{i} \|^2_{L^2(\Omega_i)} + \int_{\Omega_i} c_3 |\nabla \overline{u}_{i}|^p \, dx \Big ) \\
& \qquad + \sum_{j=1}^{|E|} \Big (\frac{d}{dt} \frac{1}{2} \| \overline{w}_{j} \|^2_{L^2(e_j)} + \int_{e_j} c_7 |\pd_x \overline{w}_j|^p \, dx \\
& \qquad \qquad + \sum_{i \, : \, \Omega_i \sim e_j}^{M} \int_{e_j} \alpha_{ij} |\overline{w}_{j}|^2 + \beta_{ij} |\overline{u}_{i}|^2 - (\alpha_{ij} + \beta_{ij}) \overline{w}_{j} \overline{u}_{i}\, dx \Big ) \\
& \qquad + \sum_{j=1}^{|E|} \sum_{k \, : \, v_k \sim e_j}^{|V|} \Big ( \big [\delta^{k}_{e_j} + \sum_{e_m \sim v_k, m \neq j} \gamma_{j \to m}^k \big ]  |\overline{w}_{j}(v_k)|^2 \\
& \qquad \qquad \qquad - \sum_{e_m \sim v_k, m \neq j} \gamma_{m \to j}^k \overline{w}_{m}(v_k) \overline{w}_{j}(v_k) \Big ) \\
& \quad \leq \sum_{j=1}^{|E|}\sum_{k \, : \, v_k \sim e_j}^{|V|} \lambda^k_{e_j} \overline{z}_{k} \overline{w}_{j}(v_k),
\end{align*}
where we have used the coercivity of $\kappa$ and $\eta$, as well as the monotonicity of $f$ and $g$. Adopting a similar argument previously used to show $\mathcal{F}_N : \mathcal{S}_N \to \mathcal{S}_N$, we obtain the analogue of \eqref{zkW:est}:
\begin{equation}\label{diff:z:est}
\begin{aligned}
|\overline{z}_k(t)|^2 & \leq 2D_k \sum_{j \, : \, e_j \sim v_k}^{|E|} \frac{(\delta^{k}_{e_j})^2}{2 \Lambda^k} \int_0^t |\overline{W}_j(s, v_k)|^2 \, ds  \\
& \leq 2 D_k J C_S C_N \int_0^t \sum_{j=1}^{|E|} \| \overline{W}_j (s) \|_{L^2(e_j)}^2 \, ds,
\end{aligned}
\end{equation}
and the analogue of \eqref{diff:ineq}:
\begin{equation}\label{diff:uw:est}
    \begin{aligned}
& \sum_{i=1}^M \Big (\frac{d}{dt} \frac{1}{2} \| \overline{u}_{i} \|^2_{L^2(\Omega_i)} + c_3 \|\nabla \overline{u}_{i}\|_{L^p(\Omega_i)}^{p} \Big ) \\
& \qquad + \sum_{j=1}^{|E|} \Big (\frac{d}{dt} \frac{1}{2} \| \overline{w}_{j} \|^2_{L^2(e_j)} + c_7 \| \pd_x \overline{w}_{j} \|_{L^p(e_j)}^p \\
& \qquad \qquad + \sum_{i \, : \, \Omega_i \sim e_j}^{M} \frac{\beta_{ij}}{2} \| \overline{u}_{i} \|_{L^2(e_j)}^2 + \frac{1}{4}\sum_{k \, : \, v_k \sim e_j}^{|V|} \delta^{k}_{e_j} |\overline{w}_{j,N}(v_k)|^2 \Big ) \\
& \quad \leq  \sum_{j=1}^{|E|} \sum_{k \, : \, v_k \sim e_j}^{|V|} \frac{(\lambda^k_{e_j})^2}{\delta_{e_j}^{v_k}} |\overline{z}_{k}|^2  + \sum_{j=1}^{|E|} \frac{(\alpha_i + \beta_j)^2}{2 \beta_j}\| \overline{w}_{j} \|_{L^2(e_j)}^2 \\
& \quad \leq \frac{Q_2}{2} \int_0^t \sum_{j=1}^{|E|} \| \overline{W}_j(s) \|_{L^2(e_j)}^2 \, ds +  \frac{G}{2} \sum_{j=1}^{|E|} \| \overline{w}_{j} \|_{L^2(e_j)}^2,
    \end{aligned}
\end{equation}
where we used \eqref{diff:z:est} and the definitions of the constants $G$, $H$ and $J$. Applying Gronwall's inequality to \eqref{diff:uw:est} we obtain
\begin{align*}
\max_{t \in [0,T_N]} \Big ( \sum_{i=1}^{M} \| \overline{u}_i(t) \|_{L^2(\Omega_i)}^2 + \sum_{j=1}^{|E|} \| \overline{w}_j(t) \|_{L^2(e_j)}^2 \Big ) \leq C(T_N) \max_{t \in [0,T_N]} \sum_{j=1}^{|E|} \| \overline{W}_j(t) \|_{L^2(e_j)}^2.
\end{align*}
Since $\overline{W}_j \to 0$ strongly in $C^0([0,T_N]; \mathbb{W}^{(j)}_N)$ we deduce that 
\begin{align*}
u_{n,i} \to u_i & \text{ strongly in } C^0([0,T_N]; \mathbb{U}^{(i)}_N), \\
w_{n,j} \to w_j & \text{ strongly in } C^0([0,T_N];\mathbb{W}^{(j)}_N),
\end{align*}
for $1 \leq i \leq M$ and $1 \leq j \leq |E|$, leading to the continuity of $\mathcal{F}_N$. Hence, for each $n \in \N$, we infer via Schauder's fixed point theorem the existence of solutions 
\[
\hat{u}_{i,N} \in C^0([0,T_N]; \mathbb{U}^{(i)}_N), \quad \hat{w}_{j,N} \in C^0([0,T_N]; \mathbb{W}^{(j)}_N), \quad \hat{z}_{k,N} \in C^1([0,T_N]),
\]
for $1 \leq i \leq M$, $1 \leq j \leq |E|$ and $1 \leq k \leq |V|$ to the approximate system \eqref{fixedpoint:sys}.

\subsection{Uniform estimates}
In this section we drop the hat notation and derive uniform estimates in $N$ for the approximate solutions $(\bm{u}_N, \bm{w}_N, \bm{z}_N)$ to \eqref{fixedpoint:sys}. The symbol $C$ denotes positive constants independent of $N$. Testing \eqref{fp:u} with $u_{i,N}$, \eqref{fp:w} with $w_{j,N}$ and \eqref{fp:z} with $z_{k,N}$, summing over $1 \leq i \leq M$, $1 \leq j \leq |E|$, $1 \leq k \leq |V|$, and employing the coercivity of $\kappa$ and $\eta$, as well as the monotonicity of $f$ and $g$ yields
\begin{equation}\label{unifest:1}
\begin{aligned}
&  \sum_{k=1}^{|V|} \frac{1}{2} \frac{d}{dt} |z_{k,N}|^2  + \sum_{i=1}^M \Big (\frac{d}{dt} \frac{1}{2} \| u_{i,N} \|^2_{L^2(\Omega_i)} + c_3 \|\nabla u_{i,N} \|_{L^p(\Omega_i)}^p \Big ) \\
& \qquad + \sum_{j=1}^{|E|} \Big (\frac{d}{dt} \frac{1}{2} \| w_{j,N} \|^2_{L^2(e_j)} + c_7 \| \pd_x w_{j,N} \|_{L^p(e_j)}^p \\
& \qquad \qquad + \sum_{i \, : \, \Omega_i \sim e_j}^{M} \int_{e_j} \alpha_{ij} |w_{j,N}|^2 + \beta_{ij} |u_{i,N}|^2 - (\alpha_{ij} + \beta_{ij}) w_{j,N} u_{i,N}\, dx \Big ) \\
& \qquad  + \sum_{k=1}^{|V|} \sum_{j \neq m \, : \, e_j, e_m \sim v_k}^{|E|} \Big [ (\gamma^k_{j \to m} w_{j,N}(v_k) - \gamma^k_{m \to j} w_{m,N}(v_k) )(w_{j,N}(v_k) - w_{m,N}(v_k)) \Big ] \\
& \qquad + \sum_{j=1}^{|E|} \sum_{k \, : \, v_k \sim e_j}^{|V|} (\delta^{k}_{e_j} w_{j,N}(v_k)  - \lambda^k_{e_j} z_{k,N})(w_{j,N}(v_k) - z_{k,N}) \\ 
& \quad \leq 0.
\end{aligned}
\end{equation}
The third line of \eqref{unifest:1} can be handled as before:
\begin{equation}\label{rhs:3}
\begin{aligned}
& \sum_{j=1}^{|E|}\sum_{i \, : \, \Omega_i \sim e_j}^{M} \int_{e_j} \alpha_{ij} |w_{j,N}|^2 + \beta_{ij} |u_{i,N}|^2 - (\alpha_{ij} + \beta_{ij}) w_{j,N} u_{i,N}\, dx  \\
& \quad \geq - \frac{1}{2}G \sum_{j=1}^{|E|} \| w_{j,N} \|_{L^2(e_j)}^2 + \sum_{j=1}^{|E|} \sum_{i \, : \, \Omega_i \sim e_j}^{M} \frac{\beta_{ij}}{2} \| u_{i,N} \|_{L^2(e_j)}^2.
\end{aligned}
\end{equation}
For the fourth line of \eqref{unifest:1} we recall \eqref{delta:gamma:term} and \eqref{ass:coeff} to see that
\begin{equation}\label{rhs:4}
\begin{aligned}
&  \sum_{k=1}^{|V|} \sum_{j \neq m \, : \, e_j, e_m \sim v_k}^{|E|} \Big [ (\gamma^k_{j \to m} w_{j,N} - \gamma^k_{m \to j} w_{m,N} )(w_{j,N} - w_{m,N}) \Big ] \Big \vert_{v_k} \\
& \quad = \sum_{j=1}^{|E|} \sum_{k \, : \, v_k \sim e_j}^{|V|} \Big ( \Big [ \sum_{e_m \sim v_k, m \neq j} \gamma_{j \to m}^k \Big ] |w_{j,N}(v_k)|^2 - \sum_{e_m \sim v_k, m \neq j} \gamma_{m \to j}^k w_{m,N}(v_k) w_{j,N}(v_k) \Big ) \\
& \quad \geq \sum_{j=1}^{|E|} \sum_{k\, : \, v_k \sim e_j} \frac{1}{2} \Big [\sum_{e_m \sim v_k, m \neq j} (\gamma_{j \to m}^k - \gamma_{m \to j}^k) \Big ] |w_{j,N}(v_k)|^2 \\
& \quad \geq - \frac{1}{2} \sum_{j=1}^{|E|} \sum_{k \, : \, v_k \sim e_j}^{|V|} \delta^{k}_{e_j} |w_{j,N}(v_k)|^2.
\end{aligned}
\end{equation}
Lastly, we treat the fifth line of \eqref{unifest:1} in a similar fashion to the third line of \eqref{unifest:1}:
\begin{equation}\label{rhs:5}
\begin{aligned}
& \sum_{j=1}^{|E|} \sum_{k \, : \, v_k \sim e_j}^{|V|} (\delta^{k}_{e_j} w_{j,N}(v_k)  - \lambda^k_{e_j} z_{k,N})(w_{j,N}(v_k) - z_{k,N}) \\
& \quad =  \sum_{j=1}^{|E|} \sum_{k \, : \, v_k \sim e_j}^{|V|} \delta^{k}_{e_j} |w_{j,N}(v_k)|^2 + \lambda^k_{e_j} |z_{k,N}|^2 - (\delta^{k}_{e_j} + \lambda^k_{e_j})w_{j,N}(v_k) z_{k,N} \\
& \quad \geq  \sum_{j=1}^{|E|} \sum_{k \, : \, v_k \sim e_j}^{|V|} \frac{3}{4} \delta^{k}_{e_j} |w_{j,N}(v_k)|^2 + \Big ( \lambda^k_{e_j} - \frac{(\delta^{k}_{e_j} + \lambda^k_{e_j})^2}{ \delta^{k}_{e_j}} \Big ) |z_{k,N}|^2.
\end{aligned}
\end{equation}
We then deduce from \eqref{unifest:1} that
\begin{equation}\label{unifest:1:gron}
\begin{aligned}
&  \sum_{k=1}^{|V|} \frac{1}{2} \frac{d}{dt} |z_{k,N}|^2  + \sum_{i=1}^M \Big (\frac{d}{dt} \frac{1}{2} \| u_{i,N} \|^2_{L^2(\Omega_i)} + c_3 \|\nabla u_{i,N} \|_{L^p(\Omega_i)}^p \Big ) \\
& \qquad + \sum_{j=1}^{|E|} \Big (\frac{d}{dt} \frac{1}{2} \| w_{j,N} \|^2_{L^2(e_j)} + c_7 \| \pd_x w_{j,N} \|_{L^p(e_j)}^p \Big ) \\
& \qquad + \sum_{j=1}^{|E|} \Big ( \sum_{i \, : \, \Omega_i \sim e_j}^{M} \frac{\beta_{ij}}{2} \| u_{i,N} \|_{L^2(e_j)}^2 + \sum_{k \, : \, v_k \sim e_j}^{|V|} \frac{\delta^{k}_{e_j}}{4} |w_{j,N}(v_k)|^2 \Big ) \\
& \quad \leq C \sum_{j=1}^{|E|} \| w_{j,N} \|_{L^2(e_j)}^2 + C \sum_{k=1}^{|V|} |z_{k,N}|^2
\end{aligned}
\end{equation}
with positive constants $C$ independent of $N$. By Gronwall's inequality and recalling \eqref{ini:bdd}, we find that
\begin{equation}\label{unif:1}
\begin{aligned}
& \sup_{t \in (0,T)} \Big ( \sum_{i=1}^{M} \| u_{i,N}(t) \|_{L^2(\Omega_i)}^2 + \sum_{j=1}^{|E|} \| w_{j,N}(t) \|_{L^2(e_j)}^2 + \sum_{k=1}^{|V|} |z_{k,N}(t)|^2 \Big ) \\
& \quad + \int_0^T  \Big ( \sum_{i=1}^{M} \| \nabla u_{i,N} \|_{L^p(\Omega_i)}^p + \sum_{j=1}^{|E|} \| \pd_x w_{j,N} \|_{L^p(e_j)}^p \Big ) \, dt \leq C(\bm{u}_0, \bm{w}_0, \bm{z}_0, T).
\end{aligned}
\end{equation}
In particular, this uniform estimate allows us to infer that the semi-Galerkin solutions $(\bm{u}_N, \bm{w}_N)$ can be extended from $[0,T_N]$ to a uniform time interval $[0,T]$ for a fixed but arbitrary $T> 0$. Next, for arbitrary test function $\zeta_i \in W^{1,p}(\Omega_i)$ we write 
\[
\zeta_i = P_{\mathbb{U}^{(j)}_N} \zeta_i + (\mathbb{I} - P_{\mathbb{U}^{(j)}_N}) \zeta_i,
\]
and consider testing \eqref{fp:u} by $\zeta_i$, leading to 
\begin{align*}
\int_{\Omega_i} \pd_t u_{i,N} P_{\mathbb{U}^{(i)}_N} \zeta_i \, dx & = \int_{\Omega_i} - \kappa (\nabla u_{i,N}) \cdot \nabla P_{\mathbb{U}^{(i)}_N} \zeta_i - f(u_{i,N})  P_{\mathbb{U}^{(i)}_N} \zeta_i \, dx \\
& \quad + \sum_{j \, : \, e_j \sim \Omega_i}^{|E|} \int_{e_j} (\alpha_j w_{j,N} - \beta_i u_{i,N})  P_{\mathbb{U}^{(i)}_N} \zeta_i \, dx.
\end{align*}
From \eqref{ass:aniso} we see that 
\begin{align*}
& \left | \int_{\Omega_i} \kappa( \nabla u_{i,N}) \cdot \nabla  P_{\mathbb{U}^{(i)}_N} \zeta_i \, dx \right | \leq c_2 \int_{\Omega_i} | \nabla u_{i,N} |^{p-1} | \nabla P_{\mathbb{U}^{(i)}_N} \zeta_i| \, dx \\
& \quad \leq C \| \nabla u_{i,N} \|_{L^p(\Omega_i)}^{p-1} \| \nabla P_{\mathbb{U}^{(i)}_N} \zeta_i \|_{L^{p}(\Omega_i)} \leq C \| \nabla u_{i,N} \|_{L^p(\Omega_i)}^{p-1} \| \zeta_i \|_{W^{1,p}(\Omega_i)}.
\end{align*}
Meanwhile, by the trace theorem we have
\begin{align*}
\left |\int_{e_j} (\alpha_{ij} w_{j,N} - \beta_{ij} u_{i,N})  P_{\mathbb{U}^{(i)}_N} \zeta_i \, dx  \right | \leq  C \Big (\alpha_{ij} \|w_{j,N} \|_{L^2(e_j)} + \beta_{ij} \|u_{i,N} \|_{L^2(e_j)} \Big ) \| \zeta_i \|_{W^{1,p}(\Omega_i)},
\end{align*}
while by \eqref{ass:f} and H\"older's inequality, for any $z \in (1,\infty)$ if $p = 2$ and $z  = 1$ if $p > 2$ (due to the Sobolev embedding $W^{1,p}(\Omega_i) \subset L^\infty(\Omega_i)$), we have
\begin{align*}
& \left | \int_{\Omega_i} f(u_{i,N})  P_{\mathbb{U}^{(i)}_N} \zeta_i \, dx\right | \leq c_9 \int_{\Omega_i} ( 1 + |u_{i,N}|^{q}) | P_{\mathbb{U}^{(i)}_N} \zeta_i| \, dx \\
& \quad \leq C \| P_{\mathbb{U}^{(i)}_N} \zeta_i \|_{L^1(\Omega_i)} + C\| u_{i,N} \|_{L^{zq}(\Omega_i)}^{q} \|P_{\mathbb{U}^{(i)}_N} \zeta_i \|_{L^{\frac{z}{z-1}}(\Omega_i)} \\
& \quad \leq C\Big ( 1 + \| u_{i,N} \|_{L^{zq}(\Omega_i)}^q \Big ) \| \zeta_i \|_{W^{1,p}(\Omega_i)}.
\end{align*}
Collecting the above estimates yields
\begin{equation}\label{pdt:u:est:pre}
\begin{aligned}
& \int_0^T \| \pd_t u_{i,N} \|_{W^{1,p}(\Omega_i)^*}^{\frac{p}{p-1}} \, dt \\
& \quad \leq C \int_0^T \Big ( 1 + \| \nabla u_{i,N} \|_{L^p(\Omega_i)}^p + \| u_{i,N} \|_{L^{zq}(\Omega_i)}^{\frac{p}{p-1} q} + \| u_{i,N} \|_{L^2(e_j)}^2 \Big ) \, dt,
\end{aligned}
\end{equation}
where we used the uniform boundedness of $w_{j,N}$ in $L^\infty(0,T;L^2(e_j))$. From the following Gagliardo--Nirenberg inequality in two dimensions:
\begin{align*}
\| u \|_{L^s(\Omega_i)} \leq C \| u \|_{L^2(\Omega_i)}^{1 -(\frac{p}{p-1} \frac{s-2}{2s}) } \| u \|_{W^{1,p}(\Omega_i)}^{\frac{p}{p-1} \frac{s-2}{2s}}
\end{align*}
whilst recalling the uniform boundedness of $u_{i,N}$ in $L^\infty(0,T;L^2(\Omega_i))$ we see that with $s = zq$,
\begin{align}\label{GN:2}
\| u_{i,N} \|_{L^{zq}(\Omega_i)}^q \leq C \| u_{i,N} \|_{W^{1,p}(\Omega_i)}^{\frac{p}{p-1} \frac{zq-2}{2z}}.
\end{align}
Hence, for $p = 2$, choosing $z = 2$ if $q \in [0,2]$ and $z = \frac{2}{q-1}$ if $q \in (2, 2+\frac{2(p-1)^2}{p})$, we deduce that
\[
\| u_{i,N} \|_{L^{zq}(\Omega_i)}^{2q} \leq C \| u_{i,N} \|_{H^1(\Omega_i)}^2.
\]
Meanwhile, for $p \in (2,\infty)$, choosing $z = 1$ leads to
\[
\| u_{i,N} \|_{L^q(\Omega_i)}^{\frac{p}{p-1} q} \leq C \| u_{i,N} \|_{W^{1,p}(\Omega_i)}^{p},
\]
and so from \eqref{pdt:u:est:pre} we infer that
\begin{align}\label{pdt:u:est}
\int_0^T \| \pd_t u_{i,N} \|_{W^{1,p}(\Omega_i)^*}^{\frac{p}{p-1}}\, dt \leq C \int_0^T \Big (1 + \| u_{i,N} \|_{W^{1,p}(\Omega_i)}^{p} +  \| u_{i,N} \|_{L^2(e_j)}^{2} \Big ) \, dt \leq C.
\end{align}
%


\noindent 
Similarly, for arbitrary test function $\chi_j \in W^{1,p}(e_j)$ we write
\[
\chi_j = P_{\mathbb{W}^{(j)}_N} \chi_j + (\mathbb{I} - P_{\mathbb{W}^{(j)}_N}) \chi_j,
\]
and consider testing \eqref{fp:w} by $\chi_j$, leading to 
\begin{align*}
& \int_{e_j} \pd_t w_{j,N} P_{\mathbb{W}^{(j)}_N} \chi_j \, dx \\
& \quad =  \int_{e_j} -\eta(\pd_x w_{j,N}) \pd_x P_{\mathbb{W}^{(j)}_N} \chi_j - \Big ( g(w_{j,N}) + \sum_{m \, : \, \Omega_m \sim e_j}^{M}  (\alpha_{mj} w_{j,N} - \beta_{mj} u_{i,N}) \Big ) P_{\mathbb{W}^{(j)}_N} \chi_j\, dx \\
& \qquad - \sum_{k \, : \, v_k \sim e_j}^{|V|} \Big ( \Big [\delta^{k}_{e_j} + \sum_{e_m \sim v_k, m \neq j} \gamma_{j \to m}^k \Big ] w_{j,N} - \lambda^k_{e_j} z_{k,N} - \sum_{e_m \sim v_k, m \neq j} \gamma_{m \to j}^k w_{m,N} \Big ) P_{\mathbb{W}^{(j)}_N} \chi_j \Big \vert_{v_k}.
\end{align*}
Using \eqref{ass:edgeaniso} and the Sobolev embedding $W^{1,p}(e_j) \subset C^0(e_j)$ for $p \geq 2$, we find that
\begin{align*}
& \left |  \int_{e_j} \eta(\pd_x w_{j,N}) \pd_x P_{\mathbb{W}^{(j)}_N} \chi_j \right | \leq \int_{e_j} c_6 |\pd_x w_{j,N} |^{p-1} |\pd_x P_{\mathbb{W}^{(j)}_N} \chi_j| \, dx \\
& \quad \leq C \| \pd_x w_{j,N} \|_{L^p(e_j)}^{p-1} \| \chi_j \|_{W^{1,p}(e_j)}, \\
& \left | \int_{e_j} (\alpha_{mj} w_{j,N} - \beta_{mj} u_{m,N})  P_{\mathbb{W}^{(j)}_N} \chi_j \, dx \right | \leq C \Big ( \| w_{j,N} \|_{L^2(e_j)} + \| u_{m,N} \|_{L^2(e_j)} \Big ) \| \chi_i \|_{W^{1,p}(e_j)}, \\
& \left | w_{j,N}(v_k) P_{\mathbb{W}^{(j)}_N} \chi_j (v_k) \right | \leq C \| w_{j,N} \|_{W^{1,p}(e_j)} \| \chi_j \|_{W^{1,p}(e_j)}, \\
&  \left | \int_{e_j} g(w_{j,N})  P_{\mathbb{W}^{(j)}_N} \chi_j \, dx \right | \leq C \Big ( 1 + \| w_{j,N} \|_{L^r(e_j)}^r \Big ) \| \chi_j \|_{W^{1,p}(e_j)} .
\end{align*}
Collecting the above estimates yields
\begin{equation}\label{pdt:w:est:pre}
\begin{aligned}
& \int_0^T \| \pd_t w_{j,N} \|_{W^{1,p}(e_j)^*}^{\frac{p}{p-1}} \, dt \\
& \quad \leq C \int_0^T \Big ( 1 + \| \pd_x w_{j,N} \|_{L^p(e_j)}^{p} + \| w_{j,N} \|_{L^r(e_j)}^{\frac{p}{p-1}r} + \sum_{i=1}^{M}\| u_{i,N} \|_{L^2(e_j)}^2 \Big ) \, dt.
\end{aligned}
\end{equation}
From the following Gagliardo--Nirenberg inequality in one dimension:
\[
\| w \|_{L^r(e_j)} \leq C \| w \|_{L^2(e_j)}^{1-(\frac{2p}{3p-2} \frac{r-2}{2r})} \| w \|_{W^{1,p}(e_j)}^{\frac{2p}{3p-2} \frac{r-2}{2r}},
\]
whilst recalling the uniform boundedness of $w_{j,N}$ in $L^\infty(0,T;L^2(e_j))$ we see that 
\begin{align*}
\| w_{j,N} \|_{L^r(e_j)}^{\frac{p}{p-1}r} \leq C \| w_{j,N} \|_{W^{1,p}(e_j)}^{\frac{p}{p-1}\frac{2p}{3p-2} \frac{r-2}{2}} \leq C \| w_{j,N} \|_{W^{1,p}(e_j)}^p.
\end{align*}
Hence, from \eqref{pdt:w:est:pre} we find that 
\begin{align}\label{pdt:w:est}
\int_0^T \| \pd_t w_{j,N} \|_{W^{1,p}(e_j)^*}^{\frac{p}{p-1}} \, dt \leq C \int_0^T  \Big (1 + \| w_{j,N} \|_{W^{1,p}(e_j)}^p + \sum_{i=1}^{M} \| u_{i,N} \|_{L^2(e_j)}^2 \Big ) \, dt \leq C.
\end{align}
Lastly, we multiple \eqref{fp:z} by $\frac{d}{dt} z_{k,N}$ to see that 
\begin{align*}
 \Big | \frac{d}{dt} z_{k,N} \Big |^2 + \frac{\Lambda^k}{2} \frac{d}{dt} |z_{k,N}|^2&  = \Big (\sum_{j \, : \, e_j \sim v_k}^{|E|} \delta^k_{e_j} w_{j,N}(t, v_k) \Big ) \frac{d}{dt} z_{k,N} \\
&  \leq \frac{1}{2} \Big | \frac{d}{dt} z_{k,N} \Big |^2 + C\sum_{j \, : \, e_j \sim v_k}^{|E|} |w_{j,N}(t, v_k)|^2.
\end{align*}
Using the Sobolev embedding $W^{1,p}(e_j) \subset C^0(e_j)$ and the fact that $p \geq 2$ we see that
\begin{align*}
  \frac{1}{2} \Big | \frac{d}{dt} z_{k,N} \Big |^2 + \frac{\Lambda^k}{2} \frac{d}{dt} |z_{k,N}|^2 \leq C \sum_{j \, : \, e_j \sim v_k}^{|E|} \| w_{j,N} \|_{W^{1,p}(e_j)}^p.
\end{align*}
Integrating over $(0,T)$ allows us to infer 
\begin{align}\label{z':unif}
\int_0^T \Big | \frac{d}{dt} z_{k,N} \Big |^2 \, dt + \max_{t \in (0,T)} |z_{k,N}(t)|^2 \leq C \Big ( 1 + \int_0^T  \sum_{j=1}^{|E|} \| w_{j,N} \|_{W^{1,p}(e_j)}^p \, dt \Big ) \leq C.
\end{align}

\subsection{Compactness and passing to the limit}
From \eqref{unif:1}, \eqref{pdt:u:est},  \eqref{pdt:w:est} and \eqref{z':unif} we find that 
\begin{align*}
\{ u_{i,N} \}_{N \in \N} & \text{ is bounded in } L^\infty(0,T;L^2(\Omega_i)) \cap L^p(0,T;W^{1,p}(\Omega_i)), \\
\{ \pd_t u_{i,N} \}_{N \in \N} & \text{ is bounded in } L^{p/(p-1)}(0,T;W^{1,p}(\Omega_i)^*), \\
\{ w_{j,N} \}_{N \in \N} & \text{ is bounded in } L^\infty(0,T;L^2(e_j)) \cap L^p(0,T;W^{1,p}(e_j)), \\
\{ \pd_t w_{j,N} \}_{N \in \N} & \text{ is bounded in } L^{p/(p-1)}(0,T;W^{1,p}(e_j)^*), \\
\{ z_{k,N} \}_{N \in \N} & \text{ is bounded in } L^\infty(0,T) \cap H^1(0,T).
\end{align*}
Then, by the Banach--Alaoglu theorem, there exist limit functions $u_i$, $w_j$, $z_k$ for $1 \leq i \leq M$, $1 \leq j \leq |E|$ and $1 \leq k \leq |V|$, such that along a non-relabelled subsequence $N \to \infty$,
\begin{align*}
u_{i,N} \to u_i & \text{ weakly* in } L^\infty(0,T;L^2(\Omega_i)) \cap L^p(0,T;W^{1,p}(\Omega_i)), \\
\pd_t u_{i,N} \to \pd_t u_i & \text{ weakly in } L^{p/(p-1)}(0,T;W^{1,p}(\Omega_i)^*), \\
w_{j,N} \to w_j & \text{ weakly* in } L^\infty(0,T;L^2(e_j)) \cap L^p(0,T;W^{1,p}(e_j)), \\
\pd_t w_{j,N} \to \pd_t u_i & \text{ weakly in } L^{p/(p-1)}(0,T;W^{1,p}(e_j)^*), \\
z_{k,N} \to z_k & \text{ weakly in } H^1(0,T).
\end{align*}
Using \cite[Section 8, Corollary 4]{Simon}, we also obtain the following strong convergences
\begin{align*}
u_{i,N} \to u_i & \text{ strongly in } L^p(0,T;L^p(\Omega_i)), \\
w_{j,N} \to w_j & \text{ strongly in } L^p(0,T;L^p(e_j)),
\end{align*}
so that along a further non-relabelled subsequence, 
\begin{align*}
u_{i,N} \to u_i & \text{ a.e in } \Omega_i \times (0,T), \\
w_{j,N} \to w_j & \text{ a.e in } e_j \times (0,T).
\end{align*}
We now perform a standard argument to pass to the limit $N \to \infty$ and deduce that the limit functions $(u_i, w_j, z_k)$ is a weak solution to the model in the sense of Definition \ref{defn:weaksoln}.  We fix an index $s \in \N$ and multiply \eqref{fp:u} by $\zeta(t) \phi^{(i)}_s$ for an arbitrary $\zeta \in C^\infty_c(0,T)$, integrating over $(0,T)$ and performing integration by parts to arrive at
\begin{equation}\label{lim:Gal:uN}
\begin{aligned}
0 & = \int_0^T \int_{\Omega_i} \Big ( \pd_t u_{i,N} \phi^{(i)}_s + \kappa(\nabla u_{i,N}) \cdot \nabla \phi^{(i)}_s + f(u_{i,N}) \phi^{(i)}_s \Big ) \zeta(t) \, dx \, dt \\
& \quad - \sum_{j \, : \, e_j \sim \Omega_i}^{|E|} \int_{0}^T \int_{e_j} (\alpha_{ij} w_{j,N} - \beta_{ij} u_{i,N}) \phi^{(i)}_s \zeta(t) \, dx \, dt.
\end{aligned}
\end{equation}
Passing to the limit $N \to \infty$ in the linear terms is straightforward, and so we focus on the nonlinear terms. We use the a.e.~convergence of $u_{i,N}$ to $u_i$ and the continuity of $f$ to deduce that 
\[
f(u_{i,N}) \to f(u_i) \text{ a.e.~in } \Omega_i \times (0,T).
\]
By the Gagliardo--Nirenberg inequality we have the continuous embedding
\[
L^\infty(0,T;L^2(\Omega_i)) \cap L^p(0,T;W^{1,p}(\Omega_i)) \subset L^{2p}(0,T; L^{2p}(\Omega_i)) \cong L^{2p}(\Omega_i \times (0,T)),
\]
and thus $\{u_{i,N}\}_{n \in \N}$ and the limit function $u_i$ are bounded in $L^{2p}(\Omega_i \times (0,T))$.  Via interpolation theory of Lebesgue spaces it holds that for $z \in (p, 2p)$,
\[
\| u_{i,N} - u_i \|_{L^z(\Omega_i \times (0,T))} \leq \| u_{i,N} - u_i \|_{L^{2p}(\Omega_i \times (0,T))}^{\frac{2(z-p)}{p}} \| u_{i,N} - u_i \|_{L^{p}(\Omega_i \times (0,T))}^{\frac{2p-z}{z}} \to 0
\]
as $N \to \infty$.  Hence, we find that
\begin{align}\label{uiN:l2p}
u_{i,N} \to u_i \text{ strongly in } L^z(\Omega_i \times (0,T)), \quad \forall z < 2p.
\end{align}
Then, as $2 + \frac{2(p-1)^2}{p} < 2p$, for the exponents $q$ in the range specified by \eqref{ass:f} we pick $z \in (q, 2p)$ and apply H\"older's inequality to see that
\begin{align*}
& \int_0^T \int_{\Omega_i} |u_{i,N} - u_{i}|^q |\phi^{(i)}_s| \, dx \, dt \leq \int_0^T \| u_{i,N} - u_i \|_{L^z(\Omega_i)}^{q} \| \phi^{(i)}_s \|_{L^{\frac{z}{z-q}}(\Omega_i)} \, dt \\
& \quad \leq C\| \phi^{(i)}_s \|_{W^{1,p}(\Omega_i)} \int_0^T \| u_{i,N} - u_i \|_{L^z(\Omega_i)}^q  \, dt \\
& \quad \leq C\| \phi^{(i)}_s \|_{W^{1,p}(\Omega_i)} T^{\frac{z-q}{z}} \| u_{i,N} - u_i \|_{L^z(\Omega_i \times (0,T))}^q \to 0,
\end{align*}
as $N \to \infty$.  This implies that 
\[
|u_{i,N}|^q |\phi^{(i)}_s| \to |u_i|^q |\phi^{(i)}_s| \text{ strongly in } L^1(\Omega_i \times (0,T)),
\]
and so, by the generalized Lebesgue dominated convergence theorem and \eqref{ass:f} we find that 
\[
\int_0^T \int_{\Omega_i} f(u_{i,N}) \phi^{(i)}_s \zeta(t)  \, dx \, dt \to \int_0^T \int_{\Omega_i} f(u_i) \phi^{(i)}_s \zeta(t) \, dx \, dt.
\]
Similarly, we see that $q + 1 < 3 + \frac{2(p-1)^2}{p} \leq 2p$ for $p \geq 2$ and so by \eqref{uiN:l2p} and the generalized Lebesgue dominated convergence theorem we infer
\begin{align}\label{fuin:uin}
\int_0^T \int_{\Omega_i} f(u_{i,N}) u_{i,N} \, dx \, dt \to \int_0^T \int_{\Omega_i} f(u_i) u_i \, dx \, dt.
\end{align}
For the anisotropic diffusion term we argue as follows. From \eqref{ass:aniso} and the uniform boundedness of $\nabla u_{i,N}$ in $L^p(\Omega_i \times (0,T))$, it holds that 
\[
\| \kappa(\nabla u_{i,N}) \|_{L^{\frac{p}{p-1}}(\Omega_i \times (0,T))}^{\frac{p}{p-1}} \leq c_2 \| \nabla u_{i,N} \|_{L^p(\Omega_i \times (0,T))}^p \leq C,
\]
and so there exists a limit function $\bm{\kappa} \in L^{\frac{p}{p-1}}(\Omega_i \times (0,T))$ such that 
\[
\kappa(\nabla u_{i,N}) \to \bm{\kappa} \text{ weakly in } L^{\frac{p}{p-1}}(\Omega_i \times (0,T))
\]
as $N \to \infty$ along a non-relabelled subsequence. We now apply Minty's trick to identify $\bm{\kappa} = \kappa(\nabla u_i)$. Passing to the limit $N \to \infty$ in \eqref{lim:Gal:uN} yields
\begin{equation}\label{lim:Gal:u}
\begin{aligned}
 0 & = \int_0^T \Big ( \langle \pd_t u_i, \phi^{(i)}_s \rangle_{W^{1,p}(\Omega_i)}  + \int_{\Omega_i} \bm{\kappa} \cdot \nabla \phi^{(i)}_s + f(u_i) \phi^{(i)}_s \, dx \Big ) \zeta(t) \, dt   \\
 & \quad - \sum_{j \, : \, e_j \sim \Omega_i}^{|E|} \int_0^T \int_{e_j} (\alpha_{ij} w_j - \beta_{ij} u_i) \phi^{(i)}_s \zeta(t) \, dx \, dt.
\end{aligned}
\end{equation}
By a standard density argument we obtain from the above that 
\[
 0 = \langle \pd_t u_i, \phi^{(i)} \rangle_{W^{1,p}(\Omega_i)}  + \int_{\Omega_i} \bm{\kappa} \cdot \nabla \phi^{(i)} + f(u_i) \phi^{(i)} \, dx   - \sum_{j \, : \, e_j \sim \Omega_i}^{|E|} \int_{e_j} (\alpha_{ij} w_j - \beta_{ij} u_i) \phi^{(i)} \, dx 
\]
holding for a.e.~$t \in (0,T)$ and for all $\phi^{(i)} \in W^{1,p}(\Omega_i)$. Choosing $\phi^{(i)} = u_i$ then yields the identity
\begin{equation}\label{minty0}
\begin{aligned}
 0 & = \langle \pd_t u_i, u_i \rangle_{W^{1,p}(\Omega_i)}  + \int_{\Omega_i} \bm{\kappa} \cdot \nabla u_i + f(u_i) u_i \, dx    \\
 & \quad - \sum_{j \, : \, e_j \sim \Omega_i}^{|E|} \int_{e_j} (\alpha_{ij} w_j - \beta_{ij} u_i) u_i \, dx 
\end{aligned}
\end{equation}
Next, testing \eqref{fp:u} by $u_{i,N}$ yields the identity
\begin{equation}\label{minty1}
\begin{aligned}
\int_0^T \int_{\Omega_i} \kappa( \nabla u_{i,N}) \cdot \nabla u_{i,N} \, dx \, dt & = - \int_0^T \int_{\Omega_i} \pd_t u_{i,N} u_{i,N} + f(u_{i,N}) u_{i,N} \, dx \,dt \\
& \quad + \sum_{j \, : \, e_j \sim \Omega_i}^{|E|} \int_0^T \int_{e_j} (\alpha_{ij} w_{j,N} - \beta_{ij} u_{i,N}) u_{i,N} \, dx \, dt.
\end{aligned}    
\end{equation}
Then, for arbitrary $v \in L^p(0,T;W^{1,p}(\Omega_i))$, we have by \eqref{minty1} and the monotonicity property in \eqref{ass:aniso},
\begin{align*}
0 & \leq \int_0^T \int_{\Omega_i} (\kappa(\nabla u_{i,N}) - \kappa(\nabla v) ) \cdot \nabla (u_{i,N} - v) \, dx \, dt \\
& = \int_0^T \int_{\Omega_i} \kappa(\nabla u_{i,N}) \cdot \nabla u_{i,N} - \kappa(\nabla v) \cdot \nabla u_{i,N} - \kappa(\nabla u_{i,N}) \cdot \nabla v + \kappa(\nabla v) \cdot \nabla v \, dx \,dt  \\
& = \int_0^T \int_{\Omega_i}  - \kappa(\nabla v) \cdot \nabla u_{i,N} - \kappa(\nabla u_{i,N}) \cdot \nabla v + \kappa(\nabla v) \cdot \nabla v \, dx \,dt \\
& \quad - \int_0^T \int_{\Omega_i} (\pd_t u_{i,N} + f(u_{i,N})) u_{i,N} \, dx \, dt + \sum_{j \, : \, e_j \sim \Omega_i}^{|E|} \int_0^T \int_{e_j} (\alpha_{ij} w_{j,N} - \beta_{ij} u_{i,N}) u_{i,N} \, dx \, dt.
\end{align*}
Passing to the limit $N \to \infty$ and using \eqref{fuin:uin} leads to
\begin{equation}\label{minty2}
\begin{aligned}
0 & \leq \int_0^T \int_{\Omega_i} - \kappa(\nabla v) \cdot \nabla u_i - \bm{\kappa} \cdot \nabla v + \kappa(\nabla v) \cdot \nabla v \, dx \, dt \\
& \quad - \int_0^T \langle \pd_t u_i, u_i \rangle_{W^{1,p}(\Omega_i)} \, dt - \int_0^T \int_{\Omega_i} f(u_i) u_i \, dx \, dt \\
& \quad + \sum_{j \, : \, e_j \sim \Omega_i}^{|E|} \int_0^T \int_{e_j} (\alpha_{ij} w_{j} - \beta_{ij} u_i) u_i \, dx \, dt
\end{aligned}
\end{equation}
Comparing with \eqref{minty0} yields 
\begin{equation}\label{minty3}
\begin{aligned}
0 & \leq \int_0^T \int_{\Omega_i} \bm{\kappa} \cdot \nabla u_i - \kappa(\nabla v) \cdot \nabla u_i - \bm{\kappa} \cdot \nabla v + \kappa(\nabla v) \cdot \nabla v \, dx \, dt \\
& = \int_0^T \int_{\Omega_i} (\bm{\kappa} - \kappa(\nabla v)) \cdot \nabla (u_i - v) \, dx \, dt
\end{aligned}
\end{equation}
holding for arbitrary $v \in L^p(0,T;W^{1,p}(\Omega_i))$. Choosing $v = u_i \pm \alpha U$ for positive constant $\alpha$ and arbitrary $U \in L^p(0,T;W^{1,p}(\Omega_i))$ in \eqref{minty3} results in
\[
0 \leq \mp  \int_0^T \int_{\Omega_i} (\bm{\kappa} - \kappa(\nabla (u_i \pm \alpha U))) \cdot \nabla U \, dx \, dt.
\]
Sending $\alpha \to 0$ and using the continuity of $\kappa$ leads to
\[
0 \geq \pm \int_0^T \int_{\Omega_i} (\bm{\kappa} - \kappa(\nabla u_i)) \cdot \nabla U \, dx \, dt, \quad \forall U \in L^p(0,T;W^{1,p}(\Omega_i)),
\]
which allows us to infer that $\bm{\kappa} = \kappa(\nabla u_i)$ a.e.~in $\Omega_i \times (0,T)$. With this identification of $\bm{\kappa}$, upon returning to \eqref{lim:Gal:u} we see that
\begin{equation}\label{u:lim:eq}
\begin{aligned}
0 & =  \langle \pd_t u_i, \phi^{(i)} \rangle_{W^{1,p}(\Omega_i)}  + \int_{\Omega_i} \kappa(\nabla u_i) \cdot \nabla \phi^{(i)} + f(u_i) \phi^{(i)} \, dx \\
& \quad - \sum_{j \, : \, e_j \sim \Omega_i}^{|E|} \int_{e_j} (\alpha_{ij} w_j - \beta_{ij} u_i) \phi^{(i)} \, dx 
 \end{aligned}
\end{equation}
holds for a.e.~$t \in (0,T)$ and for all $\phi^{(i)} \in W^{1,p}(\Omega_i)$. For the edge solutions, we invoke the Gagliardo--Nirenberg inequality to deduce that
\[
L^\infty(0,T;L^2(e_j)) \cap L^p(0,T;W^{1,p}(e_j)) \subset L^{3p}(e_j \times (0,T)).
\]
Then, by a similar argument we obtain
\begin{align*}
    w_{j,N} \to w_{j} & \text{ strongly in } L^{z}(e_j \times (0,T)), \quad \forall z < 3p, 
\end{align*}
along with
\begin{align*}
\int_0^T \int_{e_j} g(w_{j,N}) \psi^{(j)}_s \zeta(t) \, dx \, dt & \to \int_0^T \int_{e_j} g(w_j) \psi^{(j)}_s \zeta(t) \, dx \, dt, \\
\int_0^T \int_{e_j} g(w_{j,N}) w_{j,N} \, dx \, dt & \to \int_0^T \int_{e_j} g(w_j) w_j \, dx \, dt 
\end{align*}
as $N \to \infty$. In a similar fashion, from \eqref{ass:edgeaniso} and the uniform boundedness of $\pd_x w_{j,N}$ in $L^p(e_j \times (0,T))$ it holds that
\[
\| \eta(\pd_x w_{j,N}) \|_{L^{\frac{p}{p-1}}(e_j \times (0,T))}^{\frac{p}{p-1}} \leq C \| \pd_x w_{j,N} \|_{L^p(e_j \times (0,T))}^p \leq C,
\]
and so there exists a limit function $\eta \in L^{\frac{p}{p-1}}(e_j \times (0,T))$ such that 
\[
\eta(\pd_x w_{j,N}) \to \eta \text{ weakly in } L^{\frac{p}{p-1}}(e_j \times (0,T))
\]
as $N \to \infty$ along a non-relablled subsequence. With Minty's trick we identify $\eta = \eta(\pd_x w_j)$, so that upon multiplying \eqref{fp:w} by $\zeta(t) \psi^{(j)}_s$ for an arbitrary $\zeta \in C^\infty_c(0,T)$, integrating over $(0,T)$ and passing to the limit $N \to \infty$, we arrive at
\begin{equation}
\begin{aligned}
0 & = \int_0^T \zeta(t) \Big ( \langle \pd_t w_{j,N}, \psi_s^{(j)} \rangle_{W^{1,p}(e_j)} + \int_{e_j} \eta(\pd_x w_j) \pd_x \psi_s^{(j)} + g(w_j) \psi^{(j)} \, dx \Big ) \, dt \\
& \quad + \int_0^T \zeta(t) \sum_{i \, : \, \Omega_i \sim e_j}^{M} \int_{e_j} (\alpha_{ij} w_j - \beta_{ij} u_i) \psi_s^{(j)} \, dx \, dt\\
& \quad + \int_0^T \zeta(t) \sum_{k \, : \, v_k \sim e_j}^{|V|} \Big ( \Big [ \delta^{k}_{e_j} + \sum_{e_m \sim v_k, m \neq j}^{|E|} \gamma_{j \to m}^k \Big ] w_j(v_k) - \lambda^k_{e_j} z_k \\
& \qquad \qquad - \sum_{e_m \sim v_k, m \neq j}^{|E|} \gamma^k_{m \to j} w_m(v_k) \Big ) \psi_s^{(j)}(v_k) \, dt.
\end{aligned}
\end{equation}
Hence, via a standard density argument, it holds that
\begin{equation}
\begin{aligned}
0 & = \langle \pd_t w_{j,N}, \psi^{(j)} \rangle_{W^{1,p}(e_j)} + \int_{e_j} \eta(\pd_x w_j) \pd_x \psi^{(j)} + g(w_j) \psi^{(j)} \, dx \\
& \quad + \sum_{i \, : \, \Omega_i \sim e_j}^{M} \int_{e_j} (\alpha_{ij} w_j - \beta_{ij} u_i) \psi^{(j)} \, dx \\
& \quad + \sum_{k \, : \, v_k \sim e_j}^{|V|} \Big ( \Big [ \delta^{k}_{e_j} + \sum_{e_m \sim v_k, m \neq j}^{|E|} \gamma_{j \to m}^k \Big ] w_j(v_k) - \lambda^k_{e_j} z_k \\
& \qquad \qquad - \sum_{e_m \sim v_k, m \neq j}^{|E|} \gamma^k_{m \to j} w_m(v_k) \Big ) \psi^{(j)}(v_k)
\end{aligned}
\end{equation}
for a.e.~$t \in (0,T)$ and for all $\psi^{(j)} \in W^{1,p}(e_j)$. Lastly, we multiply \eqref{fp:z} by an arbitrary $\zeta \in L^2(0,T)$ and integrate over $(0,T)$ to see that
\[
0 = \int_0^T \Big ( \frac{d}{dt} z_{k,N} -  \sum_{j \, : \, e_j \sim v_k} (\delta^{k}_{e_j} w_{j,N}(v_k) - \lambda^k_{e_j} z_{k,N}) \Big ) \zeta \, dt
\]
By the boundary trace theorem we have that
\[
w_{j,N} \to w_j \text{ weakly in } L^p(0,T;W^{1-\frac{1}{p},p}(v_k)),
\]
and so passing to the limit $N \to \infty$ we infer 
\[
0 = \int_0^T \Big ( \frac{d}{dt} z_k  - \sum_{j \, : \, e_j \sim v_k} (\delta^{k}_{e_j} w_{j}(v_k) - \lambda^k_{e_j} z_{k}) \Big ) \zeta \, dt,
\]
which then implies
\[
\frac{d}{dt} z_k = \sum_{j \, : \, e_j \sim v_k} (\delta^{k}_{e_j} w_{j}(v_k) - \lambda^k_{e_j} z_{k}) 
\]
holds for a.e.~$t \in (0,T)$. 

It remains to discuss the attainment of the initial conditions. For the case $p =2$ we use the continuous embedding (Lions--Magenes lemma)
\[
L^2(0,T;H^1(X)) \cap H^1(0,T;H^1(X)^*) \subset C^0([0,T];L^2(X))
\]
for $X = \Omega_i$ or $e_j$ to deduce that
\[
u_{i} \in C^0([0,T];L^2(\Omega_i)), \quad w_j \in C^0([0,T];L^2(e_j)).
\]
Then, for fixed $1 \leq s \leq N$, we test \eqref{fp:u} by $\zeta(t) \phi^{(i)}_s$ with an arbitrary $\zeta \in C^\infty(0,T)$ such that $\zeta(T) = 0$. Integrating by parts yields
\begin{align*}
0& = \int_0^T \int_{\Omega_i} - (u_{i,N} - P_{\mathbb{U}^{(i)}_N}(u_{i,0}))\phi^{(i)}_s \pd_t \zeta + \Big ( \kappa(\nabla u_{i,N}) \cdot \nabla \phi^{(i)}_s + f(u_{i,N}) \phi^{(i)}_s \Big ) \zeta \, dx \, dt \\
& \quad - \int_0^T \sum_{j \, : \, e_j \sim \Omega_i}^{|E|} \int_{e_j} (\alpha_{ij} w_{j,N} - \beta_{ij} u_{i,N}) \phi^{(i)}_s \zeta \, dx \, dt.
\end{align*}
Passing to the limit $N \to \infty$ and using that $P_{\mathbb{U}^{(i)}_N}(u_{0,i}) \to u_{0,i}$ strongly in $L^2(\Omega_i)$ as $N \to \infty$, we find
\begin{equation}\label{u:ini:1}
\begin{aligned}
0& = \int_0^T \int_{\Omega_i} - (u_{i} - u_{i,0})\phi^{(i)}_s  \zeta'(t) + \Big ( \kappa(\nabla u_{i}) \cdot \nabla \phi^{(i)}_s + f(u_{i}) \phi^{(i)}_s \Big ) \zeta \, dx \, dt \\
& \quad - \int_0^T \sum_{j \, : \, e_j \sim \Omega_i}^{|E|} \int_{e_j} (\alpha_{ij} w_{j} - \beta_{ij} u_{i}) \phi^{(i)}_s \zeta \, dx \, dt.
\end{aligned}
\end{equation}
Comparing \eqref{u:ini:1} with the following identity obtained from choosing $\phi^{(i)} = \phi^{(i)}_s$ in \eqref{u:lim:eq}, then multiplying with $\zeta \in C^\infty(0,T)$ such that $\zeta(T) = 0$:
\begin{align*}
0 & = \int_0^T \int_{\Omega_i} - (u_i - u_i(0)) \phi^{(i)}_s \zeta'(t) + \Big ( \kappa(\nabla u_i) \cdot \nabla \phi^{(i)}_s + f(u_i) \phi^{(i)}_s \Big ) \zeta \, dx \, dt \\
& \quad - \int_0^T \sum_{j \, : \, e_j \sim \Omega_i}^{|E|} \int_{e_j} (\alpha{ij} w_j - \beta_{ij} u_i) \phi^{(i)}_s \zeta \, dx \, dt,
\end{align*}
we find that
\[
\int_0^T \int_{\Omega_i} \zeta'(t) (u_{i,0} - u_i(0)) \phi^{(i)}_s \, dx \, dt = 0.
\]
By a density argument this implies that
\[
\int_0^T \int_{\Omega_i} \zeta'(t) (u_{i,0} - u_{i}(0)) \phi^{(i)} \, dx \, dt = 0
\]
for all $\phi^{(i)} \in H^1(\Omega_i)$.  Hence we deduce that $u_{i}(0) = u_{i,0}$ a.e.~in $\Omega_i$. A similar argument can be used to show $w_{j}(0) = w_{j,0}$ a.e.~in $e_j$, while by construction we have $z_k(0) = z_{k,0}$. 

For the case $p > 2$, there is no equivalent of the Lions--Magenes lemma as $W^{1,p}(X)$ is not a Hilbert space. Hence, we use the continuous embedding
\[
W^{1,\frac{p}{p-1}}(0,T;W^{1,p}(X)^*) \subset C^0([0,T];W^{1,p}(X)^*)
\]
for $X = \Omega_i$ or $e_j$ to deduce that 
\[
u_i \in C^0([0,T];W^{1,p}(\Omega_i)^*), \quad w_j \in C^0([0,T];W^{1,p}(e_j)^*).
\]
Then, by an analogous argument to the above setting with $p = 2$, we infer that
\[
\int_0^T \zeta'(t) \langle u_{i,0} - u_{i}(0),  \phi^{(i)} \rangle_{W^{1,p}(\Omega_i)} \, dt = 0
\]
for all $\phi^{(i)} \in W^{1,p}(\Omega_i)$. This shows
\[
\langle u_{i}(0), \phi^{(i)} \rangle_{W^{1,p}(\Omega_i)} = \langle u_{i,0}, \phi^{(i)} \rangle_{W^{1,p}(\Omega_i)}, \quad \forall \phi^{(i)} \in W^{1,p}(\Omega_i),
\]
and similarly we deduce
\[
\langle w_{j}(0), \psi^{(j)} \rangle_{W^{1,p}(e_j)} = \langle w_{j,0}, \psi^{(j)} \rangle_{W^{1,p}(e_j)}, \quad \forall \psi^{(j)} \in W^{1,p}(e_j).
\]
\subsection{Uniqueness and continuous dependence}
Let $\{\bm{u}_m, \bm{w}_m, \bm{z}_m\}_{m=1,2}$ denote two weak solutions to the model corresponding to two sets of initial data $\{\bm{u}_{0,m}, \bm{w}_{0,m}, \bm{z}_{0,m}\}_{m=1,2}$. Their differences
\[
\overline{u}_i = u_{i,1} - u_{i,2}, \quad \overline{w}_j = w_{j,1} - w_{j,2}, \quad \overline{z}_k = z_{k,1} - z_{k,2} 
\]
satisfy
\begin{subequations}
\begin{alignat}{2}
\label{diff:u} 0 & =  \langle \pd_t \overline{u}_i, \phi^{(i)} \rangle_{W^{1,p}(\Omega_i)}  + \int_{\Omega_i} (\kappa(\nabla u_{i,1}) - \kappa(\nabla u_{i,2})) \cdot \nabla \phi^{(i)} \, dx \\
\notag & \quad + \int_{\Omega_i} (f(u_{i,1})- f(u_{i,2})) \phi^{(i)} \, dx - \sum_{j \, : \, e_j \sim \Omega_i}^{|E|} \int_{e_j} (\alpha_{ij} \overline{w}_j - \beta_{ij} \overline{u}_i) \phi^{(i)} \, dx, \\
\label{diff:w} 0 & = \langle \pd_t \overline{w}_{j,N}, \psi^{(j)} \rangle_{W^{1,p}(e_j)} + \int_{e_j} (\eta(\pd_x w_{j,1}) - \eta(\pd_x w_{j,2})) \pd_x \psi^{(j)} \, dx \\
\notag & \quad + \int_{e_j} (g(w_{j,1}) - g(w_{j,2})) \psi^{(j)} \, dx + \sum_{i \, : \, \Omega_i \sim e_j}^{M} \int_{e_j} (\alpha_{ij} \overline{w}_j - \beta_{ij} \overline{u}_i) \psi^{(j)} \, dx \\
\notag & \quad + \sum_{k \, : \, v_k \sim e_j}^{|V|} \Big ( \Big [ \delta^{k}_{e_j} + \sum_{e_m \sim v_k, m \neq j}^{|E|} \gamma_{j \to m}^k \Big ] \overline{w}_j(v_k) - \lambda^k_{e_j} \overline{z}_k \\
\notag & \qquad \qquad - \sum_{e_m \sim v_k, m \neq j}^{|E|} \gamma^k_{m \to j} \overline{w}_m(v_k) \Big ) \psi^{(j)}(v_k), \\
\label{diff:z} 0 & = \frac{d}{dt} \overline{z}_k - \sum_{j \, : \, e_j \sim v_k} (\delta^{k}_{e_j} \overline{w}_{j}(v_k) - \lambda^k_{e_j} \overline{z}_{k}),
\end{alignat}
\end{subequations}
holding for a.e.~$t \in (0,T)$ and for all $\phi^{(i)} \in W^{1,p}(\Omega_i)$ and $\psi^{(j)} \in W^{1,p}(e_j)$.  We choose $\phi^{(i)} = \overline{u}_i$ in \eqref{diff:u}, $\psi^{(j)} = \overline{w}_j$ in \eqref{diff:w}, and test \eqref{diff:z} with $\overline{z}_k$, utilizing the monotonicity of $\kappa$, $\eta$, $f$ and $g$, as well as the estimates \eqref{rhs:3}, \eqref{rhs:4} and \eqref{rhs:5}, we obtain an inequality akin to \eqref{unifest:1:gron} for the differences:
\begin{align*}
&  \sum_{k=1}^{|V|} \frac{1}{2} \frac{d}{dt} |\overline{z}_{k}|^2  + \sum_{i=1}^M \Big (\frac{d}{dt} \frac{1}{2} \| \overline{u}_{i} \|^2_{L^2(\Omega_i)} + c_3 \|\nabla \overline{u}_{i} \|_{L^p(\Omega_i)}^p \Big ) \\
& \qquad + \sum_{j=1}^{|E|} \Big (\frac{d}{dt} \frac{1}{2} \| \overline{w}_{j} \|^2_{L^2(e_j)} + c_7 \| \pd_x \overline{w}_{j} \|_{L^p(e_j)}^p \Big ) \\
& \qquad + \sum_{j=1}^{|E|} \Big ( \sum_{i \, : \, \Omega_i \sim e_j}^{M} \frac{\beta_{ij}}{2} \| \overline{u}_{i} \|_{L^2(e_j)}^2 + \sum_{k \, : \, v_k \sim e_j}^{|V|} \frac{\delta^{k}_{e_j}}{4} |\overline{w}_{j}(v_k)|^2 \Big ) \\
& \quad \leq C \sum_{j=1}^{|E|} \| \overline{w}_{j} \|_{L^2(e_j)}^2 + C \sum_{k=1}^{|V|} |\overline{z}_{k}|^2.
\end{align*}
By Gronwall's inequality we deduce that 
\begin{align*}
& \sup_{t \in (0,T)} \Big ( \sum_{i=1}^M \| \overline{u}_i(t) \|_{L^2(\Omega_i)}^{2} + \sum_{j=1}^{|E|} \| w_{j}(t) \|_{L^2(e_j)}^2 + \sum_{k=1}^{|V|} |\overline{z}_k(t)|^2 \Big ) \\
& \qquad + \int_0^T \Big ( \sum_{i=1}^M \| \nabla u_i \|_{L^p(\Omega_i)}^p + \| \pd_x w_j \|_{L^p(e_j)}^p \Big ) \, dt \\
& \quad \leq C  \Big ( \sum_{i=1}^M \| \overline{u}_{i,0} \|_{L^2(\Omega_i)}^{2} + \sum_{j=1}^{|E|} \| w_{j,0} \|_{L^2(e_j)}^2 + \sum_{k=1}^{|V|} |\overline{z}_{k,0}|^2 \Big ).
\end{align*}
When the initial data coincide, we see that $\overline{u}_i = 0$, $\overline{w}_j = 0$ and $\overline{z}_k = 0$, which yields the uniqueness of solutions.

\subsection{Unpopulated vertices setting with symmetric inter-edge transfer}
In the setting where all vertices are unpopulated and \eqref{ass:sym:coeff} holds for the inter-edge transfer coefficients, the analysis of the model simplifies, as \eqref{weakform:alt} reduces to 
\begin{equation}\label{weakform:alt:unpop}
\begin{aligned}
0 & = \sum_{i=1}^M \Big (  \langle \pd_t u_{i},  \phi^{(i)} \rangle_{W^{1,p}(\Omega_i)} + \int_{\Omega_i} \kappa(\nabla u_{N}) \cdot \nabla \phi^{(i)} + f(u_{i}) \phi^{(i)} \, dx \Big ) \\
& \quad + \sum_{j=1}^{|E|} \Big ( \langle \pd_t w_{j},  \psi^{(j)} \rangle_{W^{1,p}(e_j)} + \int_{e_j} \eta(\pd_x w_{j}) \pd_x \psi^{(j)} + g(w_{j}) \psi^{(j)}  \, dx \\
& \qquad  \qquad + \sum_{i \,  :\, \Omega_i \sim e_j}^M \int_{e_j} (\alpha_{ij} w_{j} - \beta_{ij} u_{i})(\psi^{(j)} - \phi^{(i)}) \, dx \Big ) \\
& \quad + \sum_{k=1}^{|V|} \sum_{j \neq m \, : \, e_j, e_m \sim v_k}^{|E|} \gamma^k_{j \leftrightarrow  m}\Big [ ( w_j - w_{m} ) (\psi^{(j)} - \psi^{(m)}) \Big ] \Big \vert_{v_k}.
\end{aligned}
\end{equation}
Recalling that $E$ denotes the full metric graph, we may then appeal to a standard Galerkin approximation to deduce the well-posedness of weak solutions $\bm{u} = (u_1, \dots, u_M)$ and $w$ satisfying
\begin{align*}
u_i & \in L^\infty(0,T;L^2(\Omega_i)) \cap L^p(0,T;W^{1,p}(\Omega_i)) \cap W^{1,\frac{p}{p-1}}(0,T;W^{1,p}(\Omega_i)^*), \quad \text{ for } 1 \leq i \leq M,\\
w & \in L^\infty(0,T;L^2(E)) \cap L^p(0,T;W^{1,p}(E)) \cap W^{1,\frac{p}{p-1}}(0,T;W^{1,p}(E)^*),
\end{align*}
and \eqref{weakform:alt:unpop} holds with $w \vert_{e_j} = w_j$ for a.e.~$t \in (0,T)$ and for all $\phi^{(i)} \in W^{1,p}(\Omega_i)$ and for all $\psi \in W^{1,p}(E)$ with $\psi^{(j)} = \psi \vert_{e_j}$. The initial conditions are also attained in the sense of Definition \ref{defn:weaksoln}. 

\section{Regularity of solutions}\label{sec:reg}
\subsection{Temporal regularity}
We return to the semi-Galerkin system \eqref{fixedpoint:sys} and test \eqref{fp:u} with $\pd_t u_{i,N}$, \eqref{fp:w} with $\pd_t w_{j,N}$, so that upon summing we obtain an analogue of \eqref{fixedpoint:compact2}:
\begin{equation}\label{reg:1}
\begin{aligned}
&  \sum_{i=1}^M \Big ( \| \pd_t u_{i,N} \|_{L^2(\Omega_i)}^2 + \frac{d}{dt} \int_{\Omega_i} \widehat{\kappa} (\nabla u_{i,N}) + \widehat{f}(u_{i,N}) \, dx\Big ) \\
& \qquad + \sum_{j=1}^{|E|} \Big ( \| \pd_t w_{j,N} \|_{L^2(e_j)}^2 + \frac{d}{dt} \int_{e_j} \widehat{\eta}(\pd_x w_{j,N}) + \widehat{g}(w_{j,N}) \, dx  \\
& \qquad  + \sum_{i \, : \, \Omega_i \sim e_j}^{M} \int_{e_j} \frac{\alpha_{ij}}{2} \frac{d}{dt}|w_{j,N}|^2 + \frac{\beta_{ij}}{2} \frac{d}{dt}|u_{i,N}|^2 - \alpha_{ij} w_{j,N} \pd_t u_{i,N} - \beta_{ij} u_{i,N} \pd_t w_{j,N} \, dx \Big ) \\
& \qquad + \sum_{j=1}^{|E|} \sum_{k \, : \, v_k \sim e_j}^{|V|} \frac{d}{dt} \frac{1}{2} \Big ( \delta^{k}_{e_j} |w_{j,N}(v_k)|^2 + \sum_{e_m \sim v_k, m \neq j} \gamma_{j \leftrightarrow m}^k |w_{j,N}(v_k) - w_{m,N}(v_k)|^2  \Big )\\
& = \sum_{j=1}^{|E|}\sum_{k \, : \, v_k \sim e_j}^{|V|} \lambda^k_{e_j} z_{k,N} \pd_t w_{j,N}(v_k).
\end{aligned}
\end{equation}
Integrating \eqref{reg:1} over $(0,s)$ for arbitrary $s \in (0,T)$ yields
\begin{equation}\label{reg:2}
\begin{aligned}
& \sum_{i=1}^M \Big ( \int_{\Omega_i} \widehat{\kappa}(\nabla u_{i,N}(s)) + \widehat{f}(u_{i,N}(s)) \, dx + \int_0^s \| \pd_t u_{i,N} \|_{L^2(\Omega_i)}^2 \, dt \Big ) \\
& \qquad + \sum_{j=1}^{|E|} \Big ( \int_{e_j} \widehat{\eta}(\pd_x w_{j,N}(s)) + \widehat{g}(w_{j,N}(s)) \, dx + \int_0^s \| \pd_t w_{j,N} \|_{L^2(e_j)}^2 \, dt \Big ) \\
& \qquad + \sum_{j=1}^{|E|} \sum_{i \, : \, \Omega_i \sim e_j}^{M} \Big ( \frac{\alpha_{ij}}{2} \| w_{j,N}(s) \|_{L^2(e_j)}^2 + \frac{\beta_{ij}}{2} \| u_{i,N}(s) \|_{L^2(e_j)}^2 \Big ) \\
& \qquad + \sum_{j=1}^{|E|} \sum_{k \, : \, v_k \sim e_j}^{|V|} \frac{1}{2} \Big ( \delta^k_{e_j} |w_{j,N}(s,v_k)|^2 + \sum_{e_m \sim v_k, m \neq j}^{|V|} \gamma^k_{j \leftrightarrow m} |w_{j,N}(s,v_k) - w_{m,N}(s,v_k)|^2 \Big ) \\
& \quad = \sum_{i=1}^{M} \Big ( \int_{\Omega_i} \widehat{\kappa}(\nabla u_{i,N}(0)) + \widehat{f}(u_{i,N}(0)) \, dx \Big ) + \sum_{j=1}^{|E|} \sum_{i \, : \, \Omega_i \sim e_j}^{M} \frac{\beta_i}{2} \| u_{i,N}(0) \|_{L^2(e_j)}^2 \\
& \qquad + \sum_{j=1}^{|E|} \Big ( \int_{e_j} \widehat{\eta}(\pd_x w_{j,N}(0)) + \widehat{g}(w_{j,N}(0)) \, dx \Big ) + \sum_{j=1}^{|E|} \sum_{i \, : \, \Omega_i \sim e_j}^{M}  \frac{\alpha_{ij}}{2} \| w_{j,N}(0) \|_{L^2(e_j)}^2 \\
& \qquad + \sum_{j=1}^{|E|} \sum_{k \, : \, v_k \sim e_j} \frac{1}{2} \Big ( \delta^k_{e_j} |w_{j,N}(0,v_k)|^2 + \sum_{e_m \sim v_k, m \neq j}^{|V|} \gamma^k_{j \leftrightarrow m} |w_{j,N}(0,v_k) - w_{m,N}(0,v_k)|^2 \Big ) \\
& \qquad + \sum_{j=1}^{|E|} \sum_{i \, : \, \Omega_i \sim e_j}^{M} \Big ( \int_0^s \int_{e_j} \beta_{ij} u_{i,N} \pd_t w_{j,N} - \alpha_{ij} \pd_t w_{j,N} u_{i,N} \, dx \, dt \\
& \qquad \qquad \qquad + \int_{e_j} \alpha_{ij} \big ( w_{j,N}(s) u_{i,N}(s) - w_{j,N}(0) u_{i,N}(0) \big ) \, dx \Big ) \\
& \qquad + \sum_{j=1}^{|E|} \sum_{k \, : \, v_k \sim e_j} \lambda^k_{e_j} \Big ( z_{k,N}(s) w_{j,N}(s,v_k) - z_{k,N}(0) w_{j,N}(0,v_k) \\
& \qquad \qquad \qquad - \int_0^s \frac{d}{dt} z_{k,N} w_{j,N}(v_k) \, dt \Big ) \\
& \quad =: I_1 + \cdots + I_7,
\end{aligned}
\end{equation}
where $I_i$ denotes the $i$th line of the right-hand side. Using \eqref{ass:aniso}, \eqref{ass:edgeaniso}, \eqref{ass:f} and \eqref{ass:g} and the trace theorem, we see that 
\[
|I_1 + I_2 + I_3| \leq C \Big ( \sum_{i=1}^M \| u_{i,0} \|_{W^{1,p}(\Omega_i)}^{\max(p,q+1)} + \sum_{j=1}^{|E|} \| w_{j,0} \|_{W^{1,p}(e_j)}^{\max(p,r+1)} \Big ).
\]
For $I_4$ and $I_5$ we use H\"older's inequality and Young's inequality, as well as the uniform boundedness of $u_{i,N}$ in $L^2(0,T;L^2(e_j))$ and of $w_{j,N}$ in $L^\infty(0,T;L^2(e_j))$ to deduce that 
\begin{align*}
|I_4| & \leq  C \sum_{j=1}^{|E|} \sum_{i \, : \, \Omega_i \sim e_j}^{M} \int_0^s \| u_{i,N} \|_{L^2(e_j)}^2 \, dt  + \frac{1}{4} \sum_{j=1}^{|E|} \int_0^s \| \pd_t w_{j,N} \|_{L^2(e_j)}^2 \, dt, \\
|I_5| & \leq C+  \sum_{j=1}^{|E|} \sum_{i \, :\, \Omega_i \sim e_j}^M \frac{\beta_{ij}}{8} \| u_{i,N}(s) \|_{L^2(e_j)}^2
\end{align*}
Meanwhile, for $I_6$ we employ the Cauchy--Schwarz inequality and Young's inequality
\begin{align*}
|I_6| & \leq \frac{1}{4} \sum_{j=1}^{|E|} \int_0^s \| \pd_t w_{j,N} \|_{L^2(e_j)}^2 \, dt + C \sum_{j=1}^{|E|} \sum_{i \, : \, \Omega_i \sim e_j}^{M} \int_0^s \| u_{i,N} \|_{L^2(e_j)}^2 \, dt \\
& \quad + \sum_{j=1}^{|E|} \sum_{i \, : \, \Omega_i \sim e_j}^M \Big (\frac{\beta_{ij}}{8}\| u_{i,N}(s) \|_{L^2(e_j)}^2 + C  \| w_{j,N}(s) \|_{L^2(e_j)}^2 \Big ) \\
& \quad  + C \Big ( \sum_{i =1}^M \| u_{i,0} \|_{W^{1,p}(\Omega_i)}^p +  \sum_{j=1}^{|E|} \| w_{j,0} \|_{L^2(e_j)}^2 \Big ).
\end{align*}
For $I_7$ we have
\begin{align*}
|I_7| & \leq  \sum_{j=1}^{|E|} \sum_{k \, : \, v_k \sim e_j}^{|V|} \Big ( \frac{1}{4} \delta^k_{e_j} |w_{j,N}(s,v_k)|^2 + C\int_0^s |w_{j,N}(v_k)|^2 \, dt \Big ) + \sum_{j=1}^{|E|} \| w_{j,0} \|_{W^{1,p}(e_j)}^2 \\
& \quad + C \sum_{k=1}^{|V|} \Big (|z_{k,N}(s)|^2 + |z_{k,0}|^2 + \int_0^s \Big | \frac{d}{dt} z_{k,N} \Big |^2 \, dt \Big ).
\end{align*}
Recalling the uniform boundedness of $w_{j,N}$ in $L^\infty(0,T;L^2(e_j))$ and of $z_{k,N}$ in $H^1(0,T) \subset L^\infty(0,T)$, as well as well as the assumptions \eqref{ass:aniso}, \eqref{ass:edgeaniso}, \eqref{ass:f} and \eqref{ass:g}, we can infer from \eqref{reg:2} that
\begin{equation}\label{reg:3}
\begin{aligned}
& \sum_{i=1}^M \Big ( c_0 \| \nabla u_{i,N}(s) \|_{L^p(\Omega_i)}^p + \int_0^s \| \pd_t u_{i,N} \|_{L^2(\Omega_i)}^2 \, dt \Big ) + \sum_{j=1}^{|E|} \sum_{i \, : \, \Omega_i \sim e_j}^{M} \frac{\beta_{ij}}{4} \| u_{i,N}(s) \|_{L^2(e_j)}^2 \\
& \qquad + \sum_{j=1}^{|E|} \Big ( c_4 \| \pd_x w_{j,N}(s) \|_{L^p(e_j)}^p + \frac{1}{2} \int_0^s \| \pd_t w_{j,N} \|_{L^2(e_j)}^2 \, dt \Big ) \\
& \qquad + \sum_{j=1}^{|E|} \sum_{k \, : \, v_k \sim e_j}^{|V|} \frac{1}{4} \delta^k_{e_j} |w_{j,N}(s,v_k)|^2  \\
& \quad \leq C + C \sum_{j=1}^{|E|} \sum_{i \, : \, \Omega_i \sim e_j}^{M} \int_0^s \| u_{i,N} \|_{L^2(e_j)}^2 \, dt + C \sum_{j=1}^{|E|} \sum_{k \, : \, v_k \sim e_j}^{|V|} \int_0^s |w_{j,N}(v_k)|^2 \, dt.
\end{aligned}
\end{equation}
Applying Gronwall's inequality in integral form to \eqref{reg:3} leads to 
\begin{align*}
& \sup_{t \in (0,T)} \Big ( \sum_{i=1}^{M} \| \nabla u_{i,N}(s) \|_{L^p(\Omega_i)}^p + \| \pd_x w_{j,N}(s) \|_{L^p(e_j)}^p \Big ) \\
& \quad + \int_0^T \| \pd_t u_{j,N} \|_{L^2(\Omega_i)}^2 + \| \pd_t w_{j,N} \|_{L^2(e_j)}^2 \, dt \leq C,
\end{align*}
and so we infer that along a non-relabelled subsequence $N \to \infty$,
\begin{align*}
u_{i,N} \to u_i & \text{ weakly* in } L^\infty(0,T;W^{1,p}(\Omega_i)), \\
\pd_t u_{i,N} \to \pd_t u_{i} & \text{ weakly in } L^2(\Omega_i \times (0,T)), \\
w_{j,N} \to u_j & \text{ weakly* in } L^\infty(0,T;W^{1,p}(e_j)), \\
\pd_t w_{j,N} \to \pd_t w_j & \text{ weakly in } L^2(e_j \times (0,T)).
\end{align*}
This improved regularity for the time derivatives allows us to utilize the following embedding:
\[
u_i \in H^1(0,T;L^2(\Omega_i)) \subset C^0([0,T];L^2(\Omega_i))
\]
so that $u_i(0)$ is a well-defined function belonging to $L^2(\Omega_i)$. Hence, from the proof of Theorem \ref{thm:well} we obtain the attainment of the initial data $u_{i}(0) = u_{i,0}$ a.e.~in $\Omega_i$ even in the case $p \in (2,\infty)$.  Lastly, we return to \eqref{fp:z} and use that $w_{j,N}(v_k)$ is bounded in $L^\infty(0,T)$ to deduce that $\frac{d}{dt} z_{k,N}$ is bounded in $L^\infty(0,T)$, which yields that the limit $z_k$ belonging to $W^{1,\infty}(0,T)$.

\subsection{Boundedness of solutions}
As remarked in Remark \ref{rem:bdd} we can use the Sobolev embedding $W^{1,p}(e_j) \subset C^0(e_j)$ to deduce that $w_j \in L^\infty(0,T;L^\infty(e_j))$ for all $1 \leq j \leq |E|$. We now employ a standard comparison principle to deduce the boundedness of $u_i$. Note that for the case $p \in (2,\infty)$ we may use the Sobolev embedding $W^{1,p}(\Omega_i) \subset C^0(\Omega_i)$ for $p \in (2,\infty)$ to deduce that $u_i \in L^\infty(0,T;L^\infty(\Omega_i))$. 

We define
\[
\mathcal{M} := \max \Big ( \max_{j=1, \dots, |E|} \max_{i \, : \, \Omega_i \sim e_j} \Big \{ \frac{\alpha_{ij}}{\beta_{ij}} \| w_j \|_{L^\infty(0,T;L^\infty(e_j))} \Big \}, \max_{i = 1, \dots, M} \{ \| u_{i,0} \|_{L^\infty(\Omega_i)} \} \Big )
\]
and consider for each $i = 1, \dots, M$ the function
\[
u_i^{\mathcal{M}} := \max( u_i - \mathcal{M}, 0) = \begin{cases}
    u_i - \mathcal{M}, & \text{ if } u_i > \mathcal{M}, \\
    0, & \text{ if } u_i \leq \mathcal{M}.
\end{cases}
\]
so that 
\begin{align*}
\int_{\Omega_i} \pd_t u_i u_i^{\mathcal{M}} \, dx & = \frac{d}{dt} \frac{1}{2} \| u_i^{\mathcal{M}} \|_{L^2(\Omega_i)}^2, \\
\int_{\Omega_i} f(u_i) u_i^{\mathcal{M}} \, dx & = \int_{\Omega_i} (f(u_i) - f(\mathcal{M})) u_i^{\mathcal{M}} \, dx + \int_{\Omega_i} f(\mathcal{M}) u_i^{\mathcal{M}} \, dx \geq 0.
\end{align*}
Furthermore, via the relation
\[
\nabla u_i^{\mathcal{M}} = \begin{cases}
    \nabla u_i, & \text{ if } u_i > \mathcal{M}, \\
    0, & \text{ if } u_i < \mathcal{M},
\end{cases}
\]
we see that 
\begin{align*}
\int_{\Omega_i} \kappa(\nabla u_i) \cdot \nabla u_i^{\mathcal{M}} \, dx  = \int_{\Omega_i \cap \{ u_i > \mathcal{M}\}} \kappa(\nabla u_i^{\mathcal{M}}) \cdot \nabla u_i^{\mathcal{M}} \, dx \geq c_3 \| \nabla u_i^{\mathcal{M}} \|_{L^p(\{u_i > \mathcal{M}\})}^p \geq 0.
\end{align*}
Hence, upon testing the variational equality for $u_i$ with $u_i^{\mathcal{M}}$ and employing the above relations, we obtain
\begin{align*}
\frac{d}{dt} \frac{1}{2} \| u_i^{\mathcal{M}} \|_{L^2(\Omega_i)}^2 & \leq \sum_{j \, : \, e_j \sim \Omega_i}^{|E|} \int_{e_j} ((\alpha_{ij} w_j - \beta_{ij}\mathcal{M}) - \beta_{ij} (u_i - \mathcal{M})) u_i^{\mathcal{M}} \, dx \\
& \leq \sum_{j \, : \, e_j \sim \Omega_i}^{|E|}  \Big ( \int_{e_j} (\underbrace{\alpha_{ij} \| w_j \|_{L^\infty(e_j)} - \beta_{ij}\mathcal{M})}_{\leq 0}  u_i^{\mathcal{M}} \, dx - \int_{e_j}\beta_{ij} |u_{i}^{\mathcal{M}}|^2 \, dx \Big ) \leq 0.
\end{align*}
Using $u_i^{\mathcal{M}}(0) = \max(u_{i,0} - \mathcal{M}, 0) = 0$ we see that for all $t > 0$,
\[
\| u_i^{\mathcal{M}}(t) \|_{L^2(\Omega_i)}^2 \leq 0,
\]
which implies $u_i(t,x) \leq \mathcal{M}$ for a.e.~$x \in \Omega_i$ and for all $t > 0$. In a similar fashion we consider the function
\[
U_i^{\mathcal{M}} := \max(-u_i - \mathcal{M}, 0) = \begin{cases}
    -u_i - \mathcal{M}, & \text{ if } u_i <  - \mathcal{M}, \\
    0, & \text{ if } u_i \geq - \mathcal{M},
\end{cases}
\]
so that with the following relations
\begin{align*}
& \nabla U_i^{\mathcal{M}} = \begin{cases}
    -\nabla u_i, & \text{ if } u_i <  - \mathcal{M}, \\
    0, & \text{ if } u_i > - \mathcal{M},
\end{cases} \quad \int_{\Omega_i} \pd_t u_i U_i^{\mathcal{M}} \, dx = - \frac{1}{2} \frac{d}{dt} \| U_i^{\mathcal{M}} \|_{L^2(\Omega_i)}^2, \\
& -\int_{\Omega_i} \kappa(\nabla u_i) \cdot \nabla U_i^{\mathcal{M}} \, dx = \int_{\Omega_i \cap \{u_i < - \mathcal{M}\}} \kappa(\nabla u_i) \cdot \nabla u_i \, dx \geq c_3 \| \nabla u_i \|_{L^p(\{u_i < - \mathcal{M}\})}^p \geq 0, \\
& \int_{\Omega_i} f(u_i) U_i^{\mathcal{M}} \, dx = \int_{\Omega_i \cap \{u_i < - \mathcal{M}\}} (f(u_i) - f(-\mathcal{M})) (-u_i - \mathcal{M}) + f(-\mathcal{M}) (-u_i - \mathcal{M}) \, dx \leq 0,
\end{align*}
upon testing the variational equality for $u_i$ with $U_i^{\mathcal{M}}$ we obtain
\begin{align*}
& \frac{1}{2} \frac{d}{dt} \| U_i^{\mathcal{M}} \|_{L^2(\Omega_i)}^2 \leq \sum_{j \, : \, e_j \sim \Omega_i} \int_{e_j} (\beta_{ij} (u_i + \mathcal{M}) - \alpha_{ij} w_{j} - \beta_{ij} \mathcal{M} ) U_i^{\mathcal{M}} \, dx \\
& \quad \leq \sum_{j \, : \, e_j \sim \Omega_i} \Big ( - \beta_{ij} \| U_i^{\mathcal{M}} \|_{L^2(e_j)}^2 + \int_{e_j \cap \{u_i < - \mathcal{M}\}} (\underbrace{\beta_{ij} \mathcal{M} + \alpha_{ij} w_j}_{\geq 0}) (u_i + \mathcal{M}) \, dx \Big ) \leq 0.
\end{align*}
Using that $U_i^{\mathcal{M}}(0) = \max(-u_{i,0} - \mathcal{M},0) = 0$ we see that for all $t > 0$,
\[
\| U_i^{\mathcal{M}}(t) \|_{L^2(\Omega_i)}^2 \leq 0,
\]
which implies $u_i(t,x) \geq - \mathcal{M}$ for a.e.~$x \in \Omega_i$ and for all $t > 0$. Hence, combining with the previous upper bound, we infer the boundedness property:
\[
|u_i(t,x)| \leq \mathcal{M} \quad \text{ for a.e. } x \in \Omega_i, \quad \forall t \geq 0.
\]

\subsection{Spatial regularity}
Suppose \eqref{ass:edgeaniso} holds with $p = 2$ and $\eta$ is globally Lipschitz continuous. Then, on the edge $e_j$, the equation for $w \vert_{e_j} = w_j$ expressed in strong formulation reads as
\begin{align}\label{edge:ell:reg}
- \pd_x(\eta(\pd_x w_j)) = - \pd_t w_j - g(w_j) - \sum_{m \,: \, \Omega_m \sim e_j}^M (\alpha_j w_j - \beta_m u_m) & \text{ in } e_j.
\end{align}
From the regularities $w_j \in L^\infty(0,T;H^1(e_j)) \cap H^1(0,T;L^2(e_j)) \cap L^\infty(0,T;L^\infty(e_j))$, $u_m \in L^\infty(0,T;H^{\frac{1}{2}}(e_j))$, and $g(w_j) \in L^\infty(0,T;L^2(e_j))$, we see that the right-hand side of \eqref{edge:ell:reg} belong to $L^2(0,T;L^2(e_j))$. Invoking Theorem \ref{thm:interior:reg}, noting that the metric graph $E$ is without boundary, we obtain 
\[
\| w \|_{L^2(0,T;H^2(E))} \leq C.
\]
Next, suppose $\kappa(\bm{s}) = \frac{1}{2} |\bm{s}|^2$ holds. Then, on a convex subdomain $\Omega_i$ the equation for $u_i$ expressed in strong formulation reads as
\begin{align*}
\begin{cases}
- \Delta u_i = - \pd_t u_i - f(u_i) & \text{ in } \Omega_i, \\
\nabla u_i \cdot \bm{n} = \alpha_{ij} w - \beta_{ij} u_i & \text{ on } \pd \Omega_i.
\end{cases}
\end{align*}
Invoking classical results for the Laplacian with Robin boundary conditions on convex domains, see e.g.~the proof of Theorem 3.1.2.3 of \cite{Grisvard} with a minor modification, from the regularities $\pd_t u_i + f(u_i) \in L^2(0,T;L^2(\Omega_i))$ and $\alpha_{ij} w \in L^2(0,T;H^{1}(\pd \Omega_i))$ we deduce that 
\[
\| u_i \|_{L^2(0,T;H^2(\Omega_i))} \leq C.
\]

\section{Finite time extinction for unpopulated vertices setting}\label{sec:finitetime}
We adapt the strategy of proof in \cite[Section 5]{Antontsev}. Choosing $\phi^{(i)} = u_i$, $\psi^{(j)} = w_j$ in \eqref{weakform:alt:unpop} yields after applying \eqref{ass:aniso}, \eqref{ass:edgeaniso} and \eqref{fg:finitetime}
\begin{equation}\label{finitetime:est}
\begin{aligned}
& \frac{1}{2} \frac{d}{dt} \Big ( \sum_{i=1}^M \| u_i \|_{L^2(\Omega_i)}^2 + \sum_{j=1}^{|E|} \| w_j \|_{L^2(e_j)}^2 \Big ) \\
& \quad + \sum_{i=1}^{M} \Big ( c_3 \| \nabla u_i \|_{L^p(\Omega_i)}^p + c_{12} \| u_i \|_{L^\sigma(\Omega_i)}^{\sigma} \Big ) + \sum_{j=1}^{|E|} \Big (c_7 \| \pd_x w_{j} \|_{L^p(e_j)}^{p} + c_{13} \| w_j \|_{L^\sigma(e_j)}^{\sigma} \Big ) \leq 0,
\end{aligned}
\end{equation}
where we have neglected the non-negative interaction terms. For spatial dimensions $d \in \{1, 2\}$, we define
\[
\theta_d = \frac{2-\sigma}{2} \frac{dp}{dp + \sigma(p-d)} \in (0,1)
\]
in light of $p > 1$ and $\sigma < 2$. Then, by the Gagliardo--Nirenberg inequality, it holds that
\begin{align*}
\| u_i \|_{L^2(\Omega_i)} \leq C \| \nabla u_i \|_{L^{p}(\Omega_i)}^{\theta_2} \| u_i \|_{L^\sigma(\Omega_i)}^{1 - \theta_2} + C\| u_i \|_{L^\sigma(\Omega_i)}, \\
\| w_j \|_{L^2(e_j)} \leq C \| \pd_x w_j \|_{L^p(e_j)}^{\theta_1} \| w_j \|_{L^\sigma(e_j)}^{1-\theta_1} + C\| w_j \|_{L^\sigma(e_j)},
\end{align*}
and so using that $u_i \in L^\infty(0,T;L^2(\Omega_i))$ and $\sigma < 2$,
\begin{align*}
\| u_i \|_{L^2(\Omega_i)}^2 & \leq C \Big ( \| \nabla u_i \|_{L^p(\Omega_i)} + \| u_i \|_{L^\sigma(\Omega_i)} \Big )^{2 \theta_2} \| u_i \|_{L^\sigma(\Omega_i)}^{2(1-\theta_2)} \\
& \leq C \Big (\| \nabla u_i \|_{L^p(\Omega_i)}^p + \| u_i \|_{L^\sigma(\Omega_i)}^{\sigma} \| u_i \|_{L^\sigma(\Omega_i)}^{p-\sigma} \Big )^{\frac{2 \theta_2}{p}} \Big ( \| u_i \|_{L^\sigma(\Omega_i)}^{\sigma} \Big )^{\frac{2 (1-\theta_2)}{\sigma}} \\
& \leq C \Big (\| \nabla u_i \|_{L^p(\Omega_i)}^p + \| u_i \|_{L^\sigma(\Omega_i)}^{\sigma}  \Big )^{\frac{2 \theta_2}{p}+ \frac{2(1-\theta_2)}{\sigma}} \\
& \leq C \Big ( \sum_{i=1}^M  \Big ( \| \nabla u_i \|_{L^p(\Omega_i)}^p + \| u_i \|_{L^\sigma(\Omega_i)}^{\sigma} \Big )  \Big )^{\frac{2 \theta_2}{p}+ \frac{2(1-\theta_2)}{\sigma}}.
\end{align*}
A short calculation shows that with $\sigma < 2$ and $p \geq 1 > \frac{2\sigma}{2+\sigma}$,
\[
\frac{1}{s_2} := \frac{2 \theta_2}{p} + \frac{2(1-\theta_2)}{\sigma} = 1 + \frac{p(2-\sigma)}{p(2+\sigma) - 2 \sigma}> 1,
\]
and so 
\begin{align*}
\sum_{i=1}^M \| u_i \|_{L^2(\Omega_i)}^2  \leq CM \Big ( \sum_{i=1}^M  \Big ( \| \nabla u_i \|_{L^p(\Omega_i)}^p + \| u_i \|_{L^\sigma(\Omega_i)}^{\sigma} \Big )  \Big )^{\frac{1}{s_2}},
\end{align*}
which yields
\begin{align}\label{finitetime:1}
\Big (\sum_{i=1}^M \| u_i \|_{L^2(\Omega_i)}^2 \Big )^{s_2}  \leq C\sum_{i=1}^M  \Big ( \| \nabla u_i \|_{L^p(\Omega_i)}^p + \| u_i \|_{L^\sigma(\Omega_i)}^{\sigma} \Big ).
\end{align}
Similarly, using $w_j \in L^\infty(0,T;L^2(e_j))$ we infer
\begin{align*}
\| w_j \|_{L^2(e_j)}^2 & \leq C\Big ( \| \pd_x w_j \|_{L^p(e_j)}^p + \| w_j \|_{L^\sigma(e_j)}^\sigma \Big )^{\frac{2 \theta_1}{p} + \frac{2(1-\theta_1)}{\sigma}} \\
& \leq C \Big ( \sum_{j=1}^{|E|} \Big ( \| \pd_x w_j \|_{L^p(e_j)}^p + \| w_j \|_{L^\sigma(e_j)}^\sigma \Big )\Big )^{\frac{2 \theta_1}{p} + \frac{2(1-\theta_1)}{\sigma}}
\end{align*}
where
\[
\frac{1}{s_1} := \frac{2 \theta_1}{p} + \frac{2(1-\theta_1)}{\sigma} =1 + \frac{p(2-\sigma)}{p(1+\sigma) - \sigma} > 1.
\]
Hence, we obtain
\begin{align*}
\Big ( \sum_{j=1}^{|E|} \|w_j \|_{L^2(e_j)}^2 \Big)^{s_1} \leq C \sum_{j=1}^{|E|} \Big ( \| \pd_x w_j \|_{L^p(e_j)}^p + \| w_j \|_{L^\sigma(e_j)}^\sigma \Big ).
\end{align*}
Since $p \geq 2 > \sigma$, it holds that $s_1 < s_2$, and so we infer with $w_j \in L^\infty(0,T;L^2(e_j))$,
\begin{equation}\label{finitetime:2}
\begin{aligned}
\Big ( \sum_{j=1}^{|E|} \|w_j \|_{L^2(e_j)}^2 \Big)^{s_2}&  = \Big ( \sum_{j=1}^{|E|} \|w_j \|_{L^2(e_j)}^2 \Big)^{s_2 - s_1} \Big ( \sum_{j=1}^{|E|} \|w_j \|_{L^2(e_j)}^2 \Big)^{s_1} \\
& \leq C \sum_{j=1}^{|E|} \Big ( \| \pd_x w_j \|_{L^p(e_j)}^p + \| w_j \|_{L^\sigma(e_j)}^\sigma \Big ).
\end{aligned}
\end{equation}
Substituting the lower bounds \eqref{finitetime:1} and \eqref{finitetime:2} into \eqref{finitetime:est}, and using the inequality $(a+b)^{s_2} \leq a^{s_2} + b^{s_2}$ for non-negative constants $a$ and $b$ then yields the ordinary differential inequality:
\begin{align*}
   \frac{d}{dt} \Big ( \sum_{i=1}^M \| u_i \|_{L^2(\Omega_i)}^2 + \sum_{j=1}^{|E|} \| w_j \|_{L^2(e_j)}^2 \Big ) + C \Big ( \sum_{i=1}^M \| u_i \|_{L^2(\Omega_i)}^2 + \sum_{j=1}^{|E|} \| w_j \|_{L^2(e_j)}^2 \Big )^{s_2} \leq 0. 
\end{align*}
Setting $X(t) =  \sum_{i=1}^M \| u_i (t)\|_{L^2(\Omega_i)}^2 + \sum_{j=1}^{|E|} \| w_j (t) \|_{L^2(e_j)}^2$, the above inequality yields
\[
\frac{1}{X^{s_2}} \frac{d}{dt} X \leq -C
\]
and solving this ordinary differential inequality yields
\[
X(t) \leq (X(0)^{1-s_2} - (1-s_2) Ct)^{\frac{1}{1-s_2}}.
\]
This translates to 
\begin{align*}
& \sum_{i=1}^M \| u_i(t) \|_{L^2(\Omega_i)}^2 + \sum_{j=1}^{|E|} \| w_j(t) \|_{L^2(e_j)}^2 \\
 & \quad \leq \begin{cases}
     \Big ( \Big (\sum_{i=1}^M \| u_{i,0} \|_{L^2(\Omega_i)}^2 + \sum_{j=1}^{|E|} \| w_{j,0} \|_{L^2(e_j)}^2 \Big )^{1-s_2} - (1-s_2) Ct \Big )^{\frac{1}{1-s_2}} & \text{ for } t < t_*, \\
     0 & \text{ for }t \geq t_*,
 \end{cases}
\end{align*}
where the extinction time $t_*$ is computed as
\[
t_* = \frac{1}{C(1-s_2)} \Big ( \Big (\sum_{i=1}^M \| u_{i,0} \|_{L^2(\Omega_i)}^2 + \sum_{j=1}^{|E|} \| w_{j,0} \|_{L^2(e_j)}^2 \Big )^{1-s_2}.
\]

\section{Conclusion}
In this work we studied a mathematical model coupling ordinary differential equations at vertices with one-dimensional reaction-diffusion equations along edges and two-dimensional reaction-diffusion equations within subdomains, constituting a system exhibiting a mixed 2D-1D-0D structure. This framework facilitates the investigation of how metric graph structures embedded in two-dimensional domains influence diverse diffusion phenomena in e.g.~epidemic propagation, information dissemination and contaminant transport.

We established the existence, uniqueness and continuous dependence of weak solutions in appropriate function spaces via a semi-Galerkin appoach, along with spatial and temporal regularity properties and finite-time extinction. These enable us to relate the model to physical expectations with regards to reaction-diffusion dynamics, as well as open up future investigations in the context of optimal control.

The mathematical framework we used can be modified to include prescribed source terms on the right-hand sides of \eqref{u:domain}, \eqref{w:edge} and \eqref{pop:equ}, as well as extending to the setting of coupled reaction-diffusion systems in each subdomain $\Omega_i$, edge $e_j$ and coupled ordinary differential systems on vertices $v_k$. Hence, systems of reaction-diffusion equations with network structures, such as the epidemiological settings for modeling diffusive social phenomena \cite{GS,Nizamani} and biological settings involving microvascular networks \cite{Fritz,Koppl} can be approached with a similar strategy employed here.

\appendix
\section{Model derivation via formal asymptotic expansions}\label{sec:derivation}
In this section we derive the model equations \eqref{dom:equ}, \eqref{w:edge} and \eqref{pop:equ} by means of formally match asymptotic expansions, similarly in nature to the techniques used for identifying sharp interface limits of phase field models, see e.g.~\cite{Fife,Garcke} for further details.

\subsection{Equations on edges}
\begin{figure}[h]
    \centering
\includegraphics[width=0.7\textwidth]{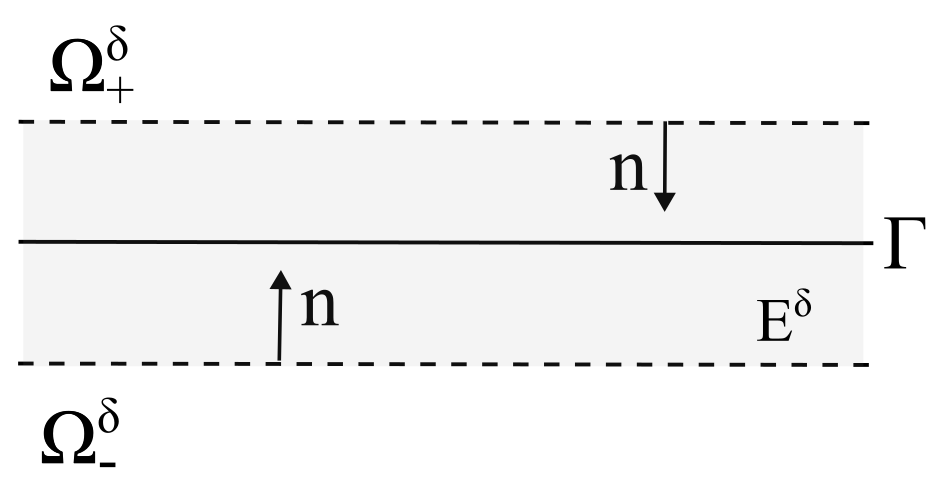}
    \caption{Schematics of the setting for the formal derivation of equations on edges: The region $E^\delta$ formally shrinks to the $x_1$-axis $\Gamma$ and becomes an edge separating the upper and lower half plane. The unit normals $\bm{n}$ on the dashed interfaces representing $\{x_2 = \pm \frac{\delta}{2}\}$ are oriented so that they point into $E^\delta$. }
    \label{fig:2Dedge}
\end{figure}

Let us consider the geometry as depicted in Figure~\ref{fig:2Dedge}, where for $\delta \in (0,1)$ we denote by $\Omega_-^\delta := \{(x_1,x_2) \in \R^2  : x_2 < -\delta/2\}$, $\Omega_+^\delta := \{ (x_1, x_2) \in \R^2 : x_2 > \delta/2\}$ and $E^\delta := \{ (x_1, x_2) \in \R^2 : | x_2 | < \delta/2\}$. Furthermore, we set $\bm{n}$ as the outer unit normal point away from $\Omega_{\pm}^\delta$ and into $E^\delta$. In this domain we pose the system of equations:
\begin{subequations}\label{app:dom:eq1}
    \begin{alignat}{3}
    \pd_t u_\delta  - \Delta u_\delta + f(u_\delta) & = 0  && \quad \text{ in } \Omega_+^\delta = \{x_2 > \delta/2\}, \\
    \pd_t v_\delta - \Delta v_\delta + f(v_\delta)  & = 0 && \quad \text{ in } \Omega_-^\delta = \{ x_2 < -\delta/2\}, \\
    \pd_t w_\delta - \Delta w_\delta + g(w_\delta)  & = 0 && \quad \text{ in } E^\delta = \{ | x_2| < \delta/2\}, 
    \end{alignat}
\end{subequations}
furnished with the boundary conditions
\begin{subequations}\label{app:dom:bc}
\begin{alignat}{2}
    -\frac{1}{\delta}  \pd_{x_2} w_\delta = \mu(v_\delta - w_\delta)=  \pd_{x_2} v_\delta & \quad \text{ on } \pd \Omega_-^\delta = \{ x_2 = -\delta/2\}, \\
    \frac{1}{\delta} \pd_{x_2} w_\delta = \theta(u_\delta - w_\delta) =  \pd_{x_2}u_\delta & \quad \text{ on } \pd \Omega_+^\delta = \{ x_2 = \delta/2\},
\end{alignat}
\end{subequations}
where $\mu$ and $\theta$ are nonnegative constants. We note that only the second component of the fluxes appear since $\bm{n} = (0,-1)^{\top}$ on the top boundary $\{x_2 = \delta/2\}$ and $\bm{n} = (0,1)^{\top}$ on the bottom boundary $\{x_2 = -\delta/2\}$.

Since the domains $\Omega_{\pm}^\delta$ change as $\delta$ varies, we introduce fixed domains $\R_+^2 = \{ x_2 > 0\}$ as the upper-half plane, $\R_-^2 = \{ x_2 < 0\}$ as the lower-half plane, along with $\Gamma := \{x_2 = 0\}$. Then, we make the ansatz that $u_\delta$ and $v_\delta$ have the following asymptotic expansions in $\delta$:
\begin{subequations}\label{outer:expand}
\begin{alignat}{2}
u_\delta(t,\x) & = u_0(t,\x) + \delta u_1(t,\x) + \delta^2 u_2(t,\x) + \cdots, \\
v_\delta(t,\x) & = v_0(t,\x) + \delta v_1(t,\x) + \delta^2 v_2(t,\x) + \cdots,
\end{alignat}
\end{subequations}
with smooth functions $u_i : [0,T] \times \R_+^2 \to \R$ and $v_i :[0,T] \times \R_-^2 \to \R$ for $i = 0, 1, 2, \dots$.  Substituting these expansions into \eqref{app:dom:eq1} and solving them order by order, we obtain to leading order $\mathcal{O}(\delta^0)$ that the pair $(u_0, v_0)$ satisfy
\begin{align*}
\pd_t u_0 - \Delta u_0 + f(u_0) = 0 & \quad \text{ in } \R_+^2, \\
\pd_t v_0 - \Delta v_0 + f(v_0) = 0 & \quad \text{ in } \R_-^2.
\end{align*}
While similar systems for higher order expansions $(u_i, v_i)$ can be derived, they are not necessary in this formal derivation of the model. It remains to supplement the above equations with suitable boundary conditions on $\Gamma$.

For the equation of $w_\delta$ defined in the domain $E^\delta$, which formally becomes $\Gamma$ as $\delta \to 0$, we introduce a change of variable $x := x_1$, $z := \frac{x_2}{\delta}$ and $W_\delta(t,x,z)$ as $w_\delta(t,x_1,x_2)$ under this change of variable, and likewise we use $U_\delta$ and $V_\delta$ to represent $u_\delta$ and $v_\delta$, respectively. Then, via the following formula involving the partial derivatives
\[
\pd_t u = \pd_t U_\delta, \quad \pd_{x_1} u = \pd_{x} U_\delta, \quad \pd_{x_2} u = \frac{1}{\delta} \pd_z U_\delta, \quad \Delta w= \pd_{xx} W_\delta + \frac{1}{\delta^2} \pd_{zz} W_\delta, 
\]
we see that \eqref{app:dom:eq1} transforms to
\begin{subequations}\label{app:sys:scaled}
\begin{alignat}{3}
\label{scaled:u} \pd_t U_\delta & = \pd_{xx} U_\delta  + \frac{1}{\delta^2} \pd_{zz} U_\delta - f(U_\delta), \\
\label{scaled:v} \pd_t V_\delta & = \pd_{xx} V_\delta + \frac{1}{\delta^2} \pd_{zz} V_\delta  -  f(V_\delta), \\
\label{scaled:w} \pd_t W_\delta & = \pd_{xx} W_\delta  + \frac{1}{\delta^2} \pd_{zz} W_\delta - g(W_\delta), 
\end{alignat}
\end{subequations}
while \eqref{app:dom:bc} transforms to 
\begin{subequations}\label{app:sys:scaled:bc}
\begin{alignat}{3}
\label{scaled:bc1} - \frac{1}{\delta^2} \pd_z W_\delta & = \mu (V_\delta - W_\delta)= \frac{1}{\delta} \pd_z V_\delta &&\quad \text{ on } \{ z = -1/2\}, \\
\label{scaled:bc2} \frac{1}{\delta^2} \pd_z W_\delta & = \theta (U_\delta - W_\delta) = \frac{1}{\delta} \pd_z U_\delta && \quad \text{ on } \{ z =1/2\},
\end{alignat}
\end{subequations}
We make the ansatz that $U_\delta$, $V_\delta$ and $W_\delta$ have the following asymptotic expansions in $\delta$:
\begin{subequations}\label{inner:expand}
\begin{alignat}{2}
U_\delta(t,x, z) & = U_0(t,x,z) + \delta U_1(t,x,z) + \delta^2 U_2(t,x,z) + \cdots, \\
V_\delta(t,x, z) & = V_0(t,x,z) + \delta V_1(t,x,z) + \delta^2 V_2(t,x,z) + \cdots, \\
W_\delta(t,x, z) & = W_0(t,x,z) + \delta W_1(t,x,z) + \delta^2 W_2(t,x,z) + \cdots,
\end{alignat}
\end{subequations}
with smooth functions $U_i, V_i, W_i : [0,T] \times \Gamma \times \R \to \R$ for $i = 0, 1, 2, \dots$
and substitute these expansions into \eqref{app:sys:scaled}-\eqref{app:sys:scaled:bc} and solve them order by order. 

To leading order $\mathcal{O}(\delta^{-2})$ in \eqref{scaled:u}, \eqref{scaled:v} and \eqref{scaled:w} we obtain that 
\[
\pd_{zz} U_0 = 0, \quad \pd_{zz} V_0 = 0, \quad \pd_{zz} W_0 = 0,
\]
while from \eqref{scaled:bc1} and \eqref{scaled:bc2} we have to leading order $\mathcal{O}(\delta^{-2})$ for $W_\delta$ and to order $\mathcal{O}(\delta^{-1})$ for $U_\delta$ and $V_\delta$:
\[
\pd_z U_0 \big \vert_{z=1/2} = 0, \quad \pd_z W_0 \big \vert_{z=1/2} = 0, \quad \pd_z V_0 \big \vert_{z=-1/2} = 0, \quad \pd_z W_0 \big \vert_{z=-1/2} = 0.
\]
Then, we see that $\pd_z U_0$, $\pd_z V_0$ and $\pd_z W_0$ must be equal to zero, and hence $U_0$, $V_0$ and $W_0$ are independent of $z$. Below we will provide some matching conditions between $U_0$ and $u_0$ that allow us to relate these two variables. To the next order $\mathcal{O}(\delta^{-1})$ in \eqref{scaled:w} we have
\[
\pd_{zz} W_1 = 0,
\]
and from \eqref{scaled:bc1} and \eqref{scaled:bc2} we have to order $\mathcal{O}(\delta^{-1})$ for $W_\delta$:
\[
\pd_z W_1 \Big \vert_{z=1/2} = 0, \quad \pd_z W_1 \Big \vert_{z=-1/2} = 0.
\]
Hence, we deduce that $\pd_z W_1$ must be equal to zero, and that $W_1$ is independent of $z$.  To second order $\mathcal{O}(\delta^0)$ in \eqref{scaled:w} we have
\begin{align}\label{app:w:equ:1}
\pd_t W_0 - \pd_{xx} W_0 - \pd_{zz} W_2 + g(W_0) = 0,
\end{align}
while from \eqref{scaled:bc1} and \eqref{scaled:bc2} to order $\mathcal{O}(\delta^0)$:
\[
\pd_z W_2 \Big \vert_{z=1/2} = \theta (U_0 - W_0), \quad  -\pd_z W_2 \Big \vert_{z=-1/2} = \mu( V_0 - W_0).
\]
Integrating \eqref{app:w:equ:1} over $z$ from $-1/2$ to $1/2$ and using that $W_0$ is independent of $z$, as well as the above boundary conditions, we obtain
\begin{align}\label{w:0}
\pd_t W_0 - \pd_{xx} W_0 + g(W_0) = \Big [\pd_z W_2 \Big ]_{z=-1/2}^{z=1/2} = \theta U_0 + \mu V_0 - (\theta+ \mu) W_0.
\end{align}
Next, to order $\mathcal{O}(\delta^0)$ in \eqref{scaled:edge:bc1} and \eqref{scaled:edge:bc2} for $U_\delta$ and $V_\delta$, we have
\begin{align}\label{bc:-1}
\mu(V_0 - W_0) = \pd_z V_1 \Big \vert_{z=-1/2}, \quad \theta (U_0 - W_0) = \pd_z U_1 \Big \vert_{z=1/2}.
\end{align}
To relate the two expansions \eqref{outer:expand} and \eqref{inner:expand}, we consider an intermediate region in which both expansions are defined. Let $\alpha \in (0,1)$ and consider a Taylor series expansion of $u_i$ about $\Gamma = \{ x_2 = 0\}$:
\[
u_i(t,x_1,x_2) = u_i(t,x_1,0) + \pd_{x_2} u_i(t,x_1,0) x_2 + \frac{1}{2} \pd_{x_2 x_2} u_i(t,x_1, 0) x_2^2 + \cdots.
\]
Evaluating at $x_2 = \delta^{\alpha}$ provides (suppressing the explicit dependence on $t$ and $x_1$)
\begin{equation}\label{matching1}
\begin{aligned}
u(\delta^{\alpha}) & = \delta^0 u_0(0^+) + \delta^\alpha \pd_{x_2} u_0(0^+) + \delta^{2 \alpha} \tfrac{1}{2} \pd_{x_2 x_2} u_0(0^+) + \cdots \\
& \quad + \delta u_1(0^+) + \delta^{\alpha+1} \pd_{x_2} u_1 (0^+) + \delta^{\alpha+2} \tfrac{1}{2} \pd_{x_2 x_2} u_1(0^+) + \cdots \\
& \quad + \delta^2 u_2(0^+) + \delta^{\alpha+2} \pd_{x_2} u_2(0^+) + \delta^{\alpha+3} \tfrac{1}{2} \pd_{x_2 x_2} u_2(0^+) + \cdots,
\end{aligned}
\end{equation}
where $0^+$ indicates we approach $\Gamma$ from $\R_+^2$. Now, we consider a similar Taylor series expansion of $U_i$ about $\Gamma$:
\[
U_i(t,x,z) = U_{i,0}(t,x,z) + U_{i,1}(t,x,z) z + U_{i,2}(t,x,z) z^2 + \cdots,
\]
and thus, evaluating at $z = \delta^{\alpha-1}$ yields (suppressing the explicit dependence on $t$ and $x$)
\begin{equation}\label{matching2}
\begin{aligned}
U(\delta^{\alpha-1}) & = \delta^0 U_{0,0} + \delta^{\alpha-1} U_{0,1} \delta^{2\alpha-2} U_{0,2} + \cdots \\
& \quad + \delta U_{1,0} + \delta^{\alpha} U_{1,1} + \delta^{2\alpha-1} U_{1,2} + \cdots \\ 
& \quad + \delta^2 U_{2,0} + \delta^{\alpha+1} U_{2,1} + \delta^{2\alpha} U_{2,2} + \cdots.
\end{aligned}
\end{equation}
We say that expansions \eqref{outer:expand} and \eqref{inner:expand} match if the coefficients to every order in $\delta$ agree.  Comparing the two series for $u$ and $U$ leads to the following relations between the coefficients $U_{k,n}$ and the derivatives $\pd_{x_2}^j u_l(0^+)$:
\begin{align*}
U_{0,0} & = u_0(0^+), \quad U_{0,i} = 0 \quad \forall i \geq 1, \\
U_{1,0} & = u_1(0^+), \quad U_{1,1} = \pd_{x_2} u_0(0^+), \quad U_{1,i} = 0 \quad \forall i \geq 2, \\
U_{2,0}& = u_2(0^+), \quad U_{2,1} = \pd_{x_2} u_1(0^+), \quad U_{2,2} = \tfrac{1}{2} \pd_{x_2 x_2} u_0(0^+), \quad U_{2,i} = 0 \quad \forall i \geq 3, 
\end{align*}
where for our current purpose we deduce that 
\begin{align}\label{matching3}
U_0 = u_0 \big \vert_\Gamma, \quad \pd_z U_1 = \pd_{x_2} u_0 \big \vert_{\Gamma}.
\end{align}
A similar argument also yields
\begin{align}\label{matching4}
V_0 = v_0 \big \vert_\Gamma, \quad \pd_z V_1 = \pd_{x_2} v_0 \big \vert_{\Gamma},
\end{align}
and so using the matching conditions \eqref{matching3} and \eqref{matching4} in \eqref{bc:-1} we obtain the boundary conditions for $u_0$ and $v_0$ on $\Gamma$, where with an abuse of notation we use $w_0$ to denote $W_0$, since $W_0$ is independent of $z$:
\begin{align*}
\nabla v_0 \cdot \bm{n} = \pd_{x_2} v_0 = \mu(w_0 - v_0), \quad \nabla u_0 \cdot \bm{n} = \pd_{x_2} u_0 = \theta (w_0 - u_0),
\end{align*}
along with the equation for $w_0$ on $\Gamma$:
\begin{align*}
\pd_t w_0 - \pd_{xx} w_0 + g(w_0) = \theta (u_0 - w_0) + \mu(v_0 - w_0.
\end{align*}
Hence, we formally obtain the following system of equations (dropping the index 0 for presentation)
\begin{equation*}
\begin{alignedat}{3}
& \pd_t u - \Delta u + f(u) = 0 && \quad \text{ in } \R_+^2, \\
&\pd_t v - \Delta v + f(v) = 0 && \quad \text{ in } \R_-^2, \\
&\nabla u \cdot \bm{n} = \theta (w - u) && \quad \text{ on } \Gamma, \\
&\nabla v \cdot \bm{n} = \mu (w - v) && \quad \text{ on } \Gamma, \\
&\pd_t w - \pd_{xx} w + g(w) + \theta(w - u)+ \mu(w - v) = 0 && \quad \text{ on } \Gamma. 
\end{alignedat}
\end{equation*}

\subsection{Equations on vertices}
\begin{figure}[h]
    \centering
\includegraphics[width=0.7\textwidth]{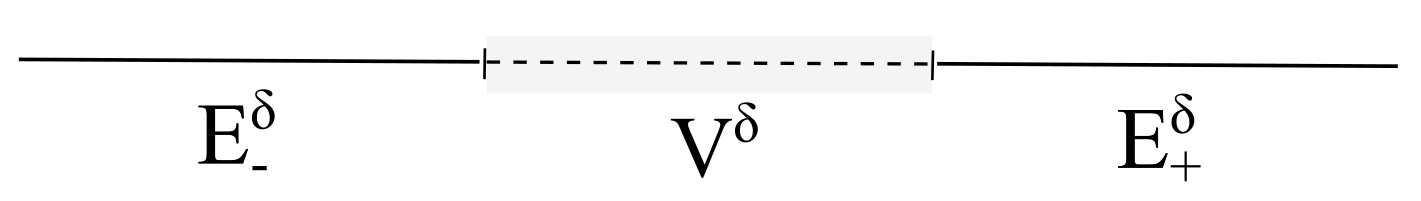}
    \caption{Schematics of the setting for the formal derivation of equations on vertices: The shaded dashed interval $V^\delta$ formally shrinks to the origin and becomes a vertex connecting the left and right edges.}
    \label{fig:1Dvertex}
\end{figure}

Let us consider the geometry as depicted in Figure \ref{fig:1Dvertex}, where for $\delta \in (0,1)$ we denote by $E_-^\delta := (-\infty,-\delta/2)$, $E_+^\delta := (\delta/2, \infty)$ and $V^\delta := (-\delta/2, \delta/2)$. We set $\bm{n}$ as the outer unit normal point away from $e_{\pm}^\delta$ and into $V^\delta$. In this domain we pose the system of equations:
\begin{subequations}\label{app:edge:eq1}
    \begin{alignat}{2}
    \pd_t u_\delta - \pd_{xx} u_\delta + g(u_\delta) = 0 & \quad \text{ in } E_+^\delta,\\
    \pd_t v_\delta - \pd_{xx} v_\delta + g(v_\delta) = 0 & \quad \text{ in } E_-^\delta, \\
    \pd_t w_\delta - \pd_{xx} w_\delta + \lambda w_\delta = 0 & \quad \text{ in } V^\delta, 
    \end{alignat}
\end{subequations}
furnished with the boundary conditions,
\begin{subequations}\label{app:edge:bc}
    \begin{alignat}{2}
   - \frac{1}{\delta} \pd_x w = \mu(v - w) =   \pd_x v \quad & \text{ at } x = -\frac{\delta}{2}, \\
    \frac{1}{\delta} \pd_x w = \theta(u - w) =  \pd_x u \quad & \text{ at } x = \frac{\delta}{2}.
    \end{alignat}
\end{subequations}

\noindent Similar to the previous section, we introduce fixed domains $\R_{>0} := \{ x > 0\}$, $\R_{<0} := \{ x < 0\}$ along with $\Gamma = \{0\}$. Then, we make the ansatz that $u_\delta$ and $v_\delta$ have the following asymptotic expansions in $\delta$:
\begin{subequations}\label{outer:expand:vert}
\begin{alignat}{2}
u_\delta(t,x) & = u_0(t,x) + \delta u_1(t,x) + \delta^2 u_2(t,x) + \cdots, \\
v_\delta(t,x) & = v_0(t,x) + \delta v_1(t,x) + \delta^2 v_2(t,x) + \cdots,
\end{alignat}
\end{subequations}
with smooth functions $u_i : [0,T] \times \R_{>0} \to \R$ and $v_i :[0,T] \times \R_{<0} \to \R$ for $i = 0, 1, 2, \dots$.  Substituting these expansions into \eqref{app:edge:eq1} and solving them order by order, we obtain to leading order $\mathcal{P}(\delta^0)$ that the pair $(u_0, v_0)$ satisfy
\begin{align*}
\pd_t u_0 - \pd_{xx} u_0 + g(u_0) = 0 & \quad \text{ in } \R_{>0}, \\
\pd_t v_0 - \pd_{xx} v_0 + g(v_0) = 0 & \quad \text{ in } \R_{<0}.
\end{align*}
Next, we use a change of variable $z := \frac{x}{\delta}$ and write $W_\delta(t,z)$ as $w_\delta(t,x)$ under this change of variable, and likewise we use $U_\delta$, and $V_\delta$ for $u_\delta$ and $v_\delta$ respectively. Then, \eqref{app:edge:eq1} transforms to
\begin{subequations}\label{app:sys:edge:scaled}
\begin{alignat}{3}
\label{scaled:edge:u} \pd_t U_\delta & = \pd_{zz} U_\delta - g(U_\delta),\\
\label{scaled:edge:v} \pd_t V_\delta & = \pd_{zz} V_\delta  - g(V_\delta), \\
\label{scaled:edge:w} \pd_t W_\delta & = \frac{1}{\delta^2} \pd_{zz} W_\delta - \lambda W_\delta,
\end{alignat}
\end{subequations}
while \eqref{app:edge:bc} transforms to
\begin{subequations}\label{app:sys:edge:scaled:bc}
\begin{alignat}{3}
\label{scaled:edge:bc1}  -\frac{1}{\delta^2} \pd_z W_\delta & = \mu (V_\delta - W_\delta)=  \frac{1}{\delta} \pd_z V_\delta  && \quad \text{ on } \{ z = -1/2\}, \\
\label{scaled:edge:bc2} \frac{1}{\delta^2} \pd_z W_\delta & = \theta (U_\delta - W_\delta) = \frac{1}{\delta} \pd_z U_\delta && \quad \text{ on } \{ z = 1/2 \}.
\end{alignat}
\end{subequations}
We make the ansatz that $U_\delta$, $V_\delta$ and $W_\delta$ have the following asymptotic expansions in $\delta$:
\begin{align*}
U_\delta(t,z) & = U_0(t,z) + \delta U_1(t,z) + \delta^2 U_2(t,z) + \cdots , \\
V_\delta(t,z) & = V_0(t,z) + \delta V_1(t,z) + \delta^2 V_2(t,z) + \cdots, \\
W_\delta(t,z) & = W_0(t,z) + \delta W_1(t,z) + \delta^2 W_2(t,z) + \cdots,
\end{align*}
with smooth functions $U_i, V_i, W_i : [0,T] \times \R \to \R$ for $i = 0, 1, 2, \dots$ and substitute these expansions into \eqref{app:edge:eq1}-\eqref{app:edge:bc} and solve them order by order. As the argument is largely similar to the previous section, let us only sketch the main points.

To leading order $\mathcal{O}(\delta^{-2})$ we infer from \eqref{app:sys:edge:scaled} that
\[
\pd_{zz} U_0 = 0, \quad \pd_{zz} V_0 = 0, \quad \pd_{zz}W_0 = 0,
\]
while to leading order in \eqref{app:sys:edge:scaled:bc} we have
\[
\pd_z U_0 \big \vert_{z=1/2} = 0, \quad \pd_z W_0 \big \vert_{z=1/2} =0 , \quad \pd_z V_0 \big \vert_{z=-1/2} = 0, \quad \pd_z W_0 \big \vert_{z=-1/2} = 0.
\]
From these we see that $\pd_{z} U_0$, $\pd_z V_0$ and $\pd_z W_0$ must be zero, and hence $U_0$, $V_0$ and $W_0$ are independent of $z$. To order $\mathcal{O}(\delta^{-1})$ in \eqref{scaled:edge:w} we have
\[
\pd_{zz} W_1 = 0,
\]
and together with the order $\mathcal{O}(\delta^{-1})$ terms in \eqref{app:sys:edge:scaled:bc} for $W_\delta$:
\[
\pd_z W_1 \big \vert_{z=1/2} = 0, \quad \pd_z W_1 \big \vert_{z=-1/2} = 0,
\]
we conclude that $W_1$ is also independent of $z$. Next, to order $\mathcal{O}(\delta^0)$ in \eqref{app:sys:edge:scaled:bc} we have
\[
\pd_z W_2 \big \vert_{z=1/2} = \theta (U_0 - W_0), \quad - \pd_z W_2 \big \vert_{z=-1/2} = \mu(V_0 - W_0).
\]
while to order $\mathcal{O}(\delta^0)$ in \eqref{scaled:edge:w} we have
\[
\frac{d}{dt} W_0 - \pd_{zz} W_2 + \lambda W_0 = 0.
\]
Integrating the above equality over $z$ from $-1/2$ to $1/2$, and using the boundary conditions as well as the fact that $W_0$ is independent of $z$, we arrive at
\[
\frac{d}{dt} W_0 + \lambda W_0 = \Big [ \pd_z W_2 \Big ]_{z=-1/2}^{z=1/2} = \theta U_0 + \mu V_0 - (\theta + \mu) W_0.
\]
Meanwhile, to order $\mathcal{O}(\delta^0)$ in \eqref{app:sys:edge:scaled:bc} for $U_\delta$ and $V_\delta$, we have
\[
\mu (V_0 - W_0) = \pd_z V_1 \big \vert_{z=-1/2}, \quad \theta (U_0 - W_0) = \pd_z U_1 \big \vert_{z=1/2}.
\]
Then, by the matching conditions, we formally obtain the following systems of equations (dropping the index $0$ for presentation and using $w$ to denote $W_0$)
\begin{equation*}
\begin{alignedat}{3}
& \pd_t u - \pd_{xx} u + g(u) = 0 && \quad \text{ in } \R_{>0}, \\
&\pd_t v - \pd_{xx} v + g(v) = 0 && \quad \text{ in } \R_{<0}, \\
& \pd_x u = \theta (w - u) && \quad \text{ at } \{x = 0\}, \\
& \pd_x v = \mu (w - v) && \quad \text{ at } \{x=0\}, \\
& \frac{d}{dt}w + \lambda w + \theta(w - u)+ \mu(w - v) = 0 && \quad \text{ at } \{x=0\}. 
\end{alignedat}
\end{equation*}

\section{Regularity theory for quasilinear elliptic PDEs}\label{sec:Garcke}
In this section we consider the following quasilinear elliptic PDE:
\begin{align}\label{quasi}
\begin{cases}
- \div (a(\nabla u)) = f & \text{ in } U, \\
a(\nabla u) \cdot \bm{n} = 0 & \text{ on } \pd U,
\end{cases}
\end{align}
where $U \subset \R^n$, $n \in \mathbb{N}$, is a bounded open subset of $\R^n$ with Lipschitz boundary, $f \in L^2(U)$, and the function $a: \R^n \to \R^n$ is Lipschitz continuous and strongly monotone, i.e. there exist positive constants $C_L$ and $C_M$ such that 
\[
|a(\bm{p}) - a(\bm{q})| \leq C_L |\bm{p} - \bm{q}|, \quad (a(\bm{p}) - a(\bm{q})) \cdot (\bm{p} - \bm{q}) \geq C_M |\bm{p} - \bm{q}|^2
\]
for all $\bm{p}, \bm{q} \in \R^n$. We say that $u \in H^1(U)$ is a weak solution to \eqref{quasi} if
\[
\int_U a(\nabla u) \cdot \nabla v \, dx = \int_U f v \, dx,  \quad \forall v \in H^1(U).
\]
Let us state the regularity theory for \eqref{quasi} encompassed in Theorems 5.2 and 5.4 of \cite{GKW} (which are stated for the more general case of a spatially-dependent anisotropy function $a(x, \nabla u)$). 

\begin{thm}[Interior regularity, Theorem 5.2 of \cite{GKW}]\label{thm:interior:reg}
Let $u \in H^1(U)$ be a weak solution to \eqref{quasi}. Then, it holds that 
\[
u \in H^2_{\mathrm{loc}}(U),
\]
and for every open subset $V \Subset U$, there exists a constant $C> 0$ depending only on $U$ and $V$ such that 
\[
\| u \|_{H^2(V)} \leq C \Big ( \| f \|_{L^2(U)} + \| u \|_{L^2(U)} \Big )
\]
holds.
\end{thm}

\begin{thm}[Boundary regularity, Theorem 5.4 of \cite{GKW}]
Suppose $U$ is of class $C^{1,1}$. Let $u \in H^1(U)$ be a weak solution to \eqref{quasi}. Then, it holds that 
\[
u \in H^2(U),
\]
and there exists a constant $C> 0$, which depends only on $U$ and $a$, such that 
\[
\| u \|_{H^2(U)} \leq C \Big ( \| f \|_{L^2(U)} + \| u \|_{L^2(U)} \Big )
\]
holds.
\end{thm}

\end{document}